\documentclass[10pt]{article}
\usepackage{authblk}
\usepackage[margin=1in]{geometry} 
\usepackage{amsmath,amsthm,amssymb}
\usepackage{enumitem}
\usepackage{soul}
\usepackage{mathtools}
\usepackage[utf8]{inputenc}
\usepackage{multirow}
\usepackage[colorlinks]{hyperref}
\usepackage{cleveref}
\usepackage{bbm}
\usepackage{tikz-cd}
\usepackage{adjustbox}
\usepackage{algorithm, algpseudocode}
\allowdisplaybreaks

\makeatletter
\DeclareFontFamily{OMX}{MnSymbolE}{}
\DeclareSymbolFont{MnLargeSymbols}{OMX}{MnSymbolE}{m}{n}
\SetSymbolFont{MnLargeSymbols}{bold}{OMX}{MnSymbolE}{b}{n}
\DeclareFontShape{OMX}{MnSymbolE}{m}{n}{
    <-6>  MnSymbolE5
    <6-7>  MnSymbolE6
    <7-8>  MnSymbolE7
    <8-9>  MnSymbolE8
    <9-10> MnSymbolE9
    <10-12> MnSymbolE10
    <12->   MnSymbolE12
}{}
\DeclareFontShape{OMX}{MnSymbolE}{b}{n}{
    <-6>  MnSymbolE-Bold5
    <6-7>  MnSymbolE-Bold6
    <7-8>  MnSymbolE-Bold7
    <8-9>  MnSymbolE-Bold8
    <9-10> MnSymbolE-Bold9
    <10-12> MnSymbolE-Bold10
    <12->   MnSymbolE-Bold12
}{}

\numberwithin{equation}{section}

\newcommand{\bftheta}{\boldsymbol{\theta}}
\newcommand{\bfsigma}{\boldsymbol{\sigma}}

\newcommand{\bfA}{\mathbf{A}}

\newcommand{\bfg}{\mathbf{g}}

\newcommand{\bfH}{\mathbf{H}}

\newcommand{\bfI}{\mathbf{I}}

\newcommand{\bfo}{\mathbf{o}}
\newcommand{\bfO}{\mathbf{O}}

\newcommand{\bft}{\mathbf{t}}

\newcommand{\bfU}{\mathbf{U}}

\newcommand{\bfV}{\mathbf{V}}

\newcommand{\bfW}{\mathbf{W}}
\newcommand{\bfx}{\mathbf{x}}
\newcommand{\bfX}{\mathbf{X}}

\newcommand{\bfz}{\mathbf{z}}

\newcommand{\bfzero}{\mathbf{0}}

\newcommand{\bbB}{\mathbb{B}}

\newcommand{\bbE}{\mathbb{E}}

\newcommand{\bbG}{\mathbb{G}}

\newcommand{\bbN}{\mathbb{N}}
\newcommand{\bbO}{\mathbb{O}}
\newcommand{\bbP}{\mathbb{P}}
\newcommand{\bbQ}{\mathbb{Q}}
\newcommand{\bbR}{\mathbb{R}}

\newcommand{\bbX}{\mathbb{X}}

\newcommand{\calA}{\mathcal{A}}
\newcommand{\calB}{\mathcal{B}}
\newcommand{\calC}{\mathcal{C}}
\newcommand{\calD}{\mathcal{D}}
\newcommand{\calE}{\mathcal{E}}
\newcommand{\calF}{\mathcal{F}}
\newcommand{\calG}{\mathcal{G}}
\newcommand{\calH}{\mathcal{H}}
\newcommand{\calI}{\mathcal{I}}

\newcommand{\calN}{\mathcal{N}}
\newcommand{\calO}{\mathcal{O}}
\newcommand{\calP}{\mathcal{P}}

\newcommand{\calR}{\mathcal{R}}

\newcommand{\calT}{\mathcal{T}}

\newcommand{\calV}{\mathcal{V}}

\newtheorem{defi}{Definition}[section]
\newtheorem{thm}{Theorem}[section]

\newtheorem{prop}{Proposition}[section]
\newtheorem{lemma}{Lemma}[section]
\newtheorem{rmk}{Remark}[section]

\newtheorem{assumption}{Assumption}

\newcommand{\RNum}[1]{\uppercase\expandafter{\romannumeral #1\relax}}

\DeclareMathOperator*{\argmin}{argmin}

\let\llangle\@undefined
\let\rrangle\@undefined
\DeclareMathDelimiter{\llangle}{\mathopen}%
                     {MnLargeSymbols}{'164}{MnLargeSymbols}{'164}
\DeclareMathDelimiter{\rrangle}{\mathclose}%
                     {MnLargeSymbols}{'171}{MnLargeSymbols}{'171}
\makeatother

\begin{document}

\title{Learning bounds for doubly-robust covariate shift adaptation}

\author{Jeonghwan Lee\footnote{E-mail: \href{mailto:jhlee97@uchicago.edu}{\texttt{jhlee97@uchicago.edu}}.} \ and Cong Ma\footnote{E-mail: \href{mailto:congm@uchicago.edu}{\texttt{congm@uchicago.edu}}.}}
\affil{Department of Statistics at the University of Chicago}

\maketitle

\begin{abstract}
Distribution shift between the training domain and the test domain poses a key challenge for modern machine learning. An extensively studied instance is the \emph{covariate shift}, where the marginal distribution of covariates differs across domains, while the conditional distribution of outcome remains the same. The doubly-robust (DR) estimator, recently introduced by \cite{kato2023double}, combines the density ratio estimation with a pilot regression model and demonstrates asymptotic normality and $\sqrt{n}$-consistency, even when the pilot estimates converge slowly. However, the prior arts has focused exclusively on deriving asymptotic results and has left open the question of non-asymptotic guarantees for the DR estimator.

\indent This paper establishes the first non-asymptotic learning bounds for the DR covariate shift adaptation. Our main contributions are two-fold: (\romannumeral 1) We establish \emph{structure-agnostic} high-probability upper bounds on the excess target risk of the DR estimator that depend only on the $L^2$-errors of the pilot estimates and the Rademacher complexity of the model class, without assuming specific procedures to obtain the pilot estimate, and (\romannumeral 2) under \emph{well-specified parameterized models}, we analyze the DR covariate shift adaptation based on modern techniques for non-asymptotic analysis of MLE, whose key terms governed by the Fisher information mismatch term between the source and target distributions. Together, these findings bridge asymptotic efficiency properties and a finite-sample out-of-distribution generalization bounds, providing a comprehensive theoretical underpinnings for the DR covariate shift adaptation.
\end{abstract}

\section{Introduction}
\label{sec:introduction}

\indent Classical supervised learning assumes that the training and test data are drawn from the same distribution \cite{vapnik2013nature, gyorfi2002distribution}. In practice, such an assumption is rarely met. For instance, credit models are typically trained on approved customers but deployed on rejected applicants; medical imaging data vary across hospitals due to differences in equipment and protocols \cite{koh2021wilds, guan2021domain}; and in natural language processing, labeled corpora such as the \emph{Wall Street Journal}, differ sharply from the domains such as \emph{arXiv} \cite{jiang2007instance}. For all these cases, distribution shift between training and test domains undermines predictive performance.

\indent A significant particular case of such a distribution shift is known as the \emph{covariate shift} \cite{shimodaira2000improving, quinonero-candela2008dataset, pan2009survey}, where the marginal distribution of covariates $X$ varies across the domains while the conditional distribution of $Y|X$ remains the same. Covariate shift is well-documented in healthcare \cite{wei2015health, hajiramezanali2018bayesian}, image classification \cite{saenko2010adapting}, remote sensing \cite{tuia2011using}, sentiment analysis \cite{blitzer2007biographies}, and speech and language processing \cite{yamada2009covariate, hassan2013acoustic, fei2015social}. 

\indent The problem of covariate shift adaptation assumes access to labeled samples from a source domain and unlabeled covariates from a target domain, with the goal of learning a predictor with a desirable performance under the target distribution. This problem has been central to the literature of transfer learning and domain adaptation \cite{sugiyama2007covariate, sugiyama2007direct, sugiyama2012machine, pan2009survey, kato2023double}, especially when the target labels are scarce or costly to obtain.

\indent A core difficulty lies in estimating the covariate density ratio between the source and target domains. The standard approach -- plugging-in an estimated covariate density ratio into an importance-weighted empirical risk minimization \cite{sugiyama2007covariate, sugiyama2007direct, sugiyama2008direct, reddi2015doubly} -- turns out to be highly sensitive to the estimation errors of the density ratio and performs poorly unless the estimator converges at a nearly parametric rate. To address this, \cite{kato2023double} suggests a \emph{doubly-robust (DR) estimator}, which augments the importance-weighting with a pilot regression model and leverages double machine learning techniques \cite{chernozhukov2017double, chernozhukov2018double, chernozhukov2022locally, chernozhukov2023simple, foster2023orthogonal}. Their results establish the asymptotic normality and $\sqrt{n}$-consistency of their DR estimator under parametric models, even when the pilot estimates converge slowly.

\indent Yet, the literature of covariate shift adaptation has centered exclusively on achieving asymptotic results. It remains unclear how the DR covariate shift adaptation performs in finite-sample regimes. This paper aims to close this gap. Our contributions can be summarized as follows:
\begin{enumerate} [label = (\roman*)]
    \item \textbf{Structure-agnostic guarantees}: We first derive the first non-asymptotic upper bounds on the excess target risk for the DR estimator, depending only on the product of the statistical rates of convergence of the pilot estimates, without assumptions on how they are obtained.
    \item \textbf{Fast rates for parameterized models:} By studying the DR estimator through the lens of modern non-asymptotic theory of maximum likelihood estimation (MLE), we prove that the estimator achieves a rate of convergence of the order $\calO \left( 1 / n \right)$ under covariate shift.
\end{enumerate}

\indent Together, these results bridge asymptotic efficiency results and a finite-sample out-of-distribution (OOD) generalization bound, providing a comprehensive theoretical underpinning of the DR covariate shift adaptation.

\subsection{Related works}
\label{subsec:related_works}

\indent We take a moment to discuss subsets of related prior works in covariate shift, doubly-robust estimation, and structure-agnostic estimation framework.

\paragraph{Covariate shift}
The study of covariate shift can be dated back to the seminal paper by \cite{shimodaira2000improving}. This paper investigates the impact of covariate shift under parametric models with the vanilla MLE and proposes the importance-weighting (IW) method, which has a remarkable improvement if the underlying regression model is mis-specified. It also establishes the asymptotic normality for a weighted version of MLE under covariate shift, but no finite-sample learning bounds are provided. Later, \cite{sugiyama2005model} further extends this work by studying an unbiased estimator under the $L^2$-generalization error. Motivated by these fundamental works, there has been a flurry of follow-up works for parametric covariate shift. \cite{mousavi2020minimax} introduces a statistical minimax framework and gives lower bounds for out-of-distribution generalization under the regression models of linear and one-hidden layer neural networks. \cite{lei2021near} takes a closer inspection on the minimax optimal estimator for fixed-design linear regression under covariate shift. \cite{zhang2022class} studies linear models under covariate shift where the learner has access to a small amount of target labels. In stark contrast, this work focuses on the covariate shift problem where the learner has no access to target labels.

\indent Beyond the cases of parametric covariate shift, \cite{cortes2010learning} investigate the IW estimator under the framework of statistical learning and provide a non-asymptotic upper bound on the excess target risk for the IW estimator. Also, there has been a strand of recent works on well-specified non-parametric models under covariate shift. \cite{kpotufe2021marginal} investigates the non-parametric classification problem over the class of H\"{o}lder continuous functions and provides a new fine-grained similarity measure. Within a focus on the class of H\"{o}lder continuous functions, \cite{pathak2022new} introduces a novel measure of distribution mismatch between the source and target domains. Under the setting of reproducing kernel Hilbert space (RKHS), \cite{ma2023optimally, gogolashvili2023importance} establish the optimal learning rates of kernel ridge regression (KRR) estimators. In particular, \cite{ma2023optimally} proves that KRR estimation using a carefully selected regularization parameter is miniax optimal provided that the covariate density ratio is uniformly bounded, and a re-weighting version of the KRR estimator using truncated covariate density ratios is minimax-optimal if the covariate density ratio has a finite second-order moment. On the other hand, \cite{wang2023pseudo} suggests the strategy of learning a predictive model using \emph{pseudo-labels}. As our final remark, over-parameterized models, such as high-dimensional models and classes of neural networks, under covariate shift has drawn increasing attention from the researchers \cite{byrd2019effect, hendrycks2019benchmarking, hendrycks2021many, tripuraneni2021overparameterization}).

\paragraph{Doubly-robust (DR) estimation}
\emph{Doubly-robust (DR) estimation} combines an outcome regression with a model for treatment or selection (e.g., the propensity score), guaranteeing its consistency if at least one is correctly specified. Its foundations lie in the seminal paper by \cite{robins1994estimation} on semi-parametric theory and influence functions, and were formalized for applications by \cite{bang2005doubly}. Some implementations include the \emph{augmented inverse propensity weighting (AIPW)} \cite{robins1994estimation, robins1995semiparametric, bang2005doubly} and \emph{target maximum likelihood estimation (TMLE)} \cite{van2006targeted, van2011targeted}, both of which leverage influence functions to correct bias. A corpus of recent studies integrate modern ML techniques for flexible nuisance estimation together with the Neyman orthogonalization and sample splitting \cite{chernozhukov2017double, chernozhukov2018double, van2018targeted, kennedy2024semiparametric} for retaining valid inference. The DR estimation framework has expanded to settings such as difference-in-differences \cite{sant2020doubly, ning2020doubly}, instrumental variables \cite{okui2012doubly,lee2023doubly}, and censored data \cite{bai2013doubly}. While the DR methods achieve robustness and potential efficiency, they require careful handling of finite-sample bias \cite{kang2007demystifying, funk2011doubly}, near-positivity violations \cite{cole2008constructing}, and model diagnostics \cite{bang2005doubly, robins1994estimation}, since the correctness of at least one nuisance estimate remains crucial.

\paragraph{Structure-agnostic estimation}
The \emph{structure-agnostic estimation} framework stands for a class of statistical methods for estimating functionals or treatment effects without assuming any parametric or structural models for the underlying data generating process. \cite{balakrishnan2023fundamental} establishes the fundamental limits for such functional estimation, characterizing the optimal rates achievable when only minimal assumptions -- such as smoothness or boundedness -- are imposed. \cite{jin2024structure} demonstrates that the DR estimators both for the average treatment effect (ATE) and the average treatment effect on the treated (ATT) attain the minimax optimal rates under the structure-agnostic estimation framework. Their findings underscore the effectiveness of the DR learning in causal inference, particularly when relying on flexible ML algorithms for nuisance estimation. \cite{jin2025s} further studies the sensitivity of structure-agnostic estimation procedures to noise, highlighting several cases where standard estimators fail to achieve normality or efficiency. Finally, \cite{bonvini2024doubly} extends the framework by formalizing the DR inference under smoothness conditions. Collectively, these recent works aim to construct a rigorous framework for statistical estimation and inference that minimizes reliance on structural assumptions while achieving near-optimal statistical guarantees.

\section{Problem formulation}
\label{sec:problem_formulation}

Let $\bbX$ denote the covariate space (a.k.a., the feature space). Consider the \emph{source distribution} $\bbP \in \Delta \left( \bbX \times \bbR \right)$ and the \emph{target distribution} $\bbQ \in \Delta \left( \bbX \times \bbR \right)$. Also, let $\bbP_{X} \in \Delta(\bbX)$ and $\bbQ_{X} \in \Delta(\bbX)$ denote by the marginal distributions of $X$ under $\bbP$ and $\bbQ$, respectively.  We further define $\bbP_{Y \mid X} : \bbX \to \Delta(\bbR)$ and $\bbQ_{Y \mid X} : \bbX \to \Delta(\bbR)$ to be the conditional laws of $Y$ given $X$ under $\bbP$ and $\bbQ$:
\[
    \bbP_{Y \mid X} \left( \cdot \mid x \right) := \bbP \left( Y \in \cdot \left| X = x \right. \right) \quad \textnormal{and} \quad  
    \bbQ_{Y \mid X} \left( \cdot \mid x \right) := \bbQ \left( Y \in \cdot \left| X = x \right. \right).
\]

\begin{assumption} [Covariate shift model]
\label{assumption:covariate_shift}
\normalfont{
For every $x \in \bbX$,
\begin{equation}
    \label{eqn:covariate_shift_model}
    \begin{split}
        \bbE_{\bbP} \left[ Y \mid  X = x \right] = \bbE_{\bbQ} \left[ Y \mid  X = x \right].
    \end{split}
\end{equation}
}
\end{assumption}

\noindent Thus, the two distributions share the same \emph{Bayes regression function} $f^* : \bbX \to \bbR$,
\[
    f^*(x) := \bbE_{\bbP} \left[ Y \mid  X = x \right] = \bbE_{\bbQ} \left[ Y \mid  X = x \right], \quad x \in \bbX.
\]
Here, we emphasize that Assumption \ref{assumption:covariate_shift} does not require $\bbP_{Y \mid X} = \bbQ_{Y \mid X}$; only their Bayes regression functions must coincide. In fact, this assumption is weaker compared to the classical covariate shift model \cite{shimodaira2000improving}, which posits a full equality of the conditional distributions.

\medskip

\paragraph{Observational data.}
We observe $n_{\bbP}$ labeled samples from the source distribution $\bbP$,
\[
    \bfO^{\bbP}_{1:n_{\bbP}} := \left( O_{i}^{\bbP} := \big( X_{i}^{\bbP}, Y_{i}^{\bbP} \big) : i \in \left[ n_{\bbP} \right] \right) \sim \bbP^{\otimes n_{\bbP}},
\]
and $n_{\bbQ}$ \emph{unlabeled} target covariates,
\[
    \bfX^{\bbQ}_{1:n_{\bbQ}} = \left( X_{j}^{\bbQ} : j \in \left[ n_{\bbQ} \right] \right) \sim \bbQ_{X}^{\otimes n_{\bbQ}}.
\]
Hence, the labels are available only in the source domain.

\medskip

\paragraph{Risk and excess risk.} 
Given a function class $\calF \subseteq (\bbX \to \bbR)$, we define the \emph{$\mu$-risk} $\calR_{\mu} : \calF \to \bbR_{+}$ by
\[
    \calR_{\mu} (f) := \bbE_{\left( X, Y \right) \sim \mu} \left[ \left\{ Y - f(X) \right\}^2 \right], \quad \mu \in \Delta \left( \bbX \times \bbR \right).
\]
Given any $\mu \in \Delta \left( \bbX \times \bbR \right)$, let $f^*_\mu \in \argmin \left\{ \calR_{\mu} (f) : f \in \calF \right\}$ denote a $\mu$-risk minimizer over the function class $\calF$. The \emph{excess $\mu$-risk} is then defined by
\begin{equation}
    \label{eqn:excess_risk}
    \begin{split}
        \calE_{\mu} (f) := \calR_{\mu} (f) - \calR_{\mu} \left( f_{\mu}^* \right), \quad f \in \calF.
    \end{split}
\end{equation}

\paragraph{Goal: covariate shift adaptation.} 
Our objective is to construct an estimator $\hat f \in \calF$ that achieves small excess $\bbQ$-risk $\calE_\bbQ(\hat f)$ with high probability.

\paragraph{Covariate density ratio.}
A central quantity in the study of covariate shift is the \emph{covariate density ratio} between the source and target distributions. We assume that the marginal distributions $\bbP_X$ and $\bbQ_X$ are all absolutely continuous with respect to a $\sigma$-finite reference measure $\mu_{\bbX}$ on $\bbX$.  
Let
\[
    p_X := \frac{\mathrm{d} \bbP_X}{\mathrm{d} \mu_{\bbX}} : \bbX \to \bbR_{+} \quad \textnormal{and} \quad q_X := \frac{\mathrm{d} \bbQ_X}{\mathrm{d} \mu_{\bbX}} : \bbX \to \bbR_{+}
\]
denote their respective densities with respect to $\mu_{\bbX}$. The covariate density ratio is then defined as
\[
    \rho^* (x) := \frac{q_X (x)}{p_X( x)}, \quad x \in \bbX,
\]
which is assumed to be finite everywhere throughout this paper.

\section{Doubly-robust (DR) covariate shift adaptation}
\label{sec:dr_covariate_shift_adaptation}

\indent Re-weighting with respect to the source distribution $\bbP$ yields an alternative expression of the $\bbQ$-risk as the $\rho^*$-weighted $\bbP$-risk:
\begin{equation}
    \label{eqn:importance_weighting}
    \begin{split}
        \calR_{\bbQ}(f) = \bbE_{\left( X, Y \right) \sim \bbP} \left[ \rho^* (X) \left\{ Y - f(X) \right\}^2 \right], \quad f \in \calF.
    \end{split}
\end{equation}

\noindent The \emph{importance-weighting (IW)} estimator  \cite{shimodaira2000improving} can be obtained by minimizing the empirical analogue of the $\rho^*$-weighted $\bbP$-risk \eqref{eqn:importance_weighting} over $\calF$. Its key limitation is the reliance on the knowledge of the \emph{unknown} covariate density ratio $\rho^* : \bbX \to \bbR_{+}$: a modified estimator obtained by plugging-in an estimate $\hat{\rho} : \bbX \to \bbR_{+}$ for the covariate density ratio $\rho^*$ might have high variance and degrade its performance unless the estimation of the covariate density ratio is sufficiently accurate.

\indent The \emph{doubly-robust (DR) covariate shift adaptation} \cite{kato2023double} augments the IW method with a pilot regression model, and then subtracts a squared-error correction term to cancel the leading error term incurred by the density ratio estimation. For any given pilot estimates $\hat{\rho} : \bbX \to \bbR_{+}$ and $\hat{f}_0 : \bbX \to \bbR$ for the covariate density ratio $\rho^* : \bbX \to \bbR_{+}$ and the shared Bayes regression function $f^* : \bbX \to \bbR$, respectively, let us define the DR empirical risk $\widehat{\calR}_{\textsf{DR}} : \calF \to \bbR$ by
\begin{equation}
    \label{eqn:dr_empirical_risk_v1}
    \begin{split}
        \widehat{\calR}_{\textsf{DR}} (f) := \ & \frac{1}{n_{\bbP}} \sum_{i=1}^{n_{\bbP}} \hat{\rho} \big( X_{i}^{\bbP} \big) \left[ \left\{ Y_{i}^{\bbP} - f \big( X_{i}^{\bbP} \big) \right\}^2 - \left\{ \hat{f}_0 \big( X_{i}^{\bbP} \big) - f \big( X_{i}^{\bbP} \big) \right\}^2 \right] \\
        &+ \frac{1}{n_{\bbQ}} \sum_{j=1}^{n_{\bbQ}} \left\{ \hat{f}_0 \big( X_{j}^{\bbQ} \big) - f \big( X_{j}^{\bbQ} \big) \right\}^2
    \end{split}
\end{equation}
and the DR estimator as
\begin{equation}
    \label{eqn:dr_estimator_v1}
    \begin{split}
        \hat{f}_{\textsf{DR}} \in \argmin \left\{ \widehat{\calR}_{\textsf{DR}} (f) : f \in \calF \right\}.
    \end{split}
\end{equation}
Intuitively, the pilot regression model terms $\hat f_0$ makes the risk \emph{orthogonal} to the first-order errors in $\hat{\rho}$ (and vice-versa), yielding stability even when the pilot estimates converge slowly.

\paragraph{Structure-agnostic estimation.}
Throughout this section, the pilot estimates $\hat{\rho} : \bbX \to \bbR_{+}$ and $\hat{f}_0 : \bbX \to \bbR$ are regarded as \emph{black-boxes}: the analysis only requires the pilot estimates to achieve certain statistical error rates, not how these estimates are obtained. This \emph{structure-agnostic estimation} framework \cite{balakrishnan2023fundamental, jin2024structure, kennedy2024minimax, bonvini2024doubly, jin2025s} reflects practice, where the pilot estimates $\hat{\rho}$ and $\hat{f}_0$ can be obtained by leveraging a growing body of modern ML methods (e.g., LASSO \cite{bicker2009simultaneous, wainwright2009sharp}, tree-based algorithms \cite{syrgkanis2020estimation, wager2018estimation}, and deep neural networks \cite{chen1999improved, hieber2020nonparametric}). Later, our finite-sample guarantees will be directly stated in terms of their estimation errors.

\section{Structure-agnostic learning bounds for DR covariate shift adaptation}
\label{sec:learning_bounds_general_models}

\indent This section aims to develop finite-sample structure-agnostic learning guarantees for the doubly-robust (DR) estimator. We first state the standing assumptions, introduce the complexity measure utilized in our analysis, and finally present a high-probability bound on the excess $\bbQ$-risk of the DR estimator \eqref{eqn:dr_estimator_v1} together with a concrete illustration based on classes of Frobenius-norm-bounded neural networks.

\indent In this section, we consider the structure-agnostic perspective that treat the given pilot estimates $\big( \hat{\rho}, \hat{f}_0 \big)$ as black-boxes; our bounds depend only on their estimation errors measured by the mean-squared error with respect to $\bbP_{X}$.

\subsection{Assumptions}
\label{subsec:assumptions}

\indent We begin by introducing the minimal assumptions under which our non-asymptotic analysis holds.

\begin{assumption} [Well-specified model]
\label{assumption:well_specified_model}
\normalfont{
$f^* \in \calF$.
}
\end{assumption}

\begin{assumption} [Uniform boundedness]
\label{assumption:uniform_boundedness}
\normalfont{
We have $\sup \left\{ \left\| f \right\|_{\infty} : f \in \calF \right\} \leq 1$ and $\left| Y \right| \leq 1$ almost surely under the source distribution $\bbP$ and the target distribution $\bbQ$.
}
\end{assumption}

\begin{assumption}
\label{assumption:black_box_estimates}
\normalfont{
The pilot estimates $\hat{\rho} : \bbX \to \bbR_{+}$ and $\hat{f}_0 : \bbX \to \bbR$ of the covariate density ratio $\rho^* : \bbX \to \bbR_{+}$ and the shared Bayes regression function $f^* \in \calF$, respectively, satisfy
\begin{equation}
    \label{eqn:uniform_boundedness_black_box_estimates}
    \left\| \hat{\rho} \right\|_{\infty} \leq C_{\textsf{dr}} < +\infty \quad \textnormal{and} \quad \left\| \hat{f}_0 \right\|_{\infty} \leq C_{\textsf{rf}} < +\infty
\end{equation}
for some universal constants $C_{\textsf{dr}}, C_{\textsf{rf}} \in \left( 0, +\infty \right)$.
}
\end{assumption}

\begin{rmk}
\label{rmk:assumptions_finite_sample_analysis}
\normalfont{
We note that the uniform boundedness assumption $\left\| \hat{\rho} \right\|_{\infty} \leq C_{\textsf{dr}} < +\infty$ on the black-box ML estimate $\hat{\rho} : \bbX \to \bbR_{+}$ is standard for the case of the \emph{bounded} ground-truth covariate density ratio $\rho^* : \bbX \to \bbR_{+}$. In particular, the estimation procedures built upon the density ratio matching under the Bregman divergence \cite{sugiyama2012density, sugiyama2012densityratio} including the \emph{least-squares importance fitting (LSIF)} \cite{kanamori2009least}, \emph{kernel mean matching (KMM)} \cite{gretton2009covariate}, \emph{kernel unconstrained LSIF (KuLSIF)} \cite{kanamori2012statistical}, \emph{Kullback-Leibler importance estimation procedure (KLIEP)} \cite{sugiyama2008direct}, logistic regression-based density ratio estimation \cite{sugiyama2012density, sugiyama2012densityratio}, and deep density ratio estimation \cite{kato2021non, zheng2022an}, typically focus on the minimization of a specific empirical risk over a \emph{uniformly bounded} hypothesis class.
}
\end{rmk}

\subsection{Uniform convergence and Rademacher complexity guarantees}
\label{subsec:rademacher_complexity_guarantees}

\indent Now, we turn our attention to analysis of the DR estimator \eqref{eqn:dr_estimator_v1} in finite-sample regimes based on uniform convergence arguments. The key complexity measure is the Rademacher complexity of the \emph{$f^*$-shifted version} of the function class $\mathcal{F} \subseteq \left( \bbX \to \bbR \right)$:
\[
    \calF^* := \left\{ f - f^* : f \in \calF \right\} \subseteq \left( \bbX \to \bbR \right).
\]
We first recall the definition of the Rademacher complexity for completeness.

\begin{defi} [Rademacher complexity]
\label{defi:rademacher_complexity}
\normalfont{
Given any function class $\calG \subseteq \left( \bbX \to \bbR \right)$, the \emph{empirical Rademacher complexity of $\calG$} with respect to $n$ sample points $\bfx_{1:n} = \left( x_1, x_2, \cdots, x_n \right) \in \bbX^n$ is
\begin{equation}
    \label{eqn:defi:rademacher_complexity_v1}
    \begin{split}
        \widehat{\calR}_{n} (\calG) \left( \bfx_{1:n} \right) :=  \bbE_{\bfsigma_{1:n} \sim \textsf{Unif} \left( \left\{ \pm 1 \right\}^n \right)} \left[ \sup \left\{ \left| \frac{1}{n} \sum_{i=1}^{n} \sigma_i g \left( x_i \right) \right|: g \in \calG \right\} \right].
    \end{split}
\end{equation}
The \emph{Rademacher complexity of $\calG$} with respect to a probability measure $\mu \in \Delta (\bbX)$ is defined by
\begin{equation}
    \label{eqn:defi:rademacher_complexity_v2}
    \begin{split}
        \calR_{n}^{\mu} (\calG) := \ & \bbE_{\bfX_{1:n} \sim \mu^{\otimes n}} \left[ \widehat{\calR}_{n} (\calG) \left( \bfX_{1:n} \right) \right] \\
        = \ & \bbE_{\left( \bfX_{1:n}, \bfsigma_{1:n} \right) \sim \mu^{\otimes n} \otimes \textsf{Unif} \left( \left\{ \pm 1 \right\}^n \right)} \left[ \sup \left\{ \left| \frac{1}{n} \sum_{i=1}^{n} \sigma_i g \left( X_i \right) \right|: g \in \calG \right\} \right].
    \end{split}
\end{equation}
}
\end{defi}

\indent With these preliminary notions in hand, one can state a structure-agnostic high-probability upper bound on the excess $\bbQ$-risk of the DR estimator \eqref{eqn:dr_estimator_v1} that depends only on the $L^2$-errors of the pilot estimates and the Rademacher complexity of $\mathcal{F}^*$ under $\bbP_X$ and $\bbQ_X$.

\begin{thm} [Structure-agnostic upper bound \RNum{1} of the DR estimator]
\label{thm:upper_bound_alg:dr_covariate_shift_adaptation_v1}
With Assumptions \ref{assumption:covariate_shift}--\ref{assumption:black_box_estimates}, the doubly-robust (DR) estimator \eqref{eqn:dr_estimator_v1} achieves the $\bbQ$-estimation error
\begin{equation}
    \label{eqn:thm:upper_bound_alg:dr_covariate_shift_adaptation_v1_v1}
    \begin{split}
        \calE_{\bbQ} \left( \hat{f}_{\textnormal{\textsf{DR}}} \right)
        = \ & \bbE_{X \sim \bbQ_{X}} \left[ \left\{ \hat{f}_{\textnormal{\textsf{DR}}} (X) - f^* (X) \right\}^2 \right] \\
        \leq \ & 4 \left\| \hat{\rho} - \rho^* \right\|_{L^2 \left( \bbX, \bbP_{X} \right)} \cdot \left\| \hat{f}_0 - f^* \right\|_{L^2 \left( \bbX, \bbP_{X} \right)} \\
        &+ 12 \left( 2 + C_{\textnormal{\textsf{rf}}} \right) \log \left( \frac{3}{\delta} \right) \left( \frac{C_{\textnormal{dr}}}{n_{\bbP}} + \frac{1}{n_{\bbQ}} \right) \\
        &+ 4 \left( 1 + C_{\textnormal{\textsf{dr}}} \right) \left( 2 + C_{\textnormal{\textsf{rf}}} \right) \sqrt{2 \log \left( \frac{3}{\delta} \right)} \left( \frac{1}{\sqrt{n_{\bbP}}} + \frac{1}{\sqrt{n_{\bbQ}}} \right) \\
        &+ 8 \left( 1 + C_{\textnormal{\textsf{dr}}} \right) \left( 2 + C_{\textnormal{\textsf{rf}}} \right) \sqrt{\log \left( \frac{3}{\delta} \right)} \left( \frac{\calR_{n_{\bbP}}^{\bbP_{X}} \left( \calF^* \right)}{\sqrt{n_{\bbP}}} + \frac{\calR_{n_{\bbQ}}^{\bbQ_{X}} \left( \calF^* \right)}{\sqrt{n_{\bbQ}}} \right) \\
        &+ 8 C_{\textnormal{\textsf{dr}}} \left( 1 + C_{\textnormal{\textsf{rf}}} \right) \calR_{n_{\bbP}}^{\bbP_{X}} \left( \calF^* \right) + 8 \left( 3 + C_{\textnormal{\textsf{rf}}} \right) \calR_{n_{\bbQ}}^{\bbQ_{X}} \left( \calF^* \right)
    \end{split}
\end{equation}
with probability at least $1 - \delta$ under the probability measure $\bbP^{\otimes n_{\bbP}} \otimes \bbQ_{X}^{\otimes n_{\bbQ}}$.
\end{thm}

\noindent The proof of Theorem \ref{thm:upper_bound_alg:dr_covariate_shift_adaptation_v1} is deferred to Appendix \ref{subsec:proof_thm:upper_bound_alg:dr_covariate_shift_adaptation_v1}. Let us make use of the notation
\[
    \begin{split}
        \textsf{Err}_{\rho} := \left\| \hat{\rho} - \rho^* \right\|_{L^2 \left( \bbX, \bbP_{X} \right)} \quad \textnormal{and} \quad \textsf{Err}_{f} := \left\| \hat{f}_0 - f^* \right\|_{L^2 \left( \bbX, \bbP_{X} \right)}.
    \end{split}
\]
The leading bias term in \eqref{eqn:thm:upper_bound_alg:dr_covariate_shift_adaptation_v1_v1} can be rewritten as the product $\textnormal{Err}_{\rho} \cdot \textnormal{Err}_{f}$. This key observation leads to the following two concrete implications:

\begin{enumerate} [label = (I\arabic*)]
    \item \textbf{Having just one good pilot estimate suffices.} By assuming either $\textsf{Err}_{\rho} = o(1)$ or $\textsf{Err}_{f} = o(1)$ as $\min \left\{ n_{\bbP}, n_{\bbQ} \right\} \to \infty$ and the remaining term is bounded, one can obtain $\textsf{Err}_{\rho} \cdot \textsf{Err}_{f} = o(1)$. Hence, the DR estimator is still consistent even when one of the pilot estimates is inaccurate; this is the finite‐sample manifestation of the \emph{double robustness} phenomenon \cite{robins1995semiparametric, robins2008higher}: the error of the one-step corrected estimators is upper bounded by a product of estimation errors of the underlying nuisance components. To put it another way, the DR covariate shift adaptation allows us to reduce the bias incurred by the estimation error of the covariate density ratio through the aforementioned double robustness property. 
    \item \textbf{Rate multiplication.} Suppose $\textsf{Err}_{\rho} = \tilde{\calO} \left( n^{- \alpha} \right)$ and $\textsf{Err}_{f} = \tilde{\calO} \left( n^{- \beta} \right)$, where $n := \min \left\{ n_{\bbP}, n_{\bbQ} \right\}$. Then, the order of their product term scales as $\tilde{\calO} \left( n^{- \left( \alpha + \beta \right)} \right)$. In contrast, the upper bound on the excess $\bbQ$-risk of the IW estimator depends additively on $\textsf{Err}_{\rho}$, and typically requires $\alpha \geq 1/2$ to be competitive. Thus, the DR estimator \eqref{eqn:dr_estimator_v1} outperforms the IW method whenever $\alpha+\beta > 1/2$.
\end{enumerate}

\indent Since we have trivial bounds $\calR_{n_{\bbP}}^{\bbP_X} \left( \calF^* \right) \leq 2$ and $\calR_{n_{\bbQ}}^{\bbQ_X} \left( \calF^* \right) \leq 2$, one can simplify the excess $\bbQ$-risk bound \eqref{eqn:thm:upper_bound_alg:dr_covariate_shift_adaptation_v1_v1} in Theorem \ref{thm:upper_bound_alg:dr_covariate_shift_adaptation_v1} of the DR estimator as follows: With Assumptions \ref{assumption:covariate_shift}--\ref{assumption:black_box_estimates} in hand, it follows that the DR estimator \eqref{eqn:dr_estimator_v1} achieves
\begin{equation}
    \label{eqn:simpification_excess_risk_dr_estimator}
    \begin{split}
        \calE_{\bbQ} \left( \hat{f}_{\textsf{DR}} \right) \lesssim \ & \left\| \hat{\rho} - \rho^* \right\|_{L^2 \left( \bbX, \bbP_{X} \right)} \cdot \left\| \hat{f}_0 - f^* \right\|_{L^2 \left( \bbX, \bbP_{X} \right)} + \sqrt{\frac{\log \left( \frac{1}{\delta} \right)}{n_{\bbP}}} + \sqrt{\frac{\log \left( \frac{1}{\delta} \right)}{n_{\bbQ}}} \\
        &+ \calR_{n_{\bbP}}^{\bbP_X} \left( \calF^* \right) + \calR_{n_{\bbQ}}^{\bbQ_X} \left( \calF^* \right)
    \end{split}
\end{equation}
with probability at least $1 - \delta$.

\subsection{An illustration with Frobenius-norm-bounded neural networks}
\label{subsec:illustration_class_neural_networks}

\indent Let $\bbX \subseteq \bbR^{n_0}$ be a bounded domain such that $\sup \left\{ \left\| \bfx \right\|_{2} : \bfx \in \bbX \right\} \leq R$ for some radius $R \in \left( 0, +\infty \right)$. We also consider a collection of $1$-Lipschitz activation functions $\left\{ \sigma_j \in \left( \bbR \to \bbR \right) : j \in \bbN \right\}$ that are \emph{positive-homogeneous} (i.e., $\sigma_{j} \left( \alpha t \right) = \alpha \sigma_j (t)$ for any $\left( \alpha, t \right) \in \bbR_{+} \times \bbR$), and that are applied \emph{element-wise}. We are mainly interested in a class of real-valued neural networks of depth $d \in \bbN$ over the domain $\bbX \subseteq \bbR^{n_0}$ defined as
\begin{equation}
    \label{eqn:class_nn_v1}
    \begin{split}
        \calH_d \left( \bbX; M_{\textsf{F}} \right) := \left\{ \textsf{NN}_{d} \left( \cdot; \bftheta \right) \in \left( \bbX \to \bbR \right): \bftheta \in \Theta \left( M_{\textsf{F}} \right) \right\},
    \end{split}
\end{equation}
where $\bftheta = \left( \bfW_1, \cdots, \bfW_d \right) \in \prod_{j=1}^{d} \bbR^{n_{j} \times n_{j-1}}$ denotes the model parameter consists of $d$ parameter matrices with $n_d = 1$, and the real-valued neural network $\textsf{NN}_{d} \left( \cdot; \bftheta \right): \bbX \to \bbR$ of depth $d$ is defined to be
\begin{equation}
    \label{eqn:class_nn_v2}
    \begin{split}
        \textsf{NN}_{d} \left( \bfx; \bftheta \right) := \bfW_{d} \sigma_{d-1} \left( \bfW_{d-1} \sigma_{d-2} \left( \cdots\sigma_{1} \left( \bfW_{1} \bfx \right) \cdots \right) \right).
    \end{split}
\end{equation}
Here, $M_{\textsf{F}} : [d] \to \bbR_{+}$ specifies upper bounds on the Frobenius norm of parameter matrices, and the parameter space $\Theta \left( M_{\textsf{F}} \right) \subseteq \prod_{j=1}^{d} \bbR^{n_{j} \times n_{j-1}}$ is given by
\[
    \begin{split}
        \Theta \left( M_{\textsf{F}} \right) := \left\{ \bftheta = \left( \bfW_1, \bfW_2, \cdots, \bfW_d \right) \in \prod_{j=1}^{d} \bbR^{n_{j} \times n_{j-1}} : \left\| \bfW_{j} \right\|_{\textsf{F}} \leq M_{\textsf{F}} (j),\ \forall j \in [d] \right\}.
    \end{split}
\]
A prominent example of the above construction are ReLU networks, where every $\sigma_{j}: \bbR \to \bbR$ corresponds to applying the ReLU activation function $\sigma (\cdot) := \max \left\{ 0, \cdot \right\}: \bbR \to \bbR_{+}$. Armed with the class $\calH_d \left( \bbX; M_{\textsf{F}} \right)$, let us now introduce the function class of our interest. Let $\eta : \bbR \to \left[ -1, 1 \right]$ be an $L$-Lipschitz bounded activation function such that $\eta (0) = 0$, and define
\begin{equation}
    \label{eqn:class_nn_v3}
    \begin{split}
        \calF := \left\{ f \left( \cdot; \bftheta \right) := \eta \circ \textsf{NN}_d \left( \cdot; \bftheta \right) \in \left( \bbX \to \left[ -1, 1 \right] \right) : \bftheta \in \Theta \left( M_{\textsf{F}} \right) \right\}.
    \end{split}
\end{equation}
For example, the \emph{inverse tangent activation} function $\frac{2}{\pi} \arctan (\cdot) : \bbR \to \left[ -1, 1 \right]$ satisfies the desired properties with $L = \frac{2}{\pi}$. One can show that the Rademacher complexity of the $f^*$-shift version of \eqref{eqn:class_nn_v3} with respect to any probability measure $\mu \in \Delta (\bbX)$ is of order $\calO \left( \frac{1}{\sqrt{n}} \right)$.

\begin{prop}
\label{prop:rademacher_complexity_neural_networks}
The Rademacher complexity of the $f^*$-shifted version of the neural network class defined as \eqref{eqn:class_nn_v3}, $\calF^* := \calF - \left\{ f^* \right\}$, with respect to any given probability measure $\mu \in \Delta (\bbX)$ is upper bounded by
\begin{equation}
    \label{eqn:prop:rademacher_complexity_neural_networks_v1}
    \begin{split}
        \calR_{n}^{\mu} \left( \calF^* \right) \leq \frac{2}{\sqrt{n}} \left\{ L R \left( 1 + \sqrt{\left( 2 \log 2 \right) d} \right) \prod_{j=1}^{d} M_{\textnormal{\textsf{F}}} (j) + \sqrt{\log 2} \right\} = \calO \left( \frac{1}{\sqrt{n}} \right).
    \end{split}
\end{equation}
\end{prop}

\noindent The proof of Proposition \ref{prop:rademacher_complexity_neural_networks} can be found in Appendix \ref{subsec:proof_prop:rademacher_complexity_neural_networks}. With Proposition \ref{prop:rademacher_complexity_neural_networks} in hand, one may conclude that the DR estimator \eqref{eqn:dr_estimator_v1} achieves the following excess $\bbQ$-risk bound when we select the hypothesis class $\calF \subseteq \left( \bbX \to \left[ -1, 1 \right] \right)$ of our interest as \eqref{eqn:class_nn_v3}: with probability at least $1 - \delta$, one has
\begin{equation}
    \label{eqn:excess_Q_risk_class_nn}
    \begin{split}
        \calE_{\bbQ} \left( \hat{f}_{\textsf{DR}} \right) \lesssim \left\| \hat{\rho} - \rho^* \right\|_{L^2 \left( \bbX, \bbP_{X} \right)} \cdot \left\| \hat{f}_0 - f^* \right\|_{L^2 \left( \bbX, \bbP_{X} \right)} + \sqrt{\frac{\log \left( \frac{1}{\delta} \right)}{n_{\bbP}}} + \sqrt{\frac{\log \left( \frac{1}{\delta} \right)}{n_{\bbQ}}}.
    \end{split}
\end{equation}

\begin{rmk}
\label{rmk:excess_Q_risk_class_nn}
\normalfont{
We now turn our attention to the following approach that utilizes the double/debiased machine learning (DML) technique \cite{chernozhukov2017double, chernozhukov2018double, chernozhukov2022locally, chernozhukov2023simple, foster2023orthogonal}: We first split the observed data $\calD := \left( \bfO_{1:n_{\bbP}}^{\bbP}, \bfX_{1:n_{\bbQ}}^{\bbQ} \right)$ into two subgroups $\calD_1$ and $\calD_2$ with the equal size, and then estimate the ground-truth covariate density ratio $\rho^* : \bbX \to \bbR_{+}$ and the common Bayes regression function $f^* \in \calF$ utilizing the first subgroup $\calD_1$ to compute a nuisance estimate $\hat{\rho} : \bbX \to \bbR_{+}$ and a pilot estimate $\hat{f}_0 : \bbX \to \bbR$. A number of results from the literature of density ratio estimation propose algorithms achieving $\left\| \hat{\rho} - \rho^* \right\|_{L^2 \left( \bbX, \bbP_{X} \right)} = \calO_{p} \left( \min \left\{ n_{\bbP}, n_{\bbQ} \right\}^{- \frac{1}{2 + \gamma}} \right)$ as $\min \left\{ n_{\bbP}, n_{\bbQ} \right\} \to \infty$ for any constant $\gamma \in (0, 2)$ \cite{kanamori2012statistical, kato2021non}. Therefore, if the pilot estimate $\hat{f}_0 : \bbX \to \bbR$ of  $f^* \in \calF$ is consistent under the source distribution $\bbP$ with a rate
\begin{equation}
    \label{eqn:rmk:excess_Q_risk_class_nn_v1}
    \begin{split}
        \left\| \hat{f}_0 - f^* \right\|_{L^2 \left( \bbX, \bbP_{X} \right)} = \calO_{p} \left( \min \left\{ n_{\bbP}, n_{\bbQ} \right\}^{- \frac{\gamma}{2 \left( 2 + \gamma \right)}} \right) \quad \textnormal{as} \quad \min \left\{ n_{\bbP}, n_{\bbQ} \right\} \to \infty,
    \end{split}
\end{equation}
then the high-probability bound \eqref{eqn:excess_Q_risk_class_nn} on the $\bbQ$-risk for the DR estimator \eqref{eqn:dr_estimator_v1} together with the class \eqref{eqn:class_nn_v3} of Frobenius-norm-bounded neural networks (constructed using the second subgroup $\calD_2$) gives
\begin{equation}
    \label{eqn:rmk:excess_Q_risk_class_nn_v2}
    \begin{split}
        \calE_{\bbQ} \left( \hat{f}_{\textsf{DR}} \right) \lesssim \sqrt{\frac{\log \left( \frac{1}{\delta} \right)}{n_{\bbP}}} + \sqrt{\frac{\log \left( \frac{1}{\delta} \right)}{n_{\bbQ}}}.
    \end{split}
\end{equation}
To summarize, as long as the pilot estimate $\hat{f}_0 : \bbX \to \bbR$ for the Bayes regression function $f^* \in \calF$ is consistent under the source distribution $\bbP$ with the rate of convergence \eqref{eqn:rmk:excess_Q_risk_class_nn_v1}, we are able to enhance the pilot estimate $\hat{f}_0 : \bbX \to \bbR$ of $f^* \in \calF$ to an estimator that achieves the rate of convergence \eqref{eqn:rmk:excess_Q_risk_class_nn_v2} even if it is not consistent under the target distribution $\bbQ$.
}
\end{rmk}

\section{Learning bounds for DR covariate shift adaptation: parametric models}
\label{sec:learning_bounds_parametric_models}

\indent This section closely examines the doubly-robust (DR) covariate shift adaptation when the underlying function class is finite-dimensional and well-specified. Our central takeaway messages are two-fold: (i) with parametric models, \emph{fast} $1/n$-type rates of convergence are attainable without assuming exact knowledge of the covariate density ratio $\rho^*$; and (ii) the DR estimator achieves these rates \emph{regardless of the statistical accuracies} of the pilot estimates $\big( \hat{\rho}, \hat{f}_0 \big)$.

\paragraph{Parametric model.}
Throughout this section, we impose Assumptions \ref{assumption:well_specified_model} and \ref{assumption:uniform_boundedness} and consider a $d$-dimensional parameterization
\begin{equation}
    \label{eqn:parametrized_hypothesis_class}
    \calF = \left\{ f \left( \cdot; \bftheta \right) \in \left( \bbX \to [-1,1] \right) : \bftheta \in \Theta \subseteq \bbR^d \right\},
\end{equation}
with the ground-truth parameter $\bftheta^* \in \Theta$ such that $f^* (\cdot)= f \left( \cdot; \bftheta^* \right) \in \calF$. For any pilot estimates $\big( \hat{\rho},\hat{f}_0 \big)$, the DR empirical risk specialized to the parameterized model \eqref{eqn:parametrized_hypothesis_class} is
\begin{equation}
    \label{eqn:dr_empirical_risk_v2}
    \begin{split}
        \widehat{\calR}_{\textsf{DR}} (\bftheta) := \ & \frac{1}{n_{\bbP}} \sum_{i=1}^{n_{\bbP}} \hat{\rho} \big( X_{i}^{\bbP} \big) \left\{ \ell \left( Y_{i}^{\bbP}, f \big( X_{i}^{\bbP}; \bftheta \big) \right) - \ell \left( \hat{f}_0 \big( X_{i}^{\bbP} \big), f \big( X_{i}^{\bbP}; \bftheta \big) \right) \right\} \\
        &+ \frac{1}{n_{\bbQ}} \sum_{j=1}^{n_{\bbQ}} \ell \left( \hat{f}_0 \big( X_{j}^{\bbQ} \big), f \big( X_{j}^{\bbQ}; \bftheta \big) \right),
    \end{split}
\end{equation}
where $\ell (a, b) := (b-a)^2$. We define the DR estimator specialized to the parametric model \eqref{eqn:parametrized_hypothesis_class} as
\begin{equation}
    \label{eqn:dr_estimator_parameterized_case}
    \begin{split}
        \hat{\bftheta}_{\textsf{DR}} \in \argmin \left\{ \widehat{\calR}_{\textsf{DR}} (\bftheta) : \bftheta \in \Theta \right\}, \quad \hat{f}_{\textsf{DR}} (\cdot) := f \left( \cdot; \hat{\bftheta}_{\textsf{DR}} \right) \in \calF.
    \end{split}
\end{equation}

\paragraph{Regularity and landscape conditions.}
We first make the following smoothness assumptions customary in classical analysis of MLE \cite{lehmann1999elements, le1956asymptotic, cramer1999mathematical, van2000asymptotic, lehmann2006theory}.

\begin{assumption} [Smoothness assumptions]
\label{assumption:smoothness}
\normalfont{
Suppose the parameter space $\Theta \subseteq \bbR^d$ is \emph{star-shaped at center $\bftheta^* \in \Theta$}, i.e., $\left[ \bftheta^* , \bftheta \right] := \left\{ \bftheta^* + \lambda \left( \bftheta - \bftheta^* \right) : \lambda \in [0, 1] \right\} \subseteq \Theta$ for all $\bftheta \in \Theta$, and
\begin{enumerate} [label = (\roman*)]
    \item For each $x \in \bbX$, the function $\bftheta \in \Theta \mapsto f \left( x; \bftheta \right) \in \left[ -1, 1 \right]$ is three-times differentiable;
    \item There exist absolute constants $\left( b_1, b_2, b_3 \right) \in \left( 0, +\infty \right)^3$ such that
    \begin{equation}
        \label{eqn:assumption:smoothness_v1}
        \begin{split}
            \left\| \nabla_{\bftheta} f \left( x; \bftheta \right) \right\|_{2} \leq b_1, \quad \left\| \nabla_{\bftheta}^2 f \left( x; \bftheta \right) \right\|_{\textsf{op}} \leq b_2, \quad \textnormal{and} \left\| \nabla_{\bftheta}^3 f \left( x; \bftheta \right) \right\|_{\textsf{op}} \leq b_3
        \end{split}
    \end{equation}
    for every $\left( x, \bftheta \right) \in \bbX \times \Theta$.
\end{enumerate}
}
\end{assumption}

\begin{assumption} [Benign landscape of the DR empirical risk]
\label{assumption:simple_landscape_empirical_DR_risk}
\normalfont{
For any realization $\left( \bfO^{\bbP}_{1:n_{\bbP}}, \bfX^{\bbQ}_{1:n_{\bbQ}} \right) \in \bbO^{n_{\bbP}} \times \bbX^{n_{\bbQ}}$, the DR empirical risk $\widehat{\calR}_{\textsf{DR}} : \Theta \to \bbR$ attains a unique local minimum, which is also the global minimum.
}
\end{assumption}

\noindent Here, we note that Assumption \ref{assumption:simple_landscape_empirical_DR_risk} is satisfied, for example, if the population version of the DR empirical risk \eqref{eqn:dr_empirical_risk_v2} is strongly convex in an open neighborhood of $\bftheta^* \in \Theta$, and the Hessian of the DR empirical risk \eqref{eqn:dr_empirical_risk_v2} uniformly concentrates on that neighborhood.

\indent Now, we are ready to establish an improved structure-agnostic learning bound for the DR estimator \eqref{eqn:dr_estimator_parameterized_case} for parametrized hypothesis classes \eqref{eqn:parametrized_hypothesis_class}, which leads to faster rates of convergence. For convenience, let us first recall the classical notion of Fisher information, which plays a critical role as a key quantity to measure the difficulty of parameter estimation. The \emph{$\mu$-Fisher information matrix evaluated at $\bftheta \in \Theta$} is defined as
\begin{equation}
    \label{eqn:fisher_information_matrix}
    \begin{split}
        \calI_{\mu} (\bftheta) := \bbE_{( X, Y ) \sim \mu} \left[ \nabla_{\bftheta}^2 \ell \left( Y, f \left( X; \bftheta \right) \right) \right], \quad \bftheta \in \Theta,
    \end{split}
\end{equation}
where $\mu \in \left\{ \bbP, \bbQ \right\}$ and $\ell : \bbR \times \bbR \to \bbR_{+}$ is the squared error loss. Then, one can easily observe that
\begin{equation}
    \label{eqn:fisher_information_matrix_ground_truth_parameter}
    \begin{split}
        \calI_{\mu} \left( \bftheta^* \right) := 2 \bbE_{X \sim \mu_{X}} \left[ \nabla_{\bftheta} f \left( X; \bftheta^* \right) \left\{ \nabla_{\bftheta} f \left( X; \bftheta^* \right) \right\}^{\top} \right],
    \end{split}
\end{equation}
where $\mu_{X} (\cdot) := \mu \left( \cdot \times \bbR \right) \in \Delta (\bbX)$ refers to the covariate marginal distribution of $\mu \in \left\{ \bbP, \bbQ \right\}$. We now present our main result of this section, whose proof is deferred to Section \ref{subsec:proof_thm:upper_bound_alg:dr_covariate_shift_adaptation_v2}:

\begin{thm} [Informal, see Theorem \ref{thm:detailed_upper_bound_alg:dr_covariate_shift_adaptation_v2}]
\label{thm:upper_bound_alg:dr_covariate_shift_adaptation_v2}
With the parametrized function class \eqref{eqn:parametrized_hypothesis_class}, under Assumptions \ref{assumption:covariate_shift}--\ref{assumption:simple_landscape_empirical_DR_risk}, there is an absolute constant $K \in \left( 0, +\infty \right)$ such that with probability at least $1 - 8 \delta$ under the probability measure $\bbP^{\otimes n_{\bbP}} \otimes \bbQ_{X}^{\otimes n_{\bbQ}}$, 
\begin{equation}
    \label{eqn:thm:upper_bound_alg:dr_covariate_shift_adaptation_v2_v1}
    \begin{split}
        &\calE_{\bbQ} \left( \hat{\bftheta}_{\textnormal{\textsf{DR}}} \right) = \bbE_{X \sim \bbQ_{X}} \left[ \left\{ f \left( X; \hat{\bftheta}_{\textnormal{\textsf{DR}}} \right) - f^* (X) \right\}^2 \right] \\
        \leq \ & 18 K^2 \left( 1 + C_{\textnormal{\textsf{dr}}} \right)^2 \left( 1 + C_{\textnormal{\textsf{rf}}} \right)^2 \log \left( \frac{d}{\delta} \right) \left[ \frac{\textnormal{\textsf{Trace}} \left\{ \calI_{\bbP} \left( \bftheta^* \right) \calI_{\bbQ}^{-1} \left( \bftheta^* \right) \right\}}{n_{\bbP}} + \frac{d}{n_{\bbQ}} \right],
    \end{split}
\end{equation}
provided that $\min \left\{ n_{\bbP}, n_{\bbQ} \right\} \geq \overline{\kappa} \cdot \calN^* \log \left( \frac{d}{\delta} \right)$ for some absolute constant $\overline{\kappa} \in \left( 0, +\infty \right)$, where 
\[
    \calN^* = \textnormal{poly} \left( d, \left\| \calI_{\bbQ}^{-1} \left( \bftheta^* \right) \right\|_{\textnormal{\textsf{op}}}, \left\| \calI_{\bbQ}^{-1} \left( \bftheta^* \right) \calI_{\bbP} \left( \bftheta^* \right) \calI_{\bbQ}^{-1} \left( \bftheta^* \right) \right\|_{\textnormal{\textsf{op}}} \right).
\]
\end{thm}

\paragraph{Interpretations \& key implications}
Theorem \ref{thm:upper_bound_alg:dr_covariate_shift_adaptation_v2} shows that, for well-specified parametric models, the DR estimator \eqref{eqn:dr_estimator_parameterized_case} achieves a \emph{fast} and \emph{instance-dependent} upper bound on the excess $\bbQ$-risk that decouples the contributions of the source and target samples to the bound:
\[
    \frac{\textsf{Trace} \left\{ \calI_{\bbP} \left( \bftheta^* \right) \calI_{\bbQ}^{-1} \left( \bftheta^* \right) \right\}}{n_{\bbP}} : \textnormal{contributed by the source data}, \quad \frac{d}{n_{\bbQ}} : \textnormal{contributed by the target data},
\]
up to logarithmic factors. Here, the trace factor $\textsf{Trace} \big\{ \calI_{\bbP} \big( \bftheta^* \big) \calI_{\bbQ}^{-1} \big( \bftheta^* \big) \big\}$ quantifies the \emph{Fisher information mismatch} between $\bbP$ and $\bbQ$, and is the only way in which covariate shift affects the leading constant. Notably, the excess $\bbQ$-risk bound in Theorem \ref{thm:upper_bound_alg:dr_covariate_shift_adaptation_v2} holds \emph{without} access to the ground-truth covariate density ratio $\rho^*$, and is \emph{independent of the statistical accuracies of the pilot estimates} $\big( \hat{\rho}, \hat{f}_0 \big)$.

\indent We also discuss some appealing attributes of DR covariate shift adaptation and its fast $1/n$-type convergence guarantee \eqref{eqn:thm:upper_bound_alg:dr_covariate_shift_adaptation_v2_v1} for well-specified parametric models provided in Theorem \ref{thm:upper_bound_alg:dr_covariate_shift_adaptation_v2}:
\begin{itemize}
    \item \emph{Fast rates of convergence under covariate shift without knowing $\rho^*$}: The excess $\bbQ$-risk bound \eqref{eqn:thm:upper_bound_alg:dr_covariate_shift_adaptation_v2_v1} of the DR estimator \eqref{eqn:dr_estimator_parameterized_case} matches the fast $1/n$-rate behavior, where $n := \min \left\{ n_{\bbP}, n_{\bbQ} \right\}$, which is known to be achievable in parametric models, yet it does so \emph{without} requiring an exact knowledge (or a consistent estimate) of the covariate density ratio $\rho^* : \bbX \to \bbR_{+}$.
    \item \emph{Pilot-agnostic tightness of the rates of convergence}: The rate of convergence for the DR estimator \eqref{eqn:dr_estimator_parameterized_case} obtained from the excess $\bbQ$-risk bound \eqref{eqn:thm:upper_bound_alg:dr_covariate_shift_adaptation_v2_v1} does not degrade with the quality of given pilot estimates $\big( \hat{\rho}, \hat{f}_0 \big)$; any black-box pilot estimates suffice.
    \item \emph{No boundedness assumption on the covariate density ratio $\rho^*$}: Unlike the prior works on covariate shift (e.g., \cite{cortes2010learning, ma2023optimally}), we make no boundedness assumptions on the true covariate density ratio $\rho^* : \bbX \to \bbR$, broadening applicability of our results.
\end{itemize}

\indent It would be worth pointing out the trace factor $\textsf{Trace} \big\{ \calI_{\bbP} \big( \bftheta^* \big) \calI_{\bbQ}^{-1} \big( \bftheta^* \big) \big\}$, which is different from the trace factors that appears in the excess $\bbQ$-risk bounds for the vanilla MLE and the weighted MLE of \cite{ge2024maximum}. However, on the closer look, \cite{ge2024maximum} assumes the boundedness of the covariate density ratio $\rho^*$, under which their excess $\bbQ$-risk bound for the weighted MLE (see \emph{Theorem 5.2} therein) can be translated to the same trace factor as in the bound \eqref{eqn:thm:upper_bound_alg:dr_covariate_shift_adaptation_v2_v1} of Theorem \ref{thm:upper_bound_alg:dr_covariate_shift_adaptation_v2}.

\section{Discussion}
\label{sec:discussion}

\indent This paper establishes the first finite-sample guarantees for doubly-robust (DR) covariate shift adaptation, complementing the prior asymptotic analysis \cite{kato2023double} and clarifying the role of pilot estimates, sample allocation, and parametric modeling for the Bayes regression function. The structure-agnostic upper bound \eqref{eqn:thm:upper_bound_alg:dr_covariate_shift_adaptation_v1_v1} of the DR estimator \eqref{eqn:dr_estimator_v1} shows that the leading bias term scales as the product of statistical error rates for the pilot estimates, providing a non-asymptotic demonstration of the celebrated double robustness phenomenon \cite{robins1995semiparametric, robins2008higher}: one consistent pilot estimate suffices to obtain the consistency of the one-step corrected estimators, and the joint improvement leads us to multiplicative gains. The decomposition of the DR empirical risk \eqref{eqn:dr_empirical_risk_v1} underscores how the labeled source samples primarily benefit the pilot regression model, while the unlabeled target covariates strengthen the effect of the pilot estimate for the covariate density ratio, offering practical guidance on data collection under budget constraints in the target domain. Within well-specified parametric models, our analysis of the DR estimator \eqref{eqn:dr_estimator_parameterized_case} via modern techniques for finite-sample analysis of MLE yields a non-asymptotic fast $1/n$-type convergence guarantee, which is independent of the statistical accuracies of pilot black-box estimates. In this result, the difficulty of learning a predictive model under covariate shift is quantified by the Fisher information mismatch term between the source and target distributions. Together, the findings in this paper demonstrate that the DR covariate shift adaptation combines asymptotic efficiency results with strong finite-sample out-of-distribution generalization bounds.

\section*{Acknowledgements}

Jeonghwan Lee is partially supported by the Doctoral Overseas Scholarship from the Kwanjeong Educational Foundation. Cong Ma is partially supported by the National Science Foundation via grant DMS-2311127 and the CAREER Award DMS-2443867. 

\newpage

\bibliographystyle{plain}

\newpage

\appendix


\section{Preliminary facts}
\label{sec:preliminary_facts}

\indent In this section, let us collect a couple of useful preliminary facts that facilitates our analysis. The following contraction lemma is a modification of \emph{Theorem 4.12} of \cite{ledoux2013probability} that has been established in \cite{duchi2009probability}. See \emph{Theorem 7} therein for the proof of Lemma \ref{lemma:ledoux_talagrand_contraction}.

\begin{lemma} [The Ledoux-Talagrand contraction principle]
\label{lemma:ledoux_talagrand_contraction}
Let $f: \bbR_{+} \to \bbR_{+}$ be any non-decreasing convex function, and $\phi_{i}: \bbR \to \bbR$, $i \in [n]$, are $L$-Lipschitz continuous functions such that $\phi_{i} (0) = 0$. Then, it holds for any $T \subseteq \bbR^n$ that
\[
    \begin{split}
        \bbE_{\bfsigma_{1:n} \sim \textnormal{\textsf{Unif}} \left( \left\{ \pm 1 \right\}^n \right)} \left[ f \left\{ \frac{1}{2} \sup_{\bft_{1:n} \in T} \left| \sum_{i=1}^{n} \sigma_i \phi_i \left( t_i \right) \right| \right\} \right] \leq \bbE_{\bfsigma_{1:n} \sim \textnormal{\textsf{Unif}} \left( \left\{ \pm 1 \right\}^n \right)} \left[ f \left( L \cdot \sup_{\bft_{1:n} \in T} \left| \sum_{i=1}^{n} \sigma_i t_i \right| \right) \right].
    \end{split}
\]
In particular, if we let $f (t) = t$ for $t \in \bbR_{+}$, then we obtain
\begin{equation}
    \label{eqn:lemma:ledoux_talagrand_contraction_v1}
    \begin{split}
        \bbE_{\bfsigma_{1:n} \sim \textnormal{\textsf{Unif}} \left( \left\{ \pm 1 \right\}^n \right)} \left[ \sup_{\bft_{1:n} \in T} \left| \frac{1}{n} \sum_{i=1}^{n} \sigma_i \phi_i \left( t_i \right) \right| \right] \leq 2 L \bbE_{\bfsigma_{1:n} \sim \textnormal{\textsf{Unif}} \left( \left\{ \pm 1 \right\}^n \right)} \left[ \sup_{\bft_{1:n} \in T} \left| \frac{1}{n} \sum_{i=1}^{n} \sigma_i t_i \right| \right].
    \end{split}
\end{equation}
\end{lemma}

\noindent The following is a well-known standard deviation inequality for controlling the maxima of empirical processes; see \emph{Theorem 1.1} in \cite{klein2005concentration}.

\begin{lemma} [Classical Talagrand's concentration inequality]
\label{lemma:talagrand_concentration}
Let $\calF \subseteq \left( \bbX \to \left[ -B, B \right] \right)$ be any function class and $\bfX_{1:n} = \left( X_1, X_2, \cdots, X_n \right) \sim \bbP^{\otimes n}$ for some $\bbP \in \Delta (\bbX)$. We define
\[
    Z := \sup \left\{ \left( \widehat{\bbP} - \bbP \right) (f) := \frac{1}{n} \sum_{i=1}^{n} f \left( X_i \right) - \bbE_{X \sim \bbP} \left[ f (X) \right] : f \in \calF \right\},
\]
and $v^2 := \sup \left\{ \textnormal{\textsf{Var}}_{X \sim \bbP} \left[ f (X) \right]: f \in \calF \right\}$, where $\widehat{\bbP} := \frac{1}{n} \sum_{i=1}^{n} \delta_{X_i} \in \Delta (\bbX)$ stands for the empirical measure for the $n$ samples $\bfX_{1:n} \sim \bbP^{\otimes n}$. Then, it holds for every $x \in \bbR_{+}$ that
\begin{equation}
    \label{eqn:lemma:talagrand_concentration_v1}
    \begin{split}
        \bbP \left\{ Z > \bbE [Z] + x \right\} \leq \exp \left( - \frac{n x^2}{4 B \bbE [Z] + 2 v^2 + 3 B x} \right).
    \end{split}
\end{equation}
In particular, for any given $\delta \in (0, 1)$, it holds with probability at least $1 - \delta$ that
\begin{equation}
    \label{eqn:lemma:talagrand_concentration_v2}
    \begin{split}
        Z - \bbE [Z] \leq \frac{3 B \log \left( \frac{1}{\delta} \right)}{n} + 2 \sqrt{\frac{B \bbE [Z] \log \left( \frac{1}{\delta} \right)}{n}} + \sqrt{\frac{2 v^2 \log \left( \frac{1}{\delta} \right)}{n}}
    \end{split}
\end{equation}
under the probability measure $\bbP^{\otimes n}$.
\end{lemma}

\noindent Another key technical result is the following generic version of the Bernstein inequality for random vectors, which plays a crucial role in establishing concentration properties for the gradient of the DR empirical risk \eqref{eqn:dr_empirical_risk_v2} with respect to the parameter vector $\bftheta$. Check Lemma \ref{lemma:concentration_gradient} in Section \ref{subsec:proof_lemma:subsec:proof_thm:upper_bound_alg:dr_covariate_shift_adaptation_v2_v1} for further details.

\begin{lemma}
\label{lemma:generic_Bernstein_inequality}
Suppose $\bbP \in \Delta \left( \bbR^d \right)$ satisfies $\bbE_{\bfX \sim \bbP} \left[ \bfX \right] = \bfzero_d$ and $\calV := \bbE_{\bfX \sim \bbP} \left[ \left\| \bfX \right\|_{2}^2 \right] < +\infty$. Define
\begin{equation}
    \label{eqn:lemma:generic_Bernstein_inequality_v1}
    \begin{split}
        \calB (\alpha) := \inf \left\{ t \in \left( 0, +\infty \right): \bbE_{\bfX \sim \bbP} \left[ \exp \left\{ \left( \frac{\left\| \bfX \right\|_{2}}{t} \right)^{\alpha} \right\} \right] \leq 2 \right\}, \quad \alpha \in \left[ 1, +\infty \right),
    \end{split}
\end{equation}
and assume that $\calB (\alpha) < +\infty$ for some constant $\alpha \in \left[ 1, +\infty \right)$. Then, there exists an absolute constant $C > 0$ such that for any given $\delta \in (0, 1)$, we have
\begin{equation}
    \label{eqn:lemma:generic_Bernstein_inequality_v2}
    \begin{split}
        \left\| \frac{1}{n} \sum_{i=1}^{n} \bfX_i \right\|_{2} \leq C \left[ \sqrt{\frac{\calV \log \left( \frac{d}{\delta} \right)}{n}} + \calB (\alpha) \log^{\frac{1}{\alpha}} \left\{ \frac{\calB (\alpha)}{\sqrt{\calV}} \right\} \frac{\calV \log \left( \frac{d}{\delta} \right)}{n} \right]
    \end{split}
\end{equation}
under $\left( \bfX_{1}, \bfX_{2}, \cdots, \bfX_{n} \right) \sim \bbP^{\otimes n}$, with probability at least $1 - \delta$.
\end{lemma}

\noindent We refer to \emph{Proposition 2} in \cite{koltchinskii2011nuclear} for the proof of Lemma \ref{lemma:generic_Bernstein_inequality}.
\medskip

\indent Lastly, the following lemma gives a standard upper bound on the Rademacher complexity of finite hypothesis classes.

\begin{lemma}
\label{lemma:rademacher_complexity_finite_class}
Let $\calF \subseteq \left( \bbX \to \left[ -B, B \right] \right)$ be a finite function class, i.e., $| \calF | < +\infty$. Then, it holds that
\begin{equation}
    \label{eqn:lemma:rademacher_complexity_finite_class_v1}
    \begin{split}
        \calR_{n}^{\mu} (\calF) \leq 2 B \sqrt{\frac{\log \left( 2 |\calF| \right)}{n}},
    \end{split}
\end{equation}
for any probability measure $\mu \in \Delta (\bbX)$. 
\end{lemma}

\section{Proofs for Section \ref{sec:dr_covariate_shift_adaptation}}
\label{sec:proofs_sec:dr_covariate_shift_adaptation}

\subsection{Proof of Proposition \ref{prop:rademacher_complexity_neural_networks}}
\label{subsec:proof_prop:rademacher_complexity_neural_networks}

We first observe for any model parameter $\bftheta = \left( \bfW_1, \cdots, \bfW_d \right) \in \Theta \left( M_{\textsf{F}} \right)$ that
\[
    \begin{split}
        - \textsf{NN}_d \left( \bfx; \bftheta \right) = - \textsf{NN}_d \left( \bfx; \left( \bfW_1, \cdots, \bfW_d \right) \right) = \textsf{NN}_d \left( \bfx; \left( \bfW_1, \cdots, - \bfW_d \right) \right), \quad \forall \bfx \in \bbX,
    \end{split}
\]
together with $\left( \bfW_1, \cdots, - \bfW_d \right) \in \Theta \left( M_{\textsf{F}} \right)$. This observation implies that
\begin{equation}
    \label{eqn:proof_prop:rademacher_complexity_neural_networks_v1}
    \begin{split}
        \calH_d \left( \bbX; M_{\textsf{F}} \right) = - \calH_d \left( \bbX; M_{\textsf{F}} \right) = \left\{ - \textsf{NN}_d \left( \cdot; \bftheta \right) \in \left( \bbX \to \bbR \right) : \bftheta \in \Theta \left( M_{\textsf{F}} \right) \right\}.
    \end{split}
\end{equation}
With the observation \eqref{eqn:proof_prop:rademacher_complexity_neural_networks_v1} in hand, one can realize from \emph{Theorem 1} in \cite{golowich2018size} that
\begin{equation}
    \label{eqn:proof_prop:rademacher_complexity_neural_networks_v2}
    \begin{split}
        &\calR_{n}^{\mu} \left( \calH_d \left( \bbX; M_{\textsf{F}} \right) \right) \\
        = \ & \bbE_{\left( \bfX_{1:n}, \bfsigma_{1:n} \right) \sim \mu^{\otimes n} \otimes \textsf{Unif} \left( \left\{ \pm 1 \right\}^n \right)} \left[ \sup \left\{ \left| \frac{1}{n} \sum_{i=1}^{n} \sigma_i \textsf{NN}_d \left( X_i; \bftheta \right) \right|: \bftheta \in \Theta \left( M_{\textnormal{F}} \right) \right\} \right] \\
        \leq \ & \frac{R \left( 1 + \sqrt{\left( 2 \log 2 \right) d} \right) \prod_{j=1}^{d} M_{\textsf{F}} (j)}{\sqrt{n}}.
    \end{split}
\end{equation}
On the other hand, by virtue of the Ledoux-Talagrand contraction principle (Lemma \ref{lemma:ledoux_talagrand_contraction}), we obtain that
\begin{equation}
    \label{eqn:proof_prop:rademacher_complexity_neural_networks_v3}
    \begin{split}
        &\calR_{n}^{\mu} (\calF) \\
        = \ & \bbE_{\bfX_{1:n} \sim \mu^{\otimes n}} \left[ \widehat{\calR}_n (\calF) \left( \bfX_{1:n} \right) \right] \\
        = \ & \bbE_{\bfX_{1:n} \sim \mu^{\otimes n}} \left[ \bbE_{\bfsigma_{1:n} \sim \textsf{Unif} \left( \left\{ \pm 1 \right\}^n \right)} \left[ \sup \left\{ \left| \frac{1}{n} \sum_{i=1}^{n} \eta \left\{ \textsf{NN}_d \left( X_i ; \bftheta \right) \right\} \right|: \bftheta \in \Theta \left( M_{\textsf{F}} \right)  \right\} \right] \right] \\
        \leq \ & 2 L \cdot \bbE_{\bfX_{1:n} \sim \mu^{\otimes n}} \left[ \bbE_{\bfsigma_{1:n} \sim \textsf{Unif} \left( \left\{ \pm 1 \right\}^n \right)} \left[ \sup \left\{ \left| \frac{1}{n} \sum_{i=1}^{n} \textsf{NN}_d \left( X_i ; \bftheta \right) \right|: \bftheta \in \Theta \left( M_{\textsf{F}} \right)  \right\} \right] \right] \\
        = \ & 2 L \cdot \bbE_{\bfX_{1:n} \sim \mu^{\otimes n}} \left[ \widehat{\calR}_n \left( \calH_d \left( \bbX; M_{\textsf{F}} \right) \right) \left( \bfX_{1:n} \right) \right] \\
        = \ & 2 L \cdot \calR_{n}^{\mu} \left( \calH_d \left( \bbX; M_{\textsf{F}} \right) \right) \\
        \stackrel{\textnormal{(a)}}{\leq} \ & \frac{2 L R \left( 1 + \sqrt{\left( 2 \log 2 \right) d} \right) \prod_{j=1}^{d} M_{\textsf{F}} (j)}{\sqrt{n}},
    \end{split}
\end{equation}
where the step (a) follows from the upper bound \eqref{eqn:proof_prop:rademacher_complexity_neural_networks_v2} on the Rademacher complexity of $\calH_d \left( \bbX; M_{\textsf{F}} \right)$. Hence, one can reveal that
\[
    \begin{split}
        \calR_{n}^{\mu} \left( \calF^* \right) \stackrel{\textnormal{(b)}}{\leq} \ & \calR_{n}^{\mu} (\calF) + \calR_{n}^{\mu} \left( \left\{ f^* \right\} \right) \\
        \stackrel{\textnormal{(c)}}{\leq} \ & \frac{2 L R \left( 1 + \sqrt{\left( 2 \log 2 \right) d} \right) \prod_{j=1}^{d} M_{\textnormal{\textsf{F}}} (j)}{\sqrt{n}} + \frac{2 \sqrt{\log 2}}{\sqrt{n}},
    \end{split}
\]
which thus completes the proof of Proposition \ref{prop:rademacher_complexity_neural_networks}, where the step (b) follows by the triangle inequality and the step (c) invokes the bound \eqref{eqn:proof_prop:rademacher_complexity_neural_networks_v3} and the standard upper bound on the Rademacher complexity of finite hypothesis classes (see Lemma \ref{lemma:rademacher_complexity_finite_class} for details).

\subsection{Proof of Theorem \ref{thm:upper_bound_alg:dr_covariate_shift_adaptation_v1}}
\label{subsec:proof_thm:upper_bound_alg:dr_covariate_shift_adaptation_v1}

\indent We first provide a formal definition of the \emph{doubly-robust (DR) empirical risk} $\widehat{\calR}_{\textsf{DR}} : \bbO^{n_{\bbP}} \times \bbX^{n_{\bbQ}} \to \left( \calF \to \bbR \right)$, where
\begin{equation}
    \label{eqn:formal_def_dr_empirical_risk_v1}
    \begin{split}
        \widehat{\calR}_{\textsf{DR}} \left( \bfo_{1:n_{\bbP}}, \bfx_{1:n_{\bbQ}} \right)(f) := \ & \frac{1}{n_{\bbP}} \sum_{i=1}^{n_{\bbP}} \hat{\rho} \left( x_{i}^{\bbP} \right) \left[ \left\{ y_{i}^{\bbP} - f \left( x_{i}^{\bbP} \right) \right\}^2 - \left\{ \hat{f}_0 \left( x_{i}^{\bbP} \right) - f \left( x_{i}^{\bbP} \right) \right\}^2 \right] \\
        &+ \frac{1}{n_{\bbQ}} \sum_{j=1}^{n_{\bbQ}} \left\{ \hat{f}_0 \left( X_{j}^{\bbQ} \right) - f \left( X_{j}^{\bbQ} \right) \right\}^2,
    \end{split}
\end{equation}
and define the \emph{DR population risk} $\overline{\calR}: \calF \to \bbR$ by
\begin{equation}
    \label{eqn:subsec:proof_thm:upper_bound_alg:dr_covariate_shift_adaptation_v1_v1}
    \begin{split}
        \overline{\calR} (f) := \bbE_{\left( \bfO_{1:n_{\bbP}}^{\bbP}, \bfX_{1:n_{\bbQ}}^{\bbQ} \right) \sim \bbP^{\otimes n_{\bbP}} \otimes \bbQ_{X}^{\otimes n_{\bbQ}}} \left[ \widehat{\calR}_{\textsf{DR}} \left( \bfO_{1:n_{\bbP}}^{\bbP}, \bfX_{1:n_{\bbQ}}^{\bbQ} \right) (f) \right], \quad f \in \calF.
    \end{split}
\end{equation}
Let us note here that $\widehat{\calR}_{\textsf{DR}} = \widehat{\calR}_{\textsf{DR}} \left( \bfO_{1:n_{\bbP}}^{\bbP}, \bfX_{1:n_{\bbQ}}^{\bbQ} \right) : \calF \to \bbR$ under $\left( \bfO_{1:n_{\bbP}}^{\bbP}, \bfX_{1:n_{\bbQ}}^{\bbQ} \right) \sim \bbP^{\otimes n_{\bbP}} \otimes \bbQ_{X}^{\otimes n_{\bbQ}}$. Then, one can decompose the DR population risk $\overline{\calR}: \calF \to \bbR$ as follows:
\begin{equation}
    \label{eqn:subsec:proof_thm:upper_bound_alg:dr_covariate_shift_adaptation_v1_v2}
    \begin{split}
        &\overline{\calR} (f) \\
        = \ & \bbE_{(X, Y) \sim \bbP} \left[ \hat{\rho} (X) \left\{ Y - f (X) \right\}^2 \right] - \bbE_{X \sim \bbP_{X}} \left[ \hat{\rho} (X) \left\{ \hat{f}_0 (X) - f (X) \right\}^2 \right] \\
        &+ \bbE_{X \sim \bbQ_{X}} \left[ \left\{ \hat{f}_0 (X) - f (X) \right\}^2 \right] \\
        = \ & \bbE_{(X, Y) \sim \bbP} \left[ \rho^* (X) \left\{ Y - f (X) \right\}^2 \right] \\
        &+ \bbE_{(X, Y) \sim \bbP} \left[ \left\{ \hat{\rho} (X) - \rho^* (X) \right\} \left[ \left\{ Y - f (X) \right\}^2 - \left\{ \hat{f}_0 (X) - f (X) \right\}^2 \right] \right] \\
        \stackrel{\textnormal{(a)}}{=} \ & \calR_{\bbQ} (f) + \bbE_{(X, Y) \sim \bbP} \left[ \left\{ \hat{\rho} (X) - \rho^* (X) \right\} \left[ \left\{ Y - f (X) \right\}^2 - \left\{ \hat{f}_0 (X) - f (X) \right\}^2 \right] \right], 
    \end{split}
\end{equation}
where the step (a) follows due to the observation \eqref{eqn:importance_weighting}. The definition of the DR estimator \eqref{eqn:dr_estimator_v1} yields the following \emph{basic inequality}: $0 \leq \widehat{\calR}_{\textsf{DR}} (f) - \widehat{\calR}_{\textsf{DR}} \left( \hat{f}_{\textsf{DR}} \right)$ for every $f \in \calF$. Thus, we have
\begin{equation}
    \label{eqn:subsec:proof_thm:upper_bound_alg:dr_covariate_shift_adaptation_v1_v3}
    \begin{split}
        0 \stackrel{\textnormal{(b)}}{\leq} \ & \widehat{\calR}_{\textsf{DR}} \left( f^* \right) - \widehat{\calR}_{\textsf{DR}} \left( \hat{f}_{\textsf{DR}} \right) \\
        = \ & \left\{ \widehat{\calR}_{\textsf{DR}} \left( f^* \right) - \overline{\calR} \left( f^* \right) \right\} + \left\{ \overline{\calR} \left( f^* \right) - \calR_{\bbQ} \left( f^* \right) \right\} - \underbrace{\left\{ \calR_{\bbQ} \left( \hat{f}_{\textsf{DR}} \right) - \calR_{\bbQ} \left( f^* \right) \right\}}_{= \ \calE_{\bbQ} \left( \hat{f}_{\textsf{DR}} \right)} \\
        &- \left\{ \overline{\calR} \left( \hat{f}_{\textsf{DR}} \right) - \calR_{\bbQ} \left( \hat{f}_{\textsf{DR}} \right) \right\} - \left\{ \widehat{\calR}_{\textsf{DR}} \left( \hat{f}_{\textsf{DR}} \right) - \overline{\calR} \left( \hat{f}_{\textsf{DR}} \right) \right\},
    \end{split}
\end{equation}
where the step (b) holds by the well-specification assumption of the model $f^* \in \calF$. It follows that
\begin{equation}
    \label{eqn:subsec:proof_thm:upper_bound_alg:dr_covariate_shift_adaptation_v1_v4}
    \begin{split}
        \calE_{\bbQ} \left( \hat{f}_{\textsf{DR}} \right) = \ & \calR_{\bbQ} \left( \hat{f}_{\textsf{DR}} \right) - \calR_{\bbQ} \left( f^* \right) \\
        \leq \ & \underbrace{\left\{ \overline{\calR} \left( f^* \right) - \calR_{\bbQ} \left( f^* \right) \right\} - \left\{ \overline{\calR} \left( \hat{f}_{\textsf{DR}} \right) - \calR_{\bbQ} \left( \hat{f}_{\textsf{DR}} \right) \right\}}_{=: \ \textnormal{(T1)}} \\
        &+ \underbrace{\left\{ \widehat{\calR}_{\textsf{DR}} \left( f^* \right) - \widehat{\calR}_{\textsf{DR}} \left( \hat{f}_{\textsf{DR}} \right) \right\} - \left\{ \overline{\calR} \left( f^* \right) - \overline{\calR} \left( \hat{f}_{\textsf{DR}} \right) \right\}}_{=: \ \textnormal{(T2)}}.
    \end{split}
\end{equation}

\paragraph{Bounding the term (T1):}
With the decomposition \eqref{eqn:subsec:proof_thm:upper_bound_alg:dr_covariate_shift_adaptation_v1_v4} in hand, let us first take a closer inspection on the first term (T1).
\[
    \begin{split}
        \textnormal{(T1)} \stackrel{\textnormal{(c)}}{=} \ & \bbE_{(X, Y) \sim \bbP} \left[ \left\{ \hat{\rho} (X) - \rho^* (X) \right\} \left[ \left\{ Y - f^* (X) \right\}^2 - \left\{ \hat{f}_0 (X) - f^* (X) \right\}^2 \right. \right. \\
        &\left. \left. - \left\{ Y - \hat{f}_{\textsf{DR}} (X) \right\}^2 + \left\{ \hat{f}_0 (X) - \hat{f}_{\textsf{DR}} (X) \right\}^2 \right] \right] \\
        = \ & 2 \bbE_{(X, Y) \sim \bbP} \left[ \left\{ \hat{\rho} (X) - \rho^* (X) \right\} \left\{ Y - \hat{f}_0 (X) \right\} \left\{ \hat{f}_{\textsf{DR}} (X) - f^* (X) \right\} \right] \\
        = \ & 2 \bbE_{X \sim \bbP_{X}} \left[ \left\{ \hat{\rho} (X) - \rho^* (X) \right\} \left\{ f^* (X) - \hat{f}_0 (X) \right\} \left\{ \hat{f}_{\textsf{DR}} (X) - f^* (X) \right\} \right],
    \end{split}
\]
where the step (c) uses the decomposition \eqref{eqn:subsec:proof_thm:upper_bound_alg:dr_covariate_shift_adaptation_v1_v2} of the DR population risk $\overline{\calR}: \calF \to \bbR$. Therefore, we have
\begin{equation}
    \label{eqn:subsec:proof_thm:upper_bound_alg:dr_covariate_shift_adaptation_v1_v5}
    \begin{split}
        \textnormal{(T1)} \leq \ & 2 \left| \bbE_{X \sim \bbP_{X}} \left[ \left\{ \hat{\rho} (X) - \rho^* (X) \right\} \left\{ f^* (X) - \hat{f}_0 (X) \right\} \left\{ \hat{f}_{\textsf{DR}} (X) - f^* (X) \right\} \right] \right| \\
        \leq \ & 2 \bbE_{X \sim \bbP_{X}} \left[ \left| \hat{\rho} (X) - \rho^* (X) \right| \cdot \left| f^* (X) - \hat{f}_0 (X) \right| \cdot \left| \hat{f}_{\textsf{DR}} (X) - f^* (X) \right| \right] \\
        \stackrel{\textnormal{(d)}}{\leq} \ & 4 \left\{ \bbE_{X \sim \bbP_X} \left[ \left\{ \hat{\rho} (X) - \rho^* (X) \right\}^2 \right] \right\}^{\frac{1}{2}} \left\{ \bbE_{X \sim \bbP_X} \left[ \left\{ \hat{f}_0 (X) - f^* (X) \right\}^2 \right] \right\}^{\frac{1}{2}} \\
        = \ & 4 \left\| \hat{\rho} - \rho^* \right\|_{L^2 \left( \bbX, \bbP_X \right)} \cdot \left\| \hat{f}_0 - f^* \right\|_{L^2 \left( \bbX, \bbP_X \right)},
    \end{split}
\end{equation}
where the step (d) holds due to the Cauchy-Schwarz inequality together with the fact that $\left| \hat{f}_{\textsf{DR}} (x) - f^* (x) \right| \leq \left\| \hat{f}_{\textsf{DR}} \right\|_{\infty} + \left\| f^* \right\|_{\infty} \leq 2$, follows from Assumption \ref{assumption:uniform_boundedness}.
\medskip

\paragraph{Bounding the term (T2):} With regards to the term (T2), we utilize tools from the empirical processes theory in order to establish its upper bound. First, we observe for any $f \in \calF$ that
\begin{equation}
    \label{eqn:subsec:proof_thm:upper_bound_alg:dr_covariate_shift_adaptation_v1_v6}
    \begin{split}
        &\widehat{\calR}_{\textsf{DR}} \left( f^* \right) - \widehat{\calR}_{\textsf{DR}} (f) \\
        = \ & \frac{2}{n_{\bbP}} \sum_{i=1}^{n_{\bbP}} \hat{\rho} \left( X_{i}^{\bbP} \right) \left\{ f \left( X_{i}^{\bbP} \right) - f^* \left( X_{i}^{\bbP} \right) \right\} \left\{ Y_{i}^{\bbP} - \hat{f}_0 \left( X_{i}^{\bbP} \right) \right\} \\
        &+ \frac{2}{n_{\bbQ}} \sum_{j=1}^{n_{\bbQ}} \left\{ f \left( X_{j}^{\bbQ} \right) - f^* \left( X_{i}^{\bbQ} \right) \right\} \left\{ \hat{f}_0 \left( X_{j}^{\bbQ} \right) - f^* \left( X_{i}^{\bbQ} \right) \right\} \\
        &- \frac{1}{n_{\bbQ}} \sum_{j=1}^{n_{\bbQ}} \left\{ f \left( X_{j}^{\bbQ} \right) - f^* \left( X_{i}^{\bbQ} \right) \right\}^2.
    \end{split}
\end{equation}
With this observation in hand, it is seen that
\begin{equation}
    \label{eqn:subsec:proof_thm:upper_bound_alg:dr_covariate_shift_adaptation_v1_v7}
    \begin{split}
        \textnormal{(T2)} \leq \ & \sup \left\{ \left| \left\{ \widehat{\calR}_{\textsf{DR}} \left( f^* \right) - \widehat{\calR}_{\textsf{DR}} (f) \right\} - \left\{ \overline{\calR} \left( f^* \right) - \overline{\calR} (f) \right\} \right|: f \in \calF \right\} \\
        \stackrel{\textnormal{(e)}}{\leq} \ & 2 \sup \left\{ \left| \bbG_{n_{\bbP}}^{\bbP} (\varphi) \left( \bfO_{1:n_{\bbP}}^{\bbP} \right) \right|: \varphi \in \calF^* \right\} + 2 \sup \left\{ \left| \bbG_{n_{\bbQ}}^{\bbQ, (1)} (\varphi) \left( \bfX_{1:n_{\bbQ}}^{\bbQ} \right) \right|: \varphi \in \calF^* \right\} \\
        &+ \sup \left\{ \left| \bbG_{n_{\bbQ}}^{\bbQ, (2)} (\varphi) \left( \bfX_{1:n_{\bbQ}}^{\bbQ} \right) \right|: \varphi \in \calF^* \right\} \\
        = \ & 2 \sup \left\{ \bbG_{n_{\bbP}}^{\bbP} (\varphi) \left( \bfO_{1:n_{\bbP}}^{\bbP} \right) : \varphi \in \calF^* \cup \left( - \calF^* \right) \right\} + 2 \sup \left\{ \bbG_{n_{\bbQ}}^{\bbQ, (1)} (\varphi) \left( \bfX_{1:n_{\bbQ}}^{\bbQ} \right) : \varphi \in \calF^* \cup \left( - \calF^* \right) \right\} \\
        &+ \sup \left\{ \bbG_{n_{\bbQ}}^{\bbQ, (2)} (\varphi) \left( \bfX_{1:n_{\bbQ}}^{\bbQ} \right) : \varphi \in \calF^* \cup \left( - \calF^* \right) \right\},
    \end{split}
\end{equation}
where the step (e) follows by virtue of the triangle inequality. Here, $- \calF^* := \left\{ - \varphi : \varphi \in \calF^* \right\}$, and the functions $\bbG_{n_{\bbP}}^{\bbP} : \left( \bbX \to \bbR \right) \to \left( \bbO^{n_{\bbP}} \to \bbR \right)$ and $\left\{ \bbG_{n_{\bbQ}}^{\bbQ, (l)} : \left( \bbX \to \bbR \right) \to \left( \bbX^{n_{\bbQ}} \to \bbR \right) : l \in [2] \right\}$ are defined as
\begin{equation}
    \label{eqn:subsec:proof_thm:upper_bound_alg:dr_covariate_shift_adaptation_v1_v8}
    \begin{split}
        \bbG_{n_{\bbP}}^{\bbP} (\varphi) \left( \bfo_{1:n_{\bbP}}^{\bbP} \right) := \ & \frac{1}{n_{\bbP}} \sum_{i=1}^{n_{\bbP}} \hat{\rho} \left( x_{i}^{n_{\bbP}} \right) \varphi \left( x_{i}^{n_{\bbP}} \right) \left\{ y_{i}^{\bbP} - \hat{f}_0 \left( x_{i}^{n_{\bbP}} \right) \right\} - \bbE_{(X, Y) \sim \bbP} \left[ \hat{\rho} (X) \varphi (X) \left\{ Y - \hat{f}_0 (X) \right\} \right], \\ 
        \bbG_{n_{\bbQ}}^{\bbQ, (1)} (\varphi) \left( \bfx_{1:n_{\bbQ}}^{\bbQ} \right) := \ & \frac{1}{n_{\bbQ}} \sum_{j=1}^{n_{\bbQ}} \varphi \left( x_{j}^{\bbQ} \right) \left\{ \hat{f}_0 \left( x_{j}^{\bbQ} \right) - f^* \left( x_{j}^{\bbQ} \right) \right\} - \bbE_{X \sim \bbQ_{X}} \left[ \varphi (X) \left\{ \hat{f}_0 (X) - f^* (X) \right\} \right], \\
        \bbG_{n_{\bbQ}}^{\bbQ, (2)} (\varphi) \left( \bfx_{1:n_{\bbQ}}^{\bbQ} \right) := \ & \frac{1}{n_{\bbQ}} \sum_{j=1}^{n_{\bbQ}} \left\{ \varphi \left( x_{j}^{\bbQ} \right) \right\}^2 - \bbE_{X \sim \bbQ_{X}} \left[ \left\{ \varphi (X) \right\}^2 \right].
    \end{split}
\end{equation}
If $\bfO_{1:n_{\bbP}}^{\bbP} \sim \bbP^{\otimes n_{\bbP}}$ and $\bfX_{1:n_{\bbQ}}^{\bbQ} \sim \bbQ_{X}^{\otimes n_{\bbQ}}$, then
\begin{itemize}
    \item $\left\{ \bbG_{n_{\bbP}}^{\bbP} (\varphi) \left( \bfO_{1:n_{\bbP}}^{\bbP} \right) = \left( \widehat{\bbP} - \bbP \right) \left[ \hat{\rho} (X) \varphi(X) \left\{ Y - \hat{f}_0 (X) \right\} \right]: \varphi \in \calF^* \cup \left( - \calF^* \right) \right\}$,
    \item $\left\{ \bbG_{n_{\bbQ}}^{\bbQ, (1)} (f) \left( \bfX_{1:n_{\bbQ}}^{\bbQ} \right) = \left\{ \widehat{\bbQ}_{X} - \bbQ_{X} \right\} \left[ \varphi (X) \left\{ \hat{f}_0 (X) - f^* (X) \right\} \right]: \varphi \in \calF^* \cup \left( - \calF^* \right) \right\}$,
    \item $\left\{ \bbG_{n_{\bbQ}}^{\bbQ, (2)} (f) \left( \bfX_{1:n_{\bbQ}}^{\bbQ} \right) = \left\{ \widehat{\bbQ}_{X} - \bbQ_{X} \right\} \left[ \left\{ \varphi (X) \right\}^2 \right]: \varphi \in \calF^* \cup \left( - \calF^* \right) \right\}$
\end{itemize}
are empirical processes indexed by $\varphi \in \calF^* \cup \left( - \calF^* \right)$, where $\widehat{\bbP} \in \Delta \left( \bbX \times \bbR \right)$ and $\widehat{\bbQ}_{X} \in \Delta (\bbX)$ are the empirical distributions for the $n_{\bbP}$ labeled source samples $\bfO_{1:n_{\bbP}}^{\bbP}$ and $n_{\bbQ}$ unlabeled target samples $\bfX_{1:n_{\bbQ}}^{\bbQ}$, respectively, i.e., $\widehat{\bbP} := \frac{1}{n_{\bbP}} \sum_{i=1}^{n_{\bbP}} \delta_{\left( X_{i}^{\bbP}, Y_{i}^{\bbP} \right)}$ and $\widehat{\bbQ}_{X} := \frac{1}{n_{\bbQ}} \sum_{j=1}^{n_{\bbQ}} \delta_{X_{j}^{\bbQ}}$.

\paragraph{Control of the supremum of $\left\{ \bbG_{n_{\bbP}}^{\bbP} (\varphi) \left( \bfO_{1:n_{\bbP}}^{\bbP} \right) : \varphi \in \calF^* \cup \left( - \calF^* \right) \right\}$:}
Firstly, we are in need of a delicate control of the expectation of the supremum of the empirical process
\[
    \left\{ \bbG_{n_{\bbP}}^{\bbP} (\varphi) \left( \bfO_{1:n_{\bbP}}^{\bbP} \right) : \varphi \in \calF^* \cup \left( - \calF^* \right) \right\}.
\]
This goal can be achieved through the following lemma, whose proof is provided in Appendix \ref{subsec:proof_lemma:control_P_empirical_process_v1}.

\begin{lemma}
\label{lemma:control_P_empirical_process_v1}
The expectation of the supremum of the empirical process $\left\{ \bbG_{n_{\bbP}}^{\bbP} (\varphi) \left( \bfO_{1:n_{\bbP}}^{\bbP} \right) : \varphi \in \calF^* \cup \left( - \calF^* \right) \right\}$ is upper bounded by
\begin{equation}
    \label{eqn:lemma:control_P_empirical_process_v1_v1}
    \begin{split}
        \bbE_{\bfO_{1:n_{\bbP}}^{\bbP} \sim \bbP^{\otimes n_{\bbP}}} \left[ \sup \left\{ \bbG_{n_{\bbP}}^{\bbP} (\varphi) \left( \bfO_{1:n_{\bbP}}^{\bbP} \right) : \varphi \in \calF^* \cup \left( - \calF^* \right) \right\} \right] \leq 4 C_{\textnormal{\textsf{dr}}} \left( 1 + C_{\textnormal{\textsf{rf}}} \right) \calR_{n_{\bbP}}^{\bbP_{X}} \left( \calF^* \right).
    \end{split}
\end{equation}
\end{lemma}

\noindent We then move on to a tight control of the size of 
\[
    \begin{split}
        \left\{ \bbG_{n_{\bbP}}^{\bbP} (\varphi) \left( \bfO_{1:n_{\bbP}}^{\bbP} \right) : \varphi \in \calF^* \cup \left( - \calF^* \right) \right\} - \bbE_{\bfO_{1:n_{\bbP}}^{\bbP} \sim \bbP^{\otimes n_{\bbP}}} \left[ \left\{ \bbG_{n_{\bbP}}^{\bbP} (\varphi) \left( \bfO_{1:n_{\bbP}}^{\bbP} \right) : \varphi \in \calF^* \cup \left( - \calF^* \right) \right\} \right]
    \end{split}
\]
under $\bfO_{1:n_{\bbP}}^{\bbP} \sim \bbP^{\otimes n_{\bbP}}$. This task can be settled via the following lemma, whose proof is deferred to Appendix \ref{subsec:proof_lemma:control_P_empirical_process_v2}.

\begin{lemma}
\label{lemma:control_P_empirical_process_v2}
If $\bfO_{1:n_{\bbP}}^{\bbP} \sim \bbP^{\otimes n_{\bbP}}$, then with probability at least $1 - \delta$, we have
\begin{equation}
    \label{eqn:lemma:control_P_empirical_process_v2_v1}
    \begin{split}
        &\sup \left\{ \bbG_{n_{\bbP}}^{\bbP} (\varphi) \left( \bfO_{1:n_{\bbP}}^{\bbP} \right) : \varphi \in \calF^* \cup \left( - \calF^* \right) \right\} \\
        &- \bbE_{\bfO_{1:n_{\bbP}}^{\bbP} \sim \bbP^{\otimes n_{\bbP}}} \left[ \sup \left\{ \bbG_{n_{\bbP}}^{\bbP} (\varphi) \left( \bfO_{1:n_{\bbP}}^{\bbP} \right) : \varphi \in \calF^* \cup \left( - \calF^* \right) \right\} \right] \\
        \leq \ & \frac{6 C_{\textnormal{\textsf{dr}}} \left( 1 + C_{\textnormal{\textsf{rf}}} \right)}{n_{\bbP}} \log \left( \frac{1}{\delta} \right) + 2 C_{\textnormal{\textsf{dr}}} \left( 1 + C_{\textnormal{\textsf{rf}}} \right) \sqrt{\frac{2 \log \left( \frac{1}{\delta} \right)}{n_{\bbP}}} + 4 C_{\textnormal{\textsf{dr}}} \left( 1 + C_{\textnormal{\textsf{rf}}} \right) \sqrt{ \frac{2 \calR_{n_{\bbP}}^{\bbP_{X}} \left( \calF^* \right) \log \left( \frac{1}{\delta} \right)}{n_{\bbP}}}.
    \end{split}
\end{equation}
\end{lemma}

\indent To finish up, we first denote the right-hand side of the inequality \eqref{eqn:lemma:control_P_empirical_process_v2_v1} from Lemma \ref{lemma:control_P_empirical_process_v2} by
\begin{equation}
    \label{eqn:subsec:proof_thm:upper_bound_alg:dr_covariate_shift_adaptation_v1_v9}
    \begin{split}
        \calB_{\bbP} (\delta) := \ & \frac{6 C_{\textnormal{\textsf{dr}}} \left( 1 + C_{\textnormal{\textsf{rf}}} \right)}{n_{\bbP}} \log \left( \frac{1}{\delta} \right) + 2 C_{\textnormal{\textsf{dr}}} \left( 1 + C_{\textnormal{\textsf{rf}}} \right) \sqrt{\frac{2 \log \left( \frac{1}{\delta} \right)}{n_{\bbP}}} \\
        &+ 4 C_{\textnormal{\textsf{dr}}} \left( 1 + C_{\textnormal{\textsf{rf}}} \right) \sqrt{ \frac{2 \calR_{n_{\bbP}}^{\bbP_{X}} \left( \calF^* \right) \log \left( \frac{1}{\delta} \right)}{n_{\bbP}}}.
    \end{split}
\end{equation}
for ease of exposition. Then, with probability at least $1 - \delta$, one has
\begin{equation}
    \label{eqn:subsec:proof_thm:upper_bound_alg:dr_covariate_shift_adaptation_v1_v10}
    \begin{split}
        &\sup \left\{ \bbG_{n_{\bbP}}^{\bbP} (\varphi) \left( \bfO_{1:n_{\bbP}}^{\bbP} \right) : \varphi \in \calF^* \cup \left( - \calF^* \right) \right\} \\
        = \ & \sup \left\{ \bbG_{n_{\bbP}}^{\bbP} (\varphi) \left( \bfO_{1:n_{\bbP}}^{\bbP} \right) : \varphi \in \calF^* \cup \left( - \calF^* \right) \right\} \\
        &- \bbE_{\bfO_{1:n_{\bbP}}^{\bbP} \sim \bbP^{\otimes n_{\bbP}}} \left[ \sup \left\{ \bbG_{n_{\bbP}}^{\bbP} (\varphi) \left( \bfO_{1:n_{\bbP}}^{\bbP} \right) : \varphi \in \calF^* \cup \left( - \calF^* \right) \right\} \right] \\
        &+ \bbE_{\bfO_{1:n_{\bbP}}^{\bbP} \sim \bbP^{\otimes n_{\bbP}}} \left[ \sup \left\{ \bbG_{n_{\bbP}}^{\bbP} (\varphi) \left( \bfO_{1:n_{\bbP}}^{\bbP} \right) : \varphi \in \calF^* \cup \left( - \calF^* \right) \right\} \right] \\
        \stackrel{\textnormal{(f)}}{\leq} \ & \calB_{\bbP} (\delta) + 4 C_{\textnormal{dr}} \left( 1 + C_{\textnormal{rf}} \right) \calR_{n_{\bbP}}^{\bbP_{X}} \left( \calF^* \right), 
    \end{split}
\end{equation}
where the step (f) invokes Lemmas \ref{lemma:control_P_empirical_process_v1} and \ref{lemma:control_P_empirical_process_v2}. For simplicity, we define the following event: for $\delta \in (0, 1)$,
\begin{equation}
    \label{eqn:subsec:proof_thm:upper_bound_alg:dr_covariate_shift_adaptation_v1_v11}
    \begin{split}
        &\calE_{\bbP} (\delta) := \left\{ \left( \bfo_{1:n_{\bbP}}^{\bbP}, \bfx_{1:n_{\bbQ}}^{\bbQ} \right) \in \bbO^{n_{\bbP}} \times \bbX^{n_{\bbQ}}: \right. \\
        &\left. \sup \left\{ \bbG_{n_{\bbP}}^{\bbP} (\varphi) \left( \bfo_{1:n_{\bbP}}^{\bbP} \right) : \varphi \in \calF^* \cup \left( - \calF^* \right) \right\} \leq \calB_{\bbP} (\delta) + 4 C_{\textnormal{\textsf{dr}}} \left( 1 + C_{\textnormal{\textsf{rf}}} \right) \calR_{n_{\bbP}}^{\bbP_{X}} \left( \calF^* \right) \right\}.
    \end{split}
\end{equation}
Then, the upper bound \eqref{eqn:subsec:proof_thm:upper_bound_alg:dr_covariate_shift_adaptation_v1_v10} implies
\begin{equation}
    \label{eqn:subsec:proof_thm:upper_bound_alg:dr_covariate_shift_adaptation_v1_v12}
    \begin{split}
        \left( \bbP^{\otimes n_{\bbP}} \otimes \bbQ_{X}^{\otimes n_{\bbQ}} \right) \left\{ \calE_{\bbP} (\delta) \right\} \geq 1 - \delta \quad \textnormal{for every} \quad \delta \in (0, 1).
    \end{split}
\end{equation}

\paragraph{Control of the supremum of $\left\{ \bbG_{n_{\bbQ}}^{\bbQ, (1)} (\varphi) \left( \bfX_{1:n_{\bbQ}}^{\bbQ} \right) : \varphi \in \calF^* \cup \left( - \calF^* \right) \right\}$:}
Similar to the preceding argument for controlling the supremum of the empirical process $\left\{ \bbG_{n_{\bbP}}^{\bbP} (\varphi) \left( \bfO_{1:n_{\bbP}}^{\bbP} \right) : \varphi \in \calF^* \cup \left( - \calF^* \right) \right\}$, we first establish an upper bound on the expectation of the supremum of the empirical process
\[
    \left\{ \bbG_{n_{\bbQ}}^{\bbQ, (1)} (\varphi) \left( \bfX_{1:n_{\bbQ}}^{\bbQ} \right) : \varphi \in \calF^* \cup \left( - \calF^* \right) \right\}.
\]
We provide a desired result in the following lemma, whose proof is postponed to Appendix \ref{subsec:proof_lemma:control_Q_empirical_process_v1}.

\begin{lemma}
\label{lemma:control_Q_empirical_process_v1}
The expectation of the supremum of the empirical process $\left\{ \bbG_{n_{\bbQ}}^{\bbQ, (1)} (\varphi) \left( \bfX_{1:n_{\bbQ}}^{\bbQ} \right) : \varphi \in \calF^* \cup \left( - \calF^* \right) \right\}$ can be upper bounded by
\begin{equation}
    \label{eqn:lemma:control_Q_empirical_process_v1_v1}
    \begin{split}
        \bbE_{\bfX_{1:n_{\bbQ}}^{\bbQ} \sim \bbQ_{X}^{\otimes n_{\bbQ}}} \left[ \sup \left\{ \bbG_{n_{\bbQ}}^{\bbQ, (1)} (\varphi) \left( \bfX_{1:n_{\bbQ}}^{\bbQ} \right) : \varphi \in \calF^* \cup \left( - \calF^* \right) \right\} \right] \leq 4 \left( 1 + C_{\textnormal{\textsf{rf}}} \right) \calR_{n_{\bbQ}}^{\bbQ_{X}} \left( \calF^* \right).
    \end{split}
\end{equation}
\end{lemma}

\noindent Analogously, we now aim at a tight control of the size of
\[
    \begin{split}
        \sup \left\{ \bbG_{n_{\bbQ}}^{\bbQ, (1)} (\varphi) \left( \bfX_{1:n_{\bbQ}}^{\bbQ} \right) : \varphi \in \calF^* \cup \left( - \calF^* \right) \right\} - \bbE_{\bfX_{1:n_{\bbQ}}^{\bbQ} \sim \bbQ_{X}^{\otimes n_{\bbQ}}} \left[ \sup \left\{ \bbG_{n_{\bbQ}}^{\bbQ, (1)} (\varphi) \left( \bfX_{1:n_{\bbQ}}^{\bbQ} \right) : \varphi \in \calF^* \cup \left( - \calF^* \right) \right\} \right]
    \end{split}
\]
under the data generating process $\bfX_{1:n_{\bbQ}}^{\bbQ} \sim \bbQ_{X}^{\otimes n_{\bbQ}}$. This goal can be achieved through the following lemma, whose detailed proof is provided in Appendix \ref{subsec:proof_lemma:control_Q_empirical_process_v2}.

\begin{lemma}
\label{lemma:control_Q_empirical_process_v2}
If $\bfX_{1:n_{\bbQ}}^{\bbQ} \sim \bbQ_{X}^{\otimes n_{\bbQ}}$, then with probability at least $1 - \delta$, it holds that
\begin{equation}
    \label{eqn:lemma:control_Q_empirical_process_v2_v1}
    \begin{split}
        &\sup \left\{ \bbG_{n_{\bbQ}}^{\bbQ, (1)} (\varphi) \left( \bfX_{1:n_{\bbQ}}^{\bbQ} \right) : \varphi \in \calF^* \cup \left( - \calF^* \right) \right\} \\
        &- \bbE_{\bfX_{1:n_{\bbQ}}^{\bbQ} \sim \bbQ_{X}^{\otimes n_{\bbQ}}} \left[ \sup \left\{ \bbG_{n_{\bbQ}}^{\bbQ, (1)} (\varphi) \left( \bfX_{1:n_{\bbQ}}^{\bbQ} \right) : \varphi \in \calF^* \cup \left( - \calF^* \right) \right\} \right] \\
        \leq \ & \frac{6 \left( 1 + C_{\textnormal{\textsf{rf}}} \right)}{n_{\bbQ}} \log \left( \frac{1}{\delta} \right) + 2 \left( 1 + C_{\textnormal{\textsf{rf}}} \right) \sqrt{\frac{2 \log \left( \frac{1}{\delta} \right)}{n_{\bbQ}}} + 4 \left( 1 + C_{\textnormal{\textsf{rf}}} \right) \sqrt{\frac{2 \calR_{n_{\bbQ}}^{\bbQ_{X}} \left( \calF^* \right) \log \left( \frac{1}{\delta} \right)}{n_{\bbQ}}}.
    \end{split}
\end{equation}
\end{lemma}

\noindent For simplicity, we denote the right-hand side of the inequality \eqref{eqn:lemma:control_Q_empirical_process_v2_v1} from Lemma \ref{lemma:control_Q_empirical_process_v2} for any given $\delta \in (0, 1)$ as
\begin{equation}
    \label{eqn:subsec:proof_thm:upper_bound_alg:dr_covariate_shift_adaptation_v1_v13}
    \begin{split}
        \calB_{\bbQ}^{(1)} (\delta) := \frac{6 \left( 1 + C_{\textnormal{\textsf{rf}}} \right)}{n_{\bbQ}} \log \left( \frac{1}{\delta} \right) + 2 \left( 1 + C_{\textnormal{\textsf{rf}}} \right) \sqrt{\frac{2 \log \left( \frac{1}{\delta} \right)}{n_{\bbQ}}} + 4 \left( 1 + C_{\textnormal{\textsf{rf}}} \right) \sqrt{\frac{2 \calR_{n_{\bbQ}}^{\bbQ_{X}} \left( \calF^* \right) \log \left( \frac{1}{\delta} \right)}{n_{\bbQ}}}.
    \end{split}
\end{equation}
Then, it holds with probability greater than $1 - \delta$ that
\begin{equation}
    \label{eqn:subsec:proof_thm:upper_bound_alg:dr_covariate_shift_adaptation_v1_v14}
    \begin{split}
        &\sup \left\{ \bbG_{n_{\bbQ}}^{\bbQ, (1)} (\varphi) \left( \bfX_{1:n_{\bbQ}}^{\bbQ} \right) : \varphi \in \calF^* \cup \left( - \calF^* \right) \right\} \\
        = \ & \sup \left\{ \bbG_{n_{\bbQ}}^{\bbQ, (1)} (\varphi) \left( \bfX_{1:n_{\bbQ}}^{\bbQ} \right) : \varphi \in \calF^* \cup \left( - \calF^* \right) \right\} \\
        &- \bbE_{\bfX_{1:n_{\bbQ}}^{\bbQ} \sim \bbQ_{X}^{\otimes n_{\bbQ}}} \left[ \sup \left\{ \bbG_{n_{\bbQ}}^{\bbQ, (1)} (\varphi) \left( \bfX_{1:n_{\bbQ}}^{\bbQ} \right) : \varphi \in \calF^* \cup \left( - \calF^* \right) \right\} \right] \\
        &+ \bbE_{\bfX_{1:n_{\bbQ}}^{\bbQ} \sim \bbQ_{X}^{\otimes n_{\bbQ}}} \left[ \sup \left\{ \bbG_{n_{\bbQ}}^{\bbQ, (1)} (\varphi) \left( \bfX_{1:n_{\bbQ}}^{\bbQ} \right) : \varphi \in \calF^* \cup \left( - \calF^* \right) \right\} \right] \\
        \stackrel{\textnormal{(g)}}{\leq} \ & \calB_{\bbQ}^{(1)} (\delta) + 4 \left( 1 + C_{\textsf{rf}} \right) \calR_{n_{\bbQ}}^{\bbQ_{X}} \left( \calF^* \right), 
    \end{split}
\end{equation}
where the step (g) holds due to Lemma \ref{lemma:control_Q_empirical_process_v1} and \ref{lemma:control_Q_empirical_process_v2}. For brevity, let us define the following event: for any given $\delta \in (0, 1)$,
\begin{equation}
    \label{eqn:subsec:proof_thm:upper_bound_alg:dr_covariate_shift_adaptation_v1_v15}
    \begin{split}
        &\calE_{\bbQ}^{(1)} (\delta) := \left\{ \left( \bfo_{1:n_{\bbP}}^{\bbP}, \bfx_{1:n_{\bbQ}}^{\bbQ} \right) \in \bbO^{n_{\bbP}} \times \bbX^{n_{\bbQ}}: \right. \\
        &\left. \sup \left\{ \bbG_{n_{\bbQ}}^{\bbQ, (1)} (\varphi) \left( \bfx_{1:n_{\bbQ}}^{\bbQ} \right) : \varphi \in \calF^* \cup \left( - \calF^* \right) \right\} \leq \calB_{\bbQ}^{(1)} (\delta) + 4 \left( 1 + C_{\textsf{rf}} \right) \calR_{n_{\bbQ}}^{\bbQ_{X}} \left( \calF^* \right) \right\}.
    \end{split}
\end{equation}
Then, the upper bound \eqref{eqn:subsec:proof_thm:upper_bound_alg:dr_covariate_shift_adaptation_v1_v14} directly yields
\begin{equation}
    \label{eqn:subsec:proof_thm:upper_bound_alg:dr_covariate_shift_adaptation_v1_v16}
    \begin{split}
        \left( \bbP^{\otimes n_{\bbP}} \otimes \bbQ_{X}^{\otimes n_{\bbQ}} \right) \left\{ \calE_{\bbQ}^{(1)} (\delta) \right\} \geq 1 - \delta \quad \textnormal{for every} \quad \delta \in (0, 1).
    \end{split}
\end{equation}

\paragraph{Control of the supremum of $\left\{ \bbG_{n_{\bbQ}}^{\bbQ, (2)} (\varphi) \left( \bfX_{1:n_{\bbQ}}^{\bbQ} \right) : \varphi \in \calF^* \cup \left( - \calF^* \right) \right\}$:}
Akin to the above delicate control of the supremum of empirical processes, we bound the expectation of the supremum of the empirical process
\[
    \left\{ \bbG_{n_{\bbQ}}^{\bbQ, (2)} (\varphi) \left( \bfX_{1:n_{\bbQ}}^{\bbQ} \right) : \varphi \in \calF^* \cup \left( - \calF^* \right) \right\}.
\]
It can be developed via the following lemma, with the detailed proof postponed to Appendix \ref{subsec:proof_lemma:control_Q_empirical_process_v3}.

\begin{lemma}
\label{lemma:control_Q_empirical_process_v3}
The expectation of the supremum of the empirical process $\left\{ \bbG_{n_{\bbQ}}^{\bbQ, (2)} (\varphi) \left( \bfX_{1:n_{\bbQ}}^{\bbQ} \right) : \varphi \in \calF^* \cup \left( - \calF^* \right) \right\}$ has an upper bound
\begin{equation}
    \label{eqn:lemma:control_Q_empirical_process_v3_v1}
    \begin{split}
        \bbE_{\bfX_{1:n_{\bbQ}}^{\bbQ} \sim \bbQ_{X}^{\otimes n_{\bbQ}}} \left[ \sup \left\{ \bbG_{n_{\bbQ}}^{\bbQ, (2)} (\varphi) \left( \bfX_{1:n_{\bbQ}}^{\bbQ} \right) : \varphi \in \calF^* \cup \left( - \calF^* \right) \right\} \right] \leq 16 \cdot \calR_{n_{\bbQ}}^{\bbQ_{X}} \left( \calF^* \right).
    \end{split}
\end{equation}
\end{lemma}

\noindent As the next step, we now turn to a tight control of the size of
\[
    \begin{split}
        \sup \left\{ \bbG_{n_{\bbQ}}^{\bbQ, (2)} (\varphi) \left( \bfX_{1:n_{\bbQ}}^{\bbQ} \right) : \varphi \in \calF^* \cup \left( - \calF^* \right) \right\} - \bbE_{\bfX_{1:n_{\bbQ}}^{\bbQ} \sim \bbQ_{X}^{\otimes n_{\bbQ}}} \left[ \sup \left\{ \bbG_{n_{\bbQ}}^{\bbQ, (2)} (\varphi) \left( \bfX_{1:n_{\bbQ}}^{\bbQ} \right) : \varphi \in \calF^* \cup \left( - \calF^* \right) \right\} \right]
    \end{split}
\]
with $\bfX_{1:n_{\bbQ}}^{\bbQ} \sim \bbQ_{X}^{\otimes n_{\bbQ}}$. The following lemma takes a step forward towards this goal, whose proof is deferred to Appendix \ref{subsec:proof_lemma:control_Q_empirical_process_v4}.

\begin{lemma}
\label{lemma:control_Q_empirical_process_v4}
If $\bfX_{1:n_{\bbQ}}^{\bbQ} \sim \bbQ_{X}^{\otimes n_{\bbQ}}$, then with probability at least $1 - \delta$, it holds that
\begin{equation}
    \label{eqn:lemma:control_Q_empirical_process_v4_v1}
    \begin{split}
        &\sup \left\{ \bbG_{n_{\bbQ}}^{\bbQ, (2)} (\varphi) \left( \bfX_{1:n_{\bbQ}}^{\bbQ} \right) : \varphi \in \calF^* \cup \left( - \calF^* \right) \right\} \\
        &- \bbE_{\bfX_{1:n_{\bbQ}}^{\bbQ} \sim \bbQ_{X}^{\otimes n_{\bbQ}}} \left[ \sup \left\{ \bbG_{n_{\bbQ}}^{\bbQ, (2)} (\varphi) \left( \bfX_{1:n_{\bbQ}}^{\bbQ} \right) : \varphi \in \calF^* \cup \left( - \calF^* \right) \right\} \right] \\
        \leq \ & \frac{12}{n_{\bbQ}} \log \left( \frac{1}{\delta} \right) + 4 \sqrt{\frac{2 \log \left( \frac{1}{\delta} \right)}{n_{\bbQ}}} + 16 \sqrt{\frac{\calR_{n_{\bbQ}}^{\bbQ_{X}} \left( \calF^* \right) \log \left( \frac{1}{\delta} \right)}{n_{\bbQ}}}.
    \end{split}
\end{equation}
\end{lemma}

\noindent For ease of expression, we denote the right-hand side of the bound \eqref{eqn:lemma:control_Q_empirical_process_v4_v1} in Lemma \ref{lemma:control_Q_empirical_process_v4} by
\begin{equation}
    \label{eqn:subsec:proof_thm:upper_bound_alg:dr_covariate_shift_adaptation_v1_v17}
    \begin{split}
        \calB_{\bbQ}^{(2)} (\delta) := \frac{12}{n_{\bbQ}} \log \left( \frac{1}{\delta} \right) + 4 \sqrt{\frac{2 \log \left( \frac{1}{\delta} \right)}{n_{\bbQ}}} + 16 \sqrt{\frac{\calR_{n_{\bbQ}}^{\bbQ_{X}} \left( \calF^* \right) \log \left( \frac{1}{\delta} \right)}{n_{\bbQ}}}.
    \end{split}
\end{equation}

\noindent Then, one has with probability at least $1 - \delta$ that
\begin{equation}
    \label{eqn:subsec:proof_thm:upper_bound_alg:dr_covariate_shift_adaptation_v1_v18}
    \begin{split}
        &\sup \left\{ \bbG_{n_{\bbQ}}^{\bbQ, (2)} (\varphi) \left( \bfX_{1:n_{\bbQ}}^{\bbQ} \right) : \varphi \in \calF^* \cup \left( - \calF^* \right) \right\} \\
        = \ & \sup \left\{ \bbG_{n_{\bbQ}}^{\bbQ, (2)} (\varphi) \left( \bfX_{1:n_{\bbQ}}^{\bbQ} \right) : \varphi (\cdot) \in \calF^* \cup \left( - \calF^* \right) \right\} \\
        &- \bbE_{\bfX_{1:n_{\bbQ}}^{\bbQ} \sim \bbQ_{X}^{\otimes n_{\bbQ}}} \left[ \sup \left\{ \bbG_{n_{\bbQ}}^{\bbQ, (2)} (\varphi) \left( \bfX_{1:n_{\bbQ}}^{\bbQ} \right) : \varphi \in \calF^* \cup \left( - \calF^* \right) \right\} \right] \\
        &+ \bbE_{\bfX_{1:n_{\bbQ}}^{\bbQ} \sim \bbQ_{X}^{\otimes n_{\bbQ}}} \left[ \sup \left\{ \bbG_{n_{\bbQ}}^{\bbQ, (2)} (\varphi) \left( \bfX_{1:n_{\bbQ}}^{\bbQ} \right) : \varphi \in \calF^* \cup \left( - \calF^* \right) \right\} \right] \\
        \stackrel{\textnormal{(h)}}{\leq} \ & \calB_{\bbQ}^{(2)} (\delta) + 16 \cdot \calR_{n_{\bbQ}}^{\bbQ_{X}} \left( \calF^* \right), 
    \end{split}
\end{equation}
where the step (h) invokes Lemmas \ref{lemma:control_Q_empirical_process_v3} and \ref{lemma:control_Q_empirical_process_v4}. For the sake of simplicity, let us define the following event: for any $\delta \in (0, 1)$,
\begin{equation}
    \label{eqn:subsec:proof_thm:upper_bound_alg:dr_covariate_shift_adaptation_v1_v19}
    \begin{split}
        &\calE_{\bbQ}^{(2)} (\delta) := \left\{ \left( \bfo_{1:n_{\bbP}}^{\bbP}, \bfx_{1:n_{\bbQ}}^{\bbQ} \right) \in \bbO^{n_{\bbP}} \times \bbX^{n_{\bbQ}}: \right. \\
        &\left. \sup \left\{ \bbG_{n_{\bbQ}}^{\bbQ, (2)} (\varphi) \left( \bfx_{1:n_{\bbQ}}^{\bbQ} \right) : \varphi \in \calF^* \cup \left( - \calF^* \right) \right\} \leq \calB_{\bbQ}^{(2)} (\delta) + 16 \cdot \calR_{n_{\bbQ}}^{\bbQ_{X}} \left( \calF^* \right) \right\}.
    \end{split}
\end{equation}
Then, the upper bound \eqref{eqn:subsec:proof_thm:upper_bound_alg:dr_covariate_shift_adaptation_v1_v19} gives
\begin{equation}
    \label{eqn:subsec:proof_thm:upper_bound_alg:dr_covariate_shift_adaptation_v1_v20}
    \begin{split}
        \left( \bbP^{\otimes n_{\bbP}} \otimes \bbQ_{X}^{\otimes n_{\bbQ}} \right) \left\{ \calE_{\bbQ}^{(2)} (\delta) \right\} \geq 1 - \delta \quad \textnormal{for every} \quad \delta \in (0, 1).
    \end{split}
\end{equation}

\indent Finally, it is time to put all pieces together in order to bound the term (T2) from our main bound \eqref{eqn:subsec:proof_thm:upper_bound_alg:dr_covariate_shift_adaptation_v1_v4}. To this end, we introduce the event $\calE (\delta) := \calE_{\bbP} \left( \frac{\delta}{3} \right) \cap \calE_{\bbQ}^{(1)} \left( \frac{\delta}{3} \right) \cap \calE_{\bbQ}^{(2)} \left( \frac{\delta}{3} \right)$ for every $\delta \in (0, 1)$. By virtue of the union bound, the inequalities \eqref{eqn:subsec:proof_thm:upper_bound_alg:dr_covariate_shift_adaptation_v1_v12}, \eqref{eqn:subsec:proof_thm:upper_bound_alg:dr_covariate_shift_adaptation_v1_v16}, and \eqref{eqn:subsec:proof_thm:upper_bound_alg:dr_covariate_shift_adaptation_v1_v20} implies
\[
    \begin{split}
        \left( \bbP^{\otimes n_{\bbP}} \otimes \bbQ_{X}^{\otimes n_{\bbQ}} \right) \left\{ \calE (\delta) \right\} = 1 - \left( \bbP^{\otimes n_{\bbP}} \otimes \bbQ_{X}^{\otimes n_{\bbQ}} \right) \left\{ \left( \bbO^{n_{\bbP}} \times \bbX^{n_{\bbQ}} \right) \setminus \calE (\delta) \right\} \geq 1 - \delta.
    \end{split}
\]
On the other hand, by utilizing the definitions \eqref{eqn:subsec:proof_thm:upper_bound_alg:dr_covariate_shift_adaptation_v1_v11}, \eqref{eqn:subsec:proof_thm:upper_bound_alg:dr_covariate_shift_adaptation_v1_v15}, and \eqref{eqn:subsec:proof_thm:upper_bound_alg:dr_covariate_shift_adaptation_v1_v19} of the events $\calE_{\bbP} (\delta)$, $\calE_{\bbQ}^{(1)} (\delta)$, and $\calE_{\bbQ}^{(2)} (\delta)$, respectively, it follows on the event $\calE (\delta)$ that
\begin{equation}
    \label{eqn:subsec:proof_thm:upper_bound_alg:dr_covariate_shift_adaptation_v1_v21}
    \begin{split}
        \textnormal{(T2)} \stackrel{\textnormal{(i)}}{\leq} \ &  2 \sup \left\{ \bbG_{n_{\bbP}}^{\bbP} (\varphi) \left( \bfO_{1:n_{\bbP}}^{\bbP} \right) : \varphi \in \calF^* \cup \left( - \calF^* \right) \right\} \\
        &+ 2 \sup \left\{ \bbG_{n_{\bbQ}}^{\bbQ, (1)} (\varphi) \left( \bfX_{1:n_{\bbQ}}^{\bbQ} \right) : \varphi \in \calF^* \cup \left( - \calF^* \right) \right\} \\
        &+ \sup \left\{ \bbG_{n_{\bbQ}}^{\bbQ, (2)} (\varphi) \left( \bfX_{1:n_{\bbQ}}^{\bbQ} \right) : \varphi \in \calF^* \cup \left( - \calF^* \right) \right\} \\
        \leq \ & 2 \calB_{\bbP} \left( \frac{\delta}{3} \right) + 2 \calB_{\bbQ}^{(1)} \left( \frac{\delta}{3} \right) + \calB_{\bbQ}^{(2)} \left( \frac{\delta}{3} \right) \\
        &+ 8 C_{\textsf{dr}} \left( 1 + C_{\textsf{rf}} \right) \calR_{n_{\bbP}}^{\bbP_{X}} \left( \calF^* \right) + 8 \left( 3 + C_{\textsf{rf}} \right) \calR_{n_{\bbQ}}^{\bbQ_{X}} \left( \calF^* \right),
    \end{split}
\end{equation}
where the step (i) holds by virtue of the inequality \eqref{eqn:subsec:proof_thm:upper_bound_alg:dr_covariate_shift_adaptation_v1_v7}. Taking two bounds \eqref{eqn:subsec:proof_thm:upper_bound_alg:dr_covariate_shift_adaptation_v1_v5} and \eqref{eqn:subsec:proof_thm:upper_bound_alg:dr_covariate_shift_adaptation_v1_v21} collectively leads to the following upper bound on the excess $\bbQ$-risk of the DR estimator \eqref{eqn:dr_estimator_v1}: on the event $\calE (\delta)$, which holds with probability at least $1 - \delta$ under the probability measure $\bbP^{\otimes n_{\bbP}} \otimes \bbQ_{X}^{\otimes n_{\bbQ}}$, we obtain
\[
    \begin{split}
        \calE_{\bbQ} \left( \hat{f}_{\textsf{DR}} \right) \leq \ & \textnormal{(T1)} + \textnormal{(T2)} \\
        \leq \ & 4 \left\| \hat{\rho} - \rho^* \right\|_{L^2 \left( \bbX, \bbP_{X} \right)} \cdot \left\| \hat{f}_0 - f^* \right\|_{L^2 \left( \bbX, \bbP_{X} \right)} + 2 \calB_{\bbP} \left( \frac{\delta}{3} \right) + 2 \calB_{\bbQ}^{(1)} \left( \frac{\delta}{3} \right) + \calB_{\bbQ}^{(2)} \left( \frac{\delta}{3} \right) \\
        &+ 8 C_{\textsf{dr}} \left( 1 + C_{\textsf{rf}} \right) \calR_{n_{\bbP}}^{\bbP_{X}} \left( \calF^* \right) + 8 \left( 3 + C_{\textsf{rf}} \right) \calR_{n_{\bbQ}}^{\bbQ_{X}} \left( \calF^* \right) \\
        \leq \ & 4 \left\| \hat{\rho} - \rho^* \right\|_{L^2 \left( \bbX, \bbP_{X} \right)} \cdot \left\| \hat{f}_0 - f^* \right\|_{L^2 \left( \bbX, \bbP_{X} \right)} + 12 \left( 2 + C_{\textsf{rf}} \right) \log \left( \frac{3}{\delta} \right) \left( \frac{C_{\textsf{dr}}}{n_{\bbP}} + \frac{1}{n_{\bbQ}} \right) \\
        &+ 4 \left( 1 + C_{\textsf{dr}} \right) \left( 2 + C_{\textsf{rf}} \right) \sqrt{2 \log \left( \frac{3}{\delta} \right)} \left( \frac{1}{\sqrt{n_{\bbP}}} + \frac{1}{\sqrt{n_{\bbQ}}} \right) \\
        &+ 8 \left( 1 + C_{\textsf{dr}} \right) \left( 2 + C_{\textsf{rf}} \right) \sqrt{\log \left( \frac{3}{\delta} \right)} \left( \frac{\calR_{n_{\bbP}}^{\bbP_{X}} \left( \calF^* \right)}{\sqrt{n_{\bbP}}} + \frac{\calR_{n_{\bbQ}}^{\bbQ_{X}} \left( \calF^* \right)}{\sqrt{n_{\bbQ}}} \right) \\
        &+ 8 C_{\textsf{dr}} \left( 1 + C_{\textsf{rf}} \right) \calR_{n_{\bbP}}^{\bbP_{X}} \left( \calF^* \right) + 8 \left( 3 + C_{\textsf{rf}} \right) \calR_{n_{\bbQ}}^{\bbQ_{X}} \left( \calF^* \right),
    \end{split}
\]
and this completes the proof of Theorem \ref{thm:upper_bound_alg:dr_covariate_shift_adaptation_v1}.

\subsection{Proof of Theorem \ref{thm:upper_bound_alg:dr_covariate_shift_adaptation_v2}}
\label{subsec:proof_thm:upper_bound_alg:dr_covariate_shift_adaptation_v2}

\indent To begin with, we introduce some key universal constants to formally present our improved structure-agnostic guarantee of the DR estimator \eqref{eqn:dr_estimator_parameterized_case} for parameterized hypothesis classes $\calF \subseteq \left( \bbX \to \left[ -1, 1 \right] \right)$:
\begin{equation}
    \label{eqn:key_absolute_constants_thm:upper_bound_alg:dr_covariate_shift_adaptation_v2}
    \begin{split}
        \calB_1 := \ & 4 \left( 1 + C_{\textsf{dr}} \right) \left( 1 + C_{\textsf{rf}} \right) b_1, \\
        \calB_2 := \ & 8 \sqrt{2} \cdot \max \left\{ C_{\textsf{dr}} \left( 1 + C_{\textsf{rf}} \right) b_2, b_{1}^2 + \left( 1 + C_{\textsf{rf}} \right)^2 b_2 \right\}, \\
        \calB_3 := \ & \max \left\{ 2 \left( 1 + C_{\textsf{dr}} \right) \left( 1 + C_{\textsf{rf}} \right), 4 \right\} \cdot b_3 + 6 b_1 b_2.
    \end{split}
\end{equation}
With the above conventions, the detailed version of Theorem \ref{thm:upper_bound_alg:dr_covariate_shift_adaptation_v2} can be stated as follows:

\begin{thm} [Structure-agnostic upper bound \RNum{2} for the DR estimator]
\label{thm:detailed_upper_bound_alg:dr_covariate_shift_adaptation_v2}
With the parameterized function class \eqref{eqn:parametrized_hypothesis_class} and Assumptions \ref{assumption:covariate_shift}--\ref{assumption:simple_landscape_empirical_DR_risk}, the DR estimator \eqref{eqn:dr_estimator_parameterized_case} satisfies the following property: there is an absolute constant $K \in \left( 0, +\infty \right)$ such that, with probability at least $1 - 8 \delta$ under the probability measure $\bbP^{\otimes n_{\bbP}} \otimes \bbQ_{X}^{\otimes n_{\bbQ}}$, 
\begin{equation}
    \label{eqn:thm:detailed_upper_bound_alg:dr_covariate_shift_adaptation_v2_v1}
    \begin{split}
        &\calE_{\bbQ} \left( \hat{\bftheta}_{\textnormal{\textsf{DR}}} \right) = \bbE_{X \sim \bbQ_{X}} \left[ \left\{ f \left( X; \hat{\bftheta}_{\textnormal{\textsf{DR}}} \right) - f^* (X) \right\}^2 \right] \\
        \leq \ & 18 K^2 \left( 1 + C_{\textnormal{\textsf{dr}}} \right)^2 \left( 1 + C_{\textnormal{\textsf{rf}}} \right)^2 \log \left( \frac{d}{\delta} \right) \left[ \frac{\textnormal{\textsf{Trace}} \left\{ \calI_{\bbP} \left( \bftheta^* \right) \calI_{\bbQ}^{-1} \left( \bftheta^* \right) \right\}}{n_{\bbP}} + \frac{d}{n_{\bbQ}} \right],
    \end{split}
\end{equation}
provided that $\min \left\{ n_{\bbP}, n_{\bbQ} \right\} \geq \overline{\kappa} \cdot \calN^* \log \left( \frac{d}{\delta} \right)$, where $\calN^* := \max \left\{ \calN_1, \calN_2 \right\}$ and
\begin{align}
    \label{eqn:thm:detailed_upper_bound_alg:dr_covariate_shift_adaptation_v2_v2}
        \overline{\kappa} := \ & \max \left\{ \kappa, \left\{ 2 \calB_3 K \left( 1 + C_{\textnormal{\textsf{dr}}} \right) \left( 1 + C_{\textnormal{\textsf{rf}}} \right) \right\}^2, 18 \left\{ K \left( 1 + C_{\textnormal{\textsf{dr}}} \right) \left( 1 + C_{\textnormal{\textsf{rf}}} \right) \right\}^2 \right\}, \nonumber \\
        \kappa := \ & \max \left\{ \left( 260 \calB_2 \right)^2, \left\{ \frac{860 \calB_3}{3} K \left( 1 + C_{\textnormal{\textsf{dr}}} \right) \left( 1 + C_{\textnormal{\textsf{rf}}} \right) \right\}^2, \left\{ \frac{160 \calB_{1}^2 \calB_2}{\left( 1 + C_{\textnormal{\textsf{dr}}} \right)^2 \left( 1 + C_{\textnormal{\textsf{rf}}} \right)^2} \right\}^{\frac{2}{3}},  \right. \nonumber \\
        &\left. \left\{ \frac{640 \calB_{1}^3 \calB_3}{3 \left( 1 + C_{\textnormal{\textsf{dr}}} \right)^2 \left( 1 + C_{\textnormal{\textsf{rf}}} \right)^2} \right\}^{\frac{1}{2}}, \frac{80 \calB_{1}^2}{\left( 1 + C_{\textnormal{\textsf{dr}}} \right)^2 \left( 1 + C_{\textnormal{\textsf{rf}}} \right)^2} \right\}, \nonumber \\
        \calN_1 := \ & \left\| \calI_{\bbQ}^{-1} \left( \bftheta^* \right) \right\|_{\textnormal{\textsf{op}}}^2 \cdot \max \left\{ 1, \textnormal{\textsf{Trace}} \left\{ \calI_{\bbP} \left( \bftheta^* \right) \calI_{\bbQ}^{-2} \left( \bftheta^* \right) + \calI_{\bbQ}^{-1} \left( \bftheta^* \right) \right\}, \right. \nonumber \\ 
        &\left. \left[ \min \left\{ \textnormal{\textsf{Trace}} \left\{ \calI_{\bbP} \left( \bftheta^* \right) \calI_{\bbQ}^{-2} \left( \bftheta^* \right) \right\}, \textnormal{\textsf{Trace}} \left\{ \calI_{\bbQ}^{-1} \left( \bftheta^* \right) \right\} \right\} \right]^{- \frac{2}{3}}, \right. \\
        &\left. \left[ \min \left\{ \textnormal{\textsf{Trace}} \left\{ \calI_{\bbP} \left( \bftheta^* \right) \calI_{\bbQ}^{-2} \left( \bftheta^* \right) \right\}, \textnormal{\textsf{Trace}} \left\{ \calI_{\bbQ}^{-1} \left( \bftheta^* \right) \right\} \right\} \right]^{- \frac{1}{2}} \right\}, \nonumber \\
        \calN_2 := \ & \max \left\{ \left[ \frac{\textnormal{\textsf{Trace}} \left\{ \calI_{\bbP} \left( \bftheta^* \right) \calI_{\bbQ}^{-2} \left( \bftheta^* \right) + \calI_{\bbQ}^{-1} \left( \bftheta^* \right) \right\}}{\min \left\{ \textnormal{\textsf{Trace}} \left\{ \calI_{\bbP} \left( \bftheta^* \right) \calI_{\bbQ}^{-1} \left( \bftheta^* \right) \right\}, d \right\}} \right]^2, \frac{\left[ \textnormal{\textsf{Trace}} \left\{ \calI_{\bbP} \left( \bftheta^* \right) \calI_{\bbQ}^{-2} \left( \bftheta^* \right) + \calI_{\bbQ}^{-1} \left( \bftheta^* \right) \right\} \right]^3}{\left[ \min \left\{ \textnormal{\textsf{Trace}} \left\{ \calI_{\bbP} \left( \bftheta^* \right) \calI_{\bbQ}^{-1} \left( \bftheta^* \right) \right\}, d \right\} \right]^2}, \right. \nonumber \\
        &\left. \left[ \frac{\left\| \calI_{\bbQ}^{-1} \left( \bftheta^* \right) \right\|_{\textnormal{\textsf{op}}}^2}{\min \left\{ \textnormal{\textsf{Trace}} \left\{ \calI_{\bbP} \left( \bftheta^* \right) \calI_{\bbQ}^{-1} \left( \bftheta^* \right) \right\}, d \right\}} \right]^{\frac{2}{3}}, \left[ \frac{\left\| \calI_{\bbQ}^{-1} \left( \bftheta^* \right) \right\|_{\textnormal{\textsf{op}}}^3}{\min \left\{ \textnormal{\textsf{Trace}} \left\{ \calI_{\bbP} \left( \bftheta^* \right) \calI_{\bbQ}^{-1} \left( \bftheta^* \right) \right\}, d \right\}} \right]^{\frac{1}{2}}, \right. \nonumber \\
        &\left. \frac{\left\| \calI_{\bbQ}^{-1} \left( \bftheta^* \right) \right\|_{\textnormal{\textsf{op}}}}{\min \left\{ \textnormal{\textsf{Trace}} \left\{ \calI_{\bbP} \left( \bftheta^* \right) \calI_{\bbQ}^{-1} \left( \bftheta^* \right) \right\}, d \right\}} \right\}. \nonumber
\end{align}
\end{thm}

\noindent Towards proving Theorem \ref{thm:detailed_upper_bound_alg:dr_covariate_shift_adaptation_v2}, we first present a key technical lemma that plays a critical role in the proof. Roughly speaking, the following lemma aims to capture the distance between the DR estimate $\hat{\bftheta}_{\textsf{DR}} \in \Theta$ and the ground-truth parameter $\bftheta^* \in \Theta$ under different metrics.

\begin{lemma}
\label{lemma:subsec:proof_thm:upper_bound_alg:dr_covariate_shift_adaptation_v2_v1}
With Assumptions \ref{assumption:covariate_shift}--\ref{assumption:simple_landscape_empirical_DR_risk}, given any $\delta \in \left( 0, \frac{1}{8} \right]$ and $\left( n_{\bbP}, n_{\bbQ} \right) \in \bbN \times \bbN$ such that $\min \left\{ n_{\bbP}, n_{\bbQ} \right\} \geq \kappa \cdot \max \left\{ \calN_1, \calN_2 \right\} \log \left( \frac{d}{\delta} \right)$, where $\kappa$, $\calN_1$, and $\calN_2$ are specified as \eqref{eqn:thm:detailed_upper_bound_alg:dr_covariate_shift_adaptation_v2_v2}, the following facts hold with probability at least $1 - 8 \delta$ under $\bbP^{\otimes n_{\bbP}} \otimes \bbQ_{X}^{\otimes n_{\bbQ}}$: for some universal constant $K \in \left( 0, +\infty \right)$,
\begin{enumerate} [label = (\roman*)]
    \item we have $\hat{\bftheta}_{\textnormal{\textsf{DR}}} \in \bbB_{r (\delta)} \left( \bftheta^* \right)$, where the radius $r (\delta) \in \left( 0, +\infty \right)$ is given by
    \begin{equation}
        \label{eqn:lemma:subsec:proof_thm:upper_bound_alg:dr_covariate_shift_adaptation_v2_v1_v1}
        \begin{split}
            r (\delta) := \ & 3 K \left( 1 + C_{\textnormal{\textsf{dr}}} \right) \left( 1 + C_{\textnormal{\textsf{rf}}} \right) \sqrt{\log \left( \frac{d}{\delta} \right)} \\
            &\left[ \sqrt{\frac{\textnormal{\textsf{Trace}} \left\{ \calI_{\bbP} \left( \bftheta^* \right) \calI_{\bbQ}^{-2} \left( \bftheta^* \right) \right\}}{n_{\bbP}}} + \sqrt{\frac{\textnormal{\textsf{Trace}} \left\{ \calI_{\bbQ}^{-1} \left( \bftheta^* \right) \right\}}{n_{\bbQ}}} \right].
        \end{split}
    \end{equation}
    \item it holds that
    \begin{equation}
        \label{eqn:lemma:subsec:proof_thm:upper_bound_alg:dr_covariate_shift_adaptation_v2_v1_v2}
        \begin{split}
            &\left\| \calI_{\bbQ}^{\frac{1}{2}} \left( \bftheta^* \right) \left( \hat{\bftheta}_{\textnormal{\textsf{DR}}} - \bftheta^* \right) \right\|_{2}^2 \\
            \leq \ & 9 K^2 \left( 1 + C_{\textnormal{\textsf{dr}}} \right)^2 \left( 1 + C_{\textnormal{\textsf{rf}}} \right)^2 \log \left( \frac{d}{\delta} \right) \left[ \sqrt{\frac{\textnormal{\textsf{Trace}} \left\{ \calI_{\bbP} \left( \bftheta^* \right) \calI_{\bbQ}^{-1} \left( \bftheta^* \right) \right\}}{n_{\bbP}}} + \sqrt{\frac{d}{n_{\bbQ}}} \right]^2.
        \end{split}
    \end{equation}
\end{enumerate}
\end{lemma}

\noindent For simplicity, let $\Lambda (\delta) \subseteq \bbO^{n_{\bbP}} \times \bbX^{n_{\bbQ}}$ denote the event for which Lemma \ref{lemma:subsec:proof_thm:upper_bound_alg:dr_covariate_shift_adaptation_v2_v1} holds, which immediately implies $\left( \bbP^{\otimes n_{\bbP}} \otimes \bbQ_{X}^{\otimes n_{\bbQ}} \right) \left\{ \Lambda (\delta) \right\} \geq 1 - 8 \delta$ for any given $\delta \in \left( 0, \frac{1}{8} \right]$.
\medskip

\indent We embark on the proof of Theorem \ref{thm:detailed_upper_bound_alg:dr_covariate_shift_adaptation_v2} by doing a Taylor expansion with respect to as follows:
\begin{equation}
    \label{eqn:subsec:proof_thm:upper_bound_alg:dr_covariate_shift_adaptation_v2_v1}
    \begin{split}
        \calE_{\bbQ} \left( \hat{\bftheta}_{\textsf{DR}} \right) = \ & \bbE_{(X, Y) \sim \bbQ} \left[ \ell \left( Y, f \left( X; \hat{\bftheta}_{\textsf{DR}} \right) \right) - \ell \left( Y, f \left( X; \bftheta^* \right) \right) \right] \\
        = \ & \left\{ \bbE_{(X, Y) \sim \bbQ} \left[ \nabla_{\bftheta} \ell \left( Y, f \left( X; \bftheta^* \right) \right) \right] \right\}^{\top} \left( \hat{\bftheta}_{\textsf{DR}} - \bftheta^* \right) \\
        &+ \frac{1}{2!} \left( \hat{\bftheta}_{\textsf{DR}} - \bftheta^* \right)^{\top} \bbE_{(X, Y) \sim \bbQ} \left[ \nabla_{\bftheta}^2 \ell \left( Y, f \left( X; \bftheta^* \right) \right) \right] \left( \hat{\bftheta}_{\textsf{DR}} - \bftheta^* \right) \\
        &+ \frac{1}{3!} \left\langle \bbE_{(X, Y) \sim \bbQ} \left[ \nabla_{\bftheta}^3 \ell \left( Y, f \left( X; \tilde{\bftheta} \right) \right) \right],  \left( \hat{\bftheta}_{\textsf{DR}} - \bftheta^* \right)^{\otimes 3} \right\rangle_{\textsf{F}} \\
        \stackrel{\textnormal{(a)}}{=} \ & \frac{1}{2!} \left( \hat{\bftheta}_{\textsf{DR}} - \bftheta^* \right)^{\top} \calI_{\bbQ} \left( \bftheta^* \right) \left( \hat{\bftheta}_{\textsf{DR}} - \bftheta^* \right) \\
        &+ \frac{1}{3!} \left\langle \bbE_{(X, Y) \sim \bbQ} \left[ \nabla_{\bftheta}^3 \ell \left( Y, f \left( X; \tilde{\bftheta} \right) \right) \right],  \left( \hat{\bftheta}_{\textsf{DR}} - \bftheta^* \right)^{\otimes 3} \right\rangle_{\textsf{F}}
    \end{split}
\end{equation}
for some $\tilde{\bftheta} \in \left\{ (1 - \lambda) \bftheta^* + \lambda \hat{\bftheta}_{\textsf{DR}} : \lambda \in [0, 1] \right\}$, where the step (a) holds due to the following facts:
\begin{equation}
    \label{eqn:subsec:proof_thm:upper_bound_alg:dr_covariate_shift_adaptation_v2_v2}
    \begin{split}
        \nabla_{\bftheta} \ell \left( y, f \left( x; \bftheta \right) \right) = \ & 2 \left\{ f \left( x; \bftheta \right) - y \right\} \nabla_{\bftheta} f \left( x; \bftheta \right), \\
        \nabla_{\bftheta}^2 \ell \left( y, f \left( x; \bftheta \right) \right) = \ & 2 \nabla_{\bftheta} f \left( x; \bftheta \right) \left\{ \nabla_{\bftheta} f \left( x; \bftheta \right) \right\}^{\top} + 2 \left\{ f \left( x; \bftheta \right) - y \right\} \nabla_{\bftheta}^2 f \left( x; \bftheta \right), \\
        \nabla_{\bftheta}^3 \ell \left( y, f \left( x; \bftheta \right) \right) = \ & 2 \nabla_{\bftheta}^2 f \left( x; \bftheta \right) \otimes \nabla_{\bftheta} f \left( x; \bftheta \right) + 4 \nabla_{\bftheta} f \left( x; \bftheta \right) \otimes \nabla_{\bftheta}^2 f \left( x; \bftheta \right) \\
        &+ 2 \left\{ f \left( x; \bftheta \right) - y \right\} \nabla_{\bftheta}^3 f \left( x; \bftheta \right)
    \end{split}
\end{equation}
Thus, it follows from the equation \eqref{eqn:subsec:proof_thm:upper_bound_alg:dr_covariate_shift_adaptation_v2_v1} that
\begin{equation}
    \label{eqn:subsec:proof_thm:upper_bound_alg:dr_covariate_shift_adaptation_v2_v3}
    \begin{split}
        \calE_{\bbQ} \left( \hat{\bftheta}_{\textsf{DR}} \right) \leq \ & \frac{1}{2} \left( \hat{\bftheta}_{\textsf{DR}} - \bftheta^* \right)^{\top} \calI_{\bbQ} \left( \bftheta^* \right) \left( \hat{\bftheta}_{\textsf{DR}} - \bftheta^* \right) \\
        &+ \frac{1}{6} \left\| \bbE_{(X, Y) \sim \bbQ} \left[ \nabla_{\bftheta}^3 \ell \left( Y, f \left( X; \tilde{\bftheta} \right) \right) \right] \right\|_{\textsf{op}} \left\| \hat{\bftheta}_{\textsf{DR}} - \bftheta^* \right\|_{2}^3 \\
        \stackrel{\textnormal{(b)}}{\leq} \ & \frac{1}{2} \left( \hat{\bftheta}_{\textsf{DR}} - \bftheta^* \right)^{\top} \calI_{\bbQ} \left( \bftheta^* \right) \left( \hat{\bftheta}_{\textsf{DR}} - \bftheta^* \right) + \frac{\calB_3}{6} \left\| \hat{\bftheta}_{\textsf{DR}} - \bftheta^* \right\|_{2}^3,
    \end{split}
\end{equation}
where the step (b) holds by the observation that the operator norm $\left\| \cdot \right\|_{\textsf{op}}: \left( \bbR^d \right)^{\otimes 3} \to \bbR_{+}$ is a convex function together with Jensen's inequality and the following bound: for any $\bftheta \in \Theta$,
\[
    \begin{split}
        \left\| \nabla_{\bftheta}^3 \ell \left( Y, f \left( X; \bftheta \right) \right) \right\|_{\textsf{op}} \leq \ & 2 \left( 1 + |Y| \right) \left\| \nabla_{\bftheta}^3 f \left( X; \bftheta \right) \right\|_{\textsf{op}} + 6 \left\| \nabla_{\bftheta} f \left( X; \bftheta \right) \right\|_{2} \left\| \nabla_{\bftheta}^2 f \left( X; \bftheta \right) \right\|_{\textsf{op}} \\
        \stackrel{\textnormal{(c)}}{\leq} \ & 4 b_3 + 6 b_1 b_2 \leq \calB_3
    \end{split}
\]
$\bbQ$-almost surely, where the step (c) follows from Assumption \ref{assumption:uniform_boundedness} and the part (\romannumeral 2) of Assumption \ref{assumption:smoothness}. Therefore, while being conditioned on the event $\Lambda (\delta)$, we reach
\begin{equation}
    \label{eqn:subsec:proof_thm:upper_bound_alg:dr_covariate_shift_adaptation_v2_v4}
    \begin{split}
        &\calE_{\bbQ} \left( \hat{\bftheta}_{\textsf{DR}} \right) \\
        \leq \ & \frac{9}{2} K^2 \left( 1 + C_{\textsf{dr}} \right)^2 \left( 1 + C_{\textsf{rf}} \right)^2 \log \left( \frac{d}{\delta} \right) \left[ \sqrt{\frac{\textsf{Trace} \left\{ \calI_{\bbP} \left( \bftheta^* \right) \calI_{\bbQ}^{-1} \left( \bftheta^* \right) \right\}}{n_{\bbP}}} + \sqrt{\frac{d}{n_{\bbQ}}} \right]^2 \\
        &+ \frac{3}{2} \calB_3 K^3 \left( 1 + C_{\textsf{dr}} \right)^3 \left( 1 + C_{\textsf{rf}} \right)^3 \log^{\frac{3}{2}} \left( \frac{d}{\delta} \right) \\
        &\cdot \left[ \sqrt{\frac{\textsf{Trace} \left\{ \calI_{\bbP} \left( \bftheta^* \right) \calI_{\bbQ}^{-2} \left( \bftheta^* \right) \right\}}{n_{\bbP}}} + \sqrt{\frac{\textsf{Trace} \left\{ \calI_{\bbQ}^{-1} \left( \bftheta^* \right) \right\}}{n_{\bbQ}}} \right]^3. 
    \end{split}
\end{equation}
At this point, one can observe that if $\min \left\{ n_{\bbP}, n_{\bbQ} \right\} \geq \overline{\kappa} \cdot \max \left\{ \calN_1, \calN_2 \right\} \log \left( \frac{d}{\delta} \right)$, then
\begin{equation}
    \label{eqn:subsec:proof_thm:upper_bound_alg:dr_covariate_shift_adaptation_v2_v5}
    \begin{split}
        &\frac{3}{2} \calB_3 K^3 \left( 1 + C_{\textsf{dr}} \right)^3 \left( 1 + C_{\textsf{rf}} \right)^3 \log^{\frac{3}{2}} \left( \frac{d}{\delta} \right) \\
        &\cdot \left[ \sqrt{\frac{\textsf{Trace} \left\{ \calI_{\bbP} \left( \bftheta^* \right) \calI_{\bbQ}^{-2} \left( \bftheta^* \right) \right\}}{n_{\bbP}}} + \sqrt{\frac{\textsf{Trace} \left\{ \calI_{\bbQ}^{-1} \left( \bftheta^* \right) \right\}}{n_{\bbQ}}} \right]^3 \\
        \leq \ & \frac{9}{2} K^2 \left( 1 + C_{\textsf{dr}} \right)^2 \left( 1 + C_{\textsf{rf}} \right)^2 \log \left( \frac{d}{\delta} \right) \left[ \sqrt{\frac{\textsf{Trace} \left\{ \calI_{\bbP} \left( \bftheta^* \right) \calI_{\bbQ}^{-1} \left( \bftheta^* \right) \right\}}{n_{\bbP}}} + \sqrt{\frac{d}{n_{\bbQ}}} \right]^2.
    \end{split}
\end{equation}
Hence, by taking two pieces \eqref{eqn:subsec:proof_thm:upper_bound_alg:dr_covariate_shift_adaptation_v2_v4} and \eqref{eqn:subsec:proof_thm:upper_bound_alg:dr_covariate_shift_adaptation_v2_v5} collectively, it holds that if $\min \left\{ n_{\bbP}, n_{\bbQ} \right\} \geq \overline{\kappa} \cdot \max \left\{ \calN_1, \calN_2 \right\} \log \left( \frac{d}{\delta} \right)$, then we have on the event $\Lambda (\delta)$ that
\begin{equation}
    \label{eqn:subsec:proof_thm:upper_bound_alg:dr_covariate_shift_adaptation_v2_v6}
    \begin{split}
        &\calE_{\bbQ} \left( \hat{\bftheta}_{\textsf{DR}} \right) \\
        \leq \ & 9 K^2 \left( 1 + C_{\textsf{dr}} \right)^2 \left( 1 + C_{\textsf{rf}} \right)^2 \log \left( \frac{d}{\delta} \right) \left[ \sqrt{\frac{\textsf{Trace} \left\{ \calI_{\bbP} \left( \bftheta^* \right) \calI_{\bbQ}^{-1} \left( \bftheta^* \right) \right\}}{n_{\bbP}}} + \sqrt{\frac{d}{n_{\bbQ}}} \right]^2 \\
        \stackrel{\textnormal{(d)}}{\leq} \ & 18 K^2 \left( 1 + C_{\textsf{dr}} \right)^2 \left( 1 + C_{\textsf{rf}} \right)^2 \log \left( \frac{d}{\delta} \right) \left[ \frac{\textsf{Trace} \left\{ \calI_{\bbP} \left( \bftheta^* \right) \calI_{\bbQ}^{-1} \left( \bftheta^* \right) \right\}}{n_{\bbP}} + \frac{d}{n_{\bbQ}} \right],
    \end{split}
\end{equation}
where the step (d) invokes the Cauchy-Schwarz inequality. Since $\left( \bbP^{\otimes n_{\bbP}} \otimes \bbQ_{X}^{\otimes n_{\bbQ}} \right) \left\{ \Lambda (\delta) \right\} \geq 1 - 8 \delta$, the upper bound \eqref{eqn:subsec:proof_thm:upper_bound_alg:dr_covariate_shift_adaptation_v2_v6} on the excess $\bbQ$-risk of the DR estimator \eqref{eqn:dr_estimator_parameterized_case} holds with probability higher than $1 - 8 \delta$ under the probability measure $\bbP^{\otimes n_{\bbP}} \otimes \bbQ_{X}^{\otimes n_{\bbQ}}$, which completes the proof of Theorem \ref{thm:detailed_upper_bound_alg:dr_covariate_shift_adaptation_v2}.

\section{Proof of auxiliary lemmas for the proof of Theorem \ref{thm:upper_bound_alg:dr_covariate_shift_adaptation_v1}}
\label{sec:proof_auxiliary_lemmas_subsec:proof_thm:upper_bound_alg:dr_covariate_shift_adaptation_v1}

\subsection{Proof of Lemma \ref{lemma:control_P_empirical_process_v1}}
\label{subsec:proof_lemma:control_P_empirical_process_v1}

By following the standard symmetrization argument from empirical processes theory, we find that
\begin{align}
    \label{eqn:subsec:proof_lemma:control_P_empirical_process_v1_v1}
        &\bbE_{\bfO_{1:n_{\bbP}}^{\bbP} \sim \bbP^{\otimes n_{\bbP}}} \left[ \sup \left\{ \bbG_{n_{\bbP}}^{\bbP} (\varphi) \left( \bfO_{1:n_{\bbP}}^{\bbP} \right) : \varphi \in \calF^* \cup \left( - \calF^* \right) \right\} \right] \nonumber \\ 
        = \ & \bbE_{\bfO_{1:n_{\bbP}}^{\bbP} \sim \bbP^{\otimes n_{\bbP}}} \left[ \sup \left\{ \left| \bbG_{n_{\bbP}}^{\bbP} (\varphi) \left( \bfO_{1:n_{\bbP}}^{\bbP} \right) \right| : \varphi \in \calF^* \right\} \right] \\
        \leq \ & 2 \bbE_{\left( \bfO_{1:n_{\bbP}}^{\bbP}, \bfsigma_{1:n_{\bbP}} \right) \sim \bbP^{\otimes n_{\bbP}} \otimes \textsf{Unif} \left( \left\{ \pm 1 \right\}^{n_{\bbP}} \right)} \left[ \sup \left\{ \left| \frac{1}{n_{\bbP}} \sum_{i=1}^{n_{\bbP}} \sigma_{i} \hat{\rho} \left( X_{i}^{\bbP} \right) \varphi \left( X_{i}^{\bbP} \right) \left\{ Y_{i}^{\bbP} - \hat{f}_0 \left( X_{i}^{\bbP} \right) \right\} \right| : \varphi \in \calF^* \right\} \right]. \nonumber
\end{align}
To control the last term in the bound \eqref{eqn:subsec:proof_lemma:control_P_empirical_process_v1_v1}, we leverage the Ledoux-Talagrand contraction principle, which is formally stated in Lemma \ref{lemma:ledoux_talagrand_contraction}. To this end, define the functions $\phi_{i}^{\bbP}: \bbO^{n_{\bbP}} \to \left( \bbR \to \bbR \right)$, $i \in \left[ n_{\bbP} \right]$, by
\begin{equation}
    \label{eqn:subsec:proof_lemma:control_P_empirical_process_v1_v2}
        \phi_{i}^{\bbP} \left( \bfo_{1:n_{\bbP}}^{\bbP} \right) (t) := \hat{\rho} \left( x_{i}^{\bbP} \right) \left\{ y_{i}^{\bbP} - \hat{f}_0 \left( x_{i}^{\bbP} \right) \right\} (t), \quad \forall t \in \bbR,
\end{equation}
and let $\phi_{i}^{\bbP} := \phi_{i}^{\bbP} \left( \bfO_{1:n_{\bbP}}^{\bbP} \right) : \bbR \to \bbR$ for $\bfO_{1:n_{\bbP}}^{\bbP} \sim \bbP^{\otimes n_{\bbP}}$ for simplicity. Let $\calA_{\bbP} := \bigcap_{i=1}^{n_{\bbP}} \left\{ \left| Y_{i}^{\bbP} \right| \leq 1 \right\}$, which holds $\bbP^{\otimes n_{\bbP}}$-almost surely. Then, one can observe that the function $\phi_{i}^{\bbP}: \bbO^{n_{\bbP}} \to \left( \bbR \to \bbR \right)$ is a $C_{\textsf{dr}} \left( 1 + C_{\textsf{rf}} \right)$-Lipschitz continuous function such that $\phi_{i}^{\bbP} (0) = 0$ for every $i \in \left[ n_{\bbP} \right]$ on the event $\calA_{\bbP}$. Then, while being conditioned on $\bfO_{1:n_{\bbP}}^{\bbP}$, we obtain by virtue of Lemma \ref{lemma:ledoux_talagrand_contraction} that
\begin{align}
    \label{eqn:subsec:proof_lemma:control_P_empirical_process_v1_v3}
        &\bbE_{\bfsigma_{1:n_{\bbP}} \sim \textsf{Unif} \left( \left\{ \pm 1 \right\}^{n_{\bbP}} \right)} \left[ \sup \left\{ \left| \frac{1}{n_{\bbP}} \sum_{i=1}^{n_{\bbP}} \sigma_{i} \hat{\rho} \left( X_{i}^{\bbP} \right) \varphi \left( X_{i}^{\bbP} \right) \left\{ Y_{i}^{\bbP} - \hat{f}_0 \left( X_{i}^{\bbP} \right) \right\} \right| : \varphi \in \calF^* \right\} \cdot \mathbbm{1}_{\calA_{\bbP}} \right] \nonumber \\
        = \ & \bbE_{\bfsigma_{1:n_{\bbP}} \sim \textsf{Unif} \left( \left\{ \pm 1 \right\}^{n_{\bbP}} \right)} \left[ \sup \left\{ \left| \frac{1}{n_{\bbP}} \sum_{i=1}^{n_{\bbP}} \sigma_{i} \hat{\rho} \left( X_{i}^{\bbP} \right) \varphi \left( X_{i}^{\bbP} \right) \left\{ Y_{i}^{\bbP} - \hat{f}_0 \left( X_{i}^{\bbP} \right) \right\} \right| : \varphi \in \calF^* \right\} \right] \cdot \mathbbm{1}_{\calA_{\bbP}} \nonumber \\
        = \ & \bbE_{\bfsigma_{1:n_{\bbP}} \sim \textsf{Unif} \left( \left\{ \pm 1 \right\}^{n_{\bbP}} \right)} \left[ \sup \left\{ \left| \frac{1}{n_{\bbP}} \sum_{i=1}^{n_{\bbP}} \sigma_{i} \phi_{i}^{\bbP} \left( t_i \right) \right|: \bft_{1:n_{\bbP}} \in \calT_{\bbP} \left( \bfO_{1:n_{\bbP}}^{\bbP} \right) \right\} \right] \cdot \mathbbm{1}_{\calA_{\bbP}} \\
        \leq \ & 2 C_{\textsf{dr}} \left( 1 + C_{\textsf{rf}} \right) \bbE_{\bfsigma_{1:n_{\bbP}} \sim \textsf{Unif} \left( \left\{ \pm 1 \right\}^{n_{\bbP}} \right)} \left[ \sup \left\{ \left| \frac{1}{n_{\bbP}} \sum_{i=1}^{n_{\bbP}} \sigma_{i} t_i \right|: \bft_{1:n_{\bbP}} \in \calT_{\bbP} \left( \bfO_{1:n_{\bbP}}^{\bbP} \right) \right\} \right] \cdot \mathbbm{1}_{\calA_{\bbP}}, \nonumber
\end{align}
where $\calT_{\bbP}: \bbO^{n_{\bbP}} \to \calP \left( \bbR^{n_{\bbP}} \right)$ is defined as
\[
    \calT_{\bbP} \left( \bfo_{1:n_{\bbP}}^{\bbP} \right) := \left\{ \left( f \left( x_{i}^{\bbP} \right) - f^* \left( x_{i}^{\bbP} \right) : i \in \left[ n_{\bbP} \right] \right) : f \in \calF \right\} \subseteq \bbR^{n_{\bbP}}.
\]
Here, $\calP \left( \bbR^{n_{\bbP}} \right)$ denotes the power set of the $n_{\bbP}$-dimensional Euclidean space $\bbR^{n_{\bbP}}$. By taking the bound \eqref{eqn:subsec:proof_lemma:control_P_empirical_process_v1_v3} collectively into the inequality \eqref{eqn:subsec:proof_lemma:control_P_empirical_process_v1_v1}, we arrive at
\[
    \begin{split}
        &\bbE_{\bfO_{1:n_{\bbP}}^{\bbP} \sim \bbP^{\otimes n_{\bbP}}} \left[ \sup \left\{ \bbG_{n_{\bbP}}^{\bbP} (\varphi) \left( \bfO_{1:n_{\bbP}}^{\bbP} \right) : \varphi \in \calF^* \cup \left( - \calF^* \right) \right\} \right] \\
        \leq \ & 2 \bbE_{\left( \bfO_{1:n_{\bbP}}^{\bbP}, \bfsigma_{1:n_{\bbP}} \right) \sim \bbP^{\otimes n_{\bbP}} \otimes \textsf{Unif} \left( \left\{ \pm 1 \right\}^{n_{\bbP}} \right)} \left[ \sup \left\{ \left| \frac{1}{n_{\bbP}} \sum_{i=1}^{n_{\bbP}} \sigma_{i} \hat{\rho} \left( X_{i}^{\bbP} \right) \varphi \left( X_{i}^{\bbP} \right) \left\{ Y_{i}^{\bbP} - \hat{f}_0 \left( X_{i}^{\bbP} \right) \right\} \right| : \varphi \in \calF^* \right\} \cdot \mathbbm{1}_{\calA_{\bbP}} \right] \\
        \leq \ & 4 C_{\textsf{dr}} \left( 1 + C_{\textsf{rf}} \right) \bbE_{\bfO_{1:n_{\bbP}}^{\bbP} \sim \bbP^{\otimes n_{\bbP}}} \left[ \bbE_{\bfsigma_{1:n_{\bbP}} \sim \textsf{Unif} \left( \left\{ \pm 1 \right\}^{n_{\bbP}} \right)} \left[ \sup \left\{ \left| \frac{1}{n_{\bbP}} \sum_{i=1}^{n_{\bbP}} \sigma_{i} t_i \right|: \bft_{1:n_{\bbP}} \in \calT_{\bbP} \left( \bfO_{1:n_{\bbP}}^{\bbP} \right) \right\} \right] \cdot \mathbbm{1}_{\calA_{\bbP}} \right] \\
        = \ & 4 C_{\textsf{dr}} \left( 1 + C_{\textsf{rf}} \right) \bbE_{\bfO_{1:n_{\bbP}}^{\bbP} \sim \bbP^{\otimes n_{\bbP}}} \left[ \widehat{\calR}_{n_{\bbP}} \left( \calF^* \right) \left( \bfX_{1:n_{\bbP}}^{\bbP} \right) \cdot \mathbbm{1}_{\calA_{\bbP}} \right] \\
        = \ & 4 C_{\textsf{dr}} \left( 1 + C_{\textsf{rf}} \right) \calR_{n_{\bbP}}^{\bbP_{X}} \left( \calF^* \right), 
    \end{split}
\]
as desired.

\subsection{Proof of Lemma \ref{lemma:control_P_empirical_process_v2}}
\label{subsec:proof_lemma:control_P_empirical_process_v2}

\indent To begin with, we introduce the function class $\left\{ \theta_{\bbP} (\varphi) : \varphi \in \calF^* \cup \left( - \calF^* \right) \right\}$, where $\theta_{\bbP} (\varphi): \bbX \times \left[ -1, 1 \right] \to \bbR$ is defined to be
\[
    \theta_{\bbP} (\varphi) (x, y) := \hat{\rho} (x) \varphi (x) \left\{ y - \hat{f}_0 (x) \right\}, \quad \forall (x, y) \in \bbX \times \left[ -1, 1 \right].
\]
Recall that the event $\calA_{\bbP} = \bigcap_{i=1}^{n_{\bbP}} \left\{ \left| Y_{i}^{\bbP} \right| \leq 1 \right\}$ holds $\bbP^{\otimes n_{\bbP}}$-almost surely, i.e., $\left( \bbP^{\otimes n_{\bbP}} \right) \left( \calA_{\bbP} \right) = 1$. We thus obtain that
\begin{equation}
    \label{eqn:subsec:proof_lemma:control_P_empirical_process_v2_v1}
    \begin{split}
        &\sup \left\{ \bbG_{n_{\bbP}}^{\bbP} (\varphi) \left( \bfO_{1:n_{\bbP}}^{\bbP} \right) : \varphi \in \calF^* \cup \left( - \calF^* \right) \right\} \\
        = \ & \sup \left\{ \bbG_{n_{\bbP}}^{\bbP} (\varphi) \left( \bfO_{1:n_{\bbP}}^{\bbP} \right) : \varphi \in \calF^* \cup \left( - \calF^* \right) \right\} \cdot \mathbbm{1}_{\calA_{\bbP}} \\
        = \ & \sup \left\{ \frac{1}{n_{\bbP}} \sum_{i=1}^{n_{\bbP}} \theta_{\bbP} (\varphi) \left( O_{i}^{\bbP} \right) - \bbE_{O \sim \bbP} \left[ \theta_{\bbP} (\varphi) (O) \right] : \varphi \in \calF^* \cup \left( - \calF^* \right) \right\} \cdot \mathbbm{1}_{\calA_{\bbP}} \\
        = \ & \sup \left\{ \frac{1}{n_{\bbP}} \sum_{i=1}^{n_{\bbP}} \theta_{\bbP} (\varphi) \left( O_{i}^{\bbP} \right) - \bbE_{O \sim \bbP} \left[ \theta_{\bbP} (\varphi) (O) \right] : \varphi \in \calF^* \cup \left( - \calF^* \right) \right\}
    \end{split}
\end{equation}
$\bbP^{\otimes n}$-almost surely. In light of the equation \eqref{eqn:subsec:proof_lemma:control_P_empirical_process_v2_v1}, it can be easily seen based on the assumption \eqref{eqn:uniform_boundedness_black_box_estimates} that
\begin{enumerate} [label = (\roman*)]
    \item $\left| \theta_{\bbP} (\varphi) (x, y) \right| \leq 2 C_{\textsf{dr}} \left( 1 +C_{\textsf{rf}} \right)$ for every $\left( x, y, \varphi \right) \in \bbX \times \left[ -1, 1 \right] \times \left\{ \calF^* \cup \left( - \calF^* \right) \right\}$.
    \item $\textsf{Var}_{(X, Y) \sim \bbP} \left[ \theta_{\bbP} (\varphi) (X, Y) \right] \leq \bbE_{(X, Y) \sim \bbP} \left[ \left\{ \theta_{\bbP} (\varphi) (X, Y) \right\}^2 \right] \leq 4 C_{\textsf{dr}}^2 \left( 1 +C_{\textsf{rf}} \right)^2$ for every $\varphi \in \calF^* \cup \left( - \calF^* \right)$.
\end{enumerate}
By the classical Talagrand's concentration inequality (Lemma \ref{lemma:talagrand_concentration}) with $\left( B, v^2 \right) = \left( 2 C_{\textsf{dr}} \left( 1 +C_{\textsf{rf}} \right), 4 C_{\textsf{dr}}^2 \left( 1 +C_{\textsf{rf}} \right)^2 \right)$, it holds with probability at least $1 - \delta$ that
\begin{align}
    \label{eqn:subsec:proof_lemma:control_P_empirical_process_v2_v2}
    &\sup \left\{ \bbG_{n_{\bbP}}^{\bbP} (\varphi) \left( \bfO_{1:n_{\bbP}}^{\bbP} \right) : \varphi \in \calF^* \cup \left( - \calF^* \right) \right\} \nonumber \\
    &- \bbE_{\bfO_{1:n_{\bbP}}^{\bbP} \sim \bbP^{\otimes n_{\bbP}}} \left[ \sup \left\{ \bbG_{n_{\bbP}}^{\bbP} (\varphi) \left( \bfO_{1:n_{\bbP}}^{\bbP} \right) : \varphi \in \calF^* \cup \left( - \calF^* \right) \right\} \right] \nonumber \\
    \leq \ & \frac{6 C_{\textsf{dr}} \left( 1 + C_{\textsf{rf}} \right)}{n_{\bbP}} \log \left( \frac{1}{\delta} \right) + 2 C_{\textsf{dr}} \left( 1 + C_{\textsf{rf}} \right) \sqrt{\frac{2 \log \left( \frac{1}{\delta} \right)}{n_{\bbP}}} \\
    &+ \frac{2}{\sqrt{n_{\bbP}}} \sqrt{ 2 C_{\textsf{dr}} \left( 1 + C_{\textsf{rf}} \right) \bbE_{\bfO_{1:n_{\bbP}}^{\bbP} \sim \bbP^{\otimes n_{\bbP}}} \left[ \sup \left\{ \bbG_{n_{\bbP}}^{\bbP} (\varphi) \left( \bfO_{1:n_{\bbP}}^{\bbP} \right) : \varphi \in \calF^* \cup \left( - \calF^* \right) \right\} \right] \log \left( \frac{1}{\delta} \right)} \nonumber \\
    \stackrel{\textnormal{(a)}}{\leq} \ & \frac{6 C_{\textsf{dr}} \left( 1 + C_{\textsf{rf}} \right)}{n_{\bbP}} \log \left( \frac{1}{\delta} \right) + 2 C_{\textsf{dr}} \left( 1 + C_{\textsf{rf}} \right) \sqrt{\frac{2 \log \left( \frac{1}{\delta} \right)}{n_{\bbP}}} + 4 C_{\textsf{dr}} \left( 1 + C_{\textsf{rf}} \right) \sqrt{\frac{2 \calR_{n_{\bbP}}^{\bbP_{X}} \left( \calF^* \right) \log \left( \frac{1}{\delta} \right)}{n_{\bbP}}} \nonumber,
\end{align}
where the step (a) invokes Lemma \ref{lemma:control_P_empirical_process_v1}. We thus complete the proof of Lemma \ref{lemma:control_P_empirical_process_v2}.

\subsection{Proof of Lemma \ref{lemma:control_Q_empirical_process_v1}}
\label{subsec:proof_lemma:control_Q_empirical_process_v1}

\indent In light of the standard symmetrization argument from empirical processes theory, we reveal that
\begin{align}
    \label{eqn:subsec:proof_lemma:control_Q_empirical_process_v1_v1}
        &\bbE_{\bfX_{1:n_{\bbQ}}^{\bbQ} \sim \bbQ_{X}^{\otimes n_{\bbQ}}} \left[ \sup \left\{ \bbG_{n_{\bbQ}}^{\bbQ, (1)} (\varphi) \left( \bfX_{1:n_{\bbQ}}^{\bbQ} \right) : \varphi \in \calF^* \cup \left( - \calF^* \right) \right\} \right] \nonumber \\ 
        = \ & \bbE_{\bfX_{1:n_{\bbQ}}^{\bbQ} \sim \bbQ_{X}^{\otimes n_{\bbQ}}} \left[ \sup \left\{ \left| \bbG_{n_{\bbQ}}^{\bbQ, (1)} (\varphi) \left( \bfX_{1:n_{\bbQ}}^{\bbQ} \right) \right| : \varphi \in \calF^* \right\} \right] \\
        \leq \ & 2 \bbE_{\left( \bfX_{1:n_{\bbQ}}^{\bbQ}, \bfsigma_{1:n_{\bbQ}} \right) \sim \bbQ_{X}^{\otimes n_{\bbQ}} \otimes \textsf{Unif} \left( \left\{ \pm 1 \right\}^{n_{\bbQ}} \right)} \left[ \sup \left\{ \left| \frac{1}{n_{\bbQ}} \sum_{j=1}^{n_{\bbQ}} \sigma_{j} \varphi \left( X_{j}^{\bbQ} \right) \left\{ \hat{f}_0 \left( X_{j}^{\bbQ} \right) - f^* \left( X_{j}^{\bbQ} \right) \right\} \right| : \varphi \in \calF^* \right\} \right]. \nonumber
\end{align}
    
\noindent We now focus on a tight control of the final term of the bound \eqref{eqn:subsec:proof_lemma:control_Q_empirical_process_v1_v1} via the Ledoux-Talagrand contraction principle (Lemma \ref{lemma:ledoux_talagrand_contraction}). Towards this end, let us consider the functions $\alpha_{j}^{\bbQ} : \bbX^{n_{\bbQ}} \to \left( \left[ -2, 2 \right] \to \bbR \right)$ for $j \in \left[ n_{\bbQ} \right]$ defined as
\begin{equation}
    \label{eqn:subsec:proof_lemma:control_Q_empirical_process_v1_v2}
    \begin{split}
        \alpha_{j}^{\bbQ} \left( \bfx_{1:n_{\bbQ}}^{\bbQ} \right) (t) := \left\{ \hat{f}_0 \left( x_{j}^{\bbQ} \right) - f^* \left( x_{j}^{\bbQ} \right) \right\} t, \quad \forall t \in \left[ -2, 2 \right],
    \end{split}
\end{equation}
and simplify $\alpha_{j}^{\bbQ} := \alpha_{j}^{\bbQ} \left( \bfX_{1:n_{\bbQ}}^{\bbQ} \right): \left[ -2, 2 \right] \to \bbR$, where $\bfX_{1:n_{\bbQ}}^{\bbQ} \sim \bbQ_{X}^{n_{\bbQ}}$. Then, one can find that $\alpha_{j}^{\bbQ} : \left[ -2, 2 \right] \to \bbR$ is an $\left( 1 + C_{\textsf{rf}} \right)$-Lipschitz continuous function with $\alpha_{j}^{\bbQ} (0) = 0$ for $j \in \left[ n_{\bbQ} \right]$. By applying the Ledoux-Talagrand contraction principle (Lemma \ref{lemma:ledoux_talagrand_contraction}), while being conditioned on $\bfX_{1:n_{\bbQ}}^{\bbQ}$, we now have
\begin{equation}
    \label{eqn:subsec:proof_lemma:control_Q_empirical_process_v1_v3}
    \begin{split}
        &\bbE_{\bfsigma_{1:n_{\bbQ}} \sim \textsf{Unif} \left( \left\{ \pm 1 \right\}^{n_{\bbQ}} \right)} \left[ \sup \left\{ \left| \frac{1}{n_{\bbQ}} \sum_{j=1}^{n_{\bbQ}} \sigma_{j} \varphi \left( X_{j}^{\bbQ} \right) \left\{ \hat{f}_0 \left( X_{j}^{\bbQ} \right) - f^* \left( X_{j}^{\bbQ} \right) \right\} \right| : \varphi \in \calF^* \right\} \right] \\
        = \ & \bbE_{\bfsigma_{1:n_{\bbQ}} \sim \textsf{Unif} \left( \left\{ \pm 1 \right\}^{n_{\bbQ}} \right)} \left[ \sup \left\{ \left| \frac{1}{n_{\bbQ}} \sum_{j=1}^{n_{\bbQ}} \sigma_{j} \alpha_{j}^{\bbQ} \left( t_j \right) \right| : \bft_{1:n_{\bbQ}} \in \calT_{\bbQ} \left( \bfX_{1:n_{\bbQ}}^{\bbQ} \right) \right\} \right] \\
        \leq \ & 2 \left( 1 + C_{\textsf{rf}} \right) \bbE_{\bfsigma_{1:n_{\bbQ}} \sim \textsf{Unif} \left( \left\{ \pm 1 \right\}^{n_{\bbQ}} \right)} \left[ \sup \left\{ \left| \frac{1}{n_{\bbQ}} \sum_{j=1}^{n_{\bbQ}} \sigma_{j} t_j \right| : \bft_{1:n_{\bbQ}} \in \calT_{\bbQ} \left( \bfX_{1:n_{\bbQ}}^{\bbQ} \right) \right\} \right],
    \end{split}
\end{equation}
where $\calT_{\bbQ}: \bbX^{n_{\bbQ}} \to \calP \left( \bbR^{n_{\bbQ}} \right)$ is defined as
\begin{equation}
    \label{eqn:subsec:proof_lemma:control_Q_empirical_process_v1_v7}
    \begin{split}
        \calT_{\bbQ} \left( \bfx_{1:n_{\bbQ}}^{\bbQ} \right) := \left\{ \left( f \left( x_{j}^{\bbQ} \right) - f^* \left( x_{j}^{\bbQ} \right) : j \in \left[ n_{\bbQ} \right] \right) : f \in \calF \right\} \subseteq \bbR^{n_{\bbQ}}.
    \end{split}
\end{equation}
Here, $\calP \left( \bbR^{n_{\bbQ}} \right)$ refers to the power set of the $n_{\bbQ}$-dimensional Euclidean space $\bbR^{n_{\bbQ}}$. By taking two pieces \eqref{eqn:subsec:proof_lemma:control_Q_empirical_process_v1_v1} and \eqref{eqn:subsec:proof_lemma:control_Q_empirical_process_v1_v3} collectively, we reach
\[
    \begin{split}
        &\bbE_{\bfX_{1:n_{\bbQ}}^{\bbQ} \sim \bbQ_{X}^{\otimes n_{\bbQ}}} \left[ \sup \left\{ \bbG_{n_{\bbQ}}^{\bbQ, (1)} (\varphi) \left( \bfX_{1:n_{\bbQ}}^{\bbQ} \right) : \varphi \in \calF^* \cup \left( - \calF^* \right) \right\} \right] \\
        \leq \ & 4 \left( 1 + C_{\textsf{rf}} \right) \bbE_{\left( \bfX_{1:n_{\bbQ}}^{\bbQ}, \bfsigma_{1:n_{\bbQ}} \right) \sim \bbQ_{X}^{\otimes n_{\bbQ}} \otimes \textsf{Unif} \left( \left\{ \pm 1 \right\}^{n_{\bbQ}} \right)} \left[ \sup \left\{ \left| \frac{1}{n_{\bbQ}} \sum_{j=1}^{n_{\bbQ}} \sigma_{j} t_j \right| : \bft_{1:n_{\bbQ}} \in \calT_{\bbQ} \left( \bfX_{1:n_{\bbQ}}^{\bbQ} \right) \right\} \right] \\
        = \ & 4 \left( 1 + C_{\textsf{rf}} \right) \bbE_{\bfX_{1:n_{\bbQ}}^{\bbQ} \sim \bbQ_{X}^{\otimes n_{\bbQ}}} \left[ \widehat{\calR}_{n_{\bbQ}} \left( \calF^* \right) \left( \bfX_{1:n_{\bbQ}}^{\bbQ} \right) \right] \\
        = \ & 4 \left( 1 + C_{\textnormal{rf}} \right) \calR_{n_{\bbQ}}^{\bbQ} \left( \calF^* \right),
    \end{split}
\]
which thus completes the proof of Lemma \ref{lemma:control_Q_empirical_process_v1}.

\subsection{Proof of Lemma \ref{lemma:control_Q_empirical_process_v2}}
\label{subsec:proof_lemma:control_Q_empirical_process_v2}

\indent Similar to the proof of Lemma \ref{lemma:control_P_empirical_process_v2}, we first introduce the function class $\left\{ \theta_{\bbQ}^{(1)} (\varphi) : \varphi \in \calF^* \cup \left( - \calF^* \right) \right\}$, where $\theta_{\bbQ}^{(1)} (\varphi): \bbX \to \bbR$ is a function defined as
\[
    \theta_{\bbQ}^{(1)} (\varphi) (x) := \varphi (x) \left\{ \hat{f}_0 (x) - f^* (x) \right\}, \quad \forall x \in \bbX.
\]
Then, it follows that
\begin{equation}
    \label{eqn:subsec:proof_lemma:control_Q_empirical_process_v2_v1}
    \begin{split}
        &\sup \left\{ \bbG_{n_{\bbQ, (1)}}^{\bbQ} (\varphi) \left( \bfX_{1:n_{\bbQ}}^{\bbQ} \right) : \varphi \in \calF^* \cup \left( - \calF^* \right) \right\} \\
        = \ & \sup \left\{ \frac{1}{n_{\bbQ}} \sum_{j=1}^{n_{\bbQ}} \theta_{\bbQ}^{(1)} (\varphi) \left( X_{j}^{\bbQ} \right) - \bbE_{X \sim \bbQ_{X}} \left[ \theta_{\bbQ}^{(1)} (\varphi) (X) \right] : \varphi \in \calF^* \cup \left( - \calF^* \right) \right\}.
    \end{split}
\end{equation}
At this point, one can easily reveal based on Assumption \ref{assumption:uniform_boundedness} that
\begin{enumerate} [label = (\roman*)]
    \item $\left| \theta_{\bbQ}^{(1)} (\varphi) (x) \right| \leq 2 \left( 1 + C_{\textsf{rf}} \right)$ for every $\left( x, \varphi \right) \in \bbX \times \left\{ \calF^* \cup \left( - \calF^* \right) \right\}$.
    \item $\textsf{Var}_{X \sim \bbQ_{X}} \left[ \theta_{\bbQ}^{(1)} (\varphi) (X) \right] \leq \bbE_{X \sim \bbQ_{X}} \left[ \left\{ \theta_{\bbQ}^{(1)} (\varphi) (X) \right\}^2 \right] \leq 4 \left( 1 + C_{\textsf{rf}} \right)^2$ for all $f \in \calF \cup \left( - \calF \right)$.
\end{enumerate}
The classical Talagrand's concentration inequality (Lemma \ref{lemma:talagrand_concentration}) with parameters $\left( B, v^2 \right) = \left( 2 \left( 1 + C_{\textsf{rf}} \right), 4 \left( 1 + C_{\textsf{rf}} \right)^2 \right)$ together with the equation \eqref{eqn:subsec:proof_lemma:control_Q_empirical_process_v2_v1} tells us that with probability at least $1 - \delta$,
\begin{align}
    \label{eqn:subsec:proof_lemma:control_Q_empirical_process_v2_v2}
    &\sup \left\{ \bbG_{n_{\bbQ}}^{\bbQ, (1)} (\varphi) \left( \bfX_{1:n_{\bbQ}}^{\bbQ} \right) : \varphi \in \calF^* \cup \left( - \calF^* \right) \right\} \nonumber \\
    &- \bbE_{\bfX_{1:n_{\bbQ}}^{\bbQ} \sim \bbQ_{X}^{\otimes n_{\bbQ}}} \left[ \sup \left\{ \bbG_{n_{\bbQ}}^{\bbQ, (1)} (\varphi) \left( \bfX_{1:n_{\bbQ}}^{\bbQ} \right) : \varphi \in \calF^* \cup \left( - \calF^* \right) \right\} \right] \nonumber \\
    \leq \ & \frac{6 \left( 1 + C_{\textsf{rf}} \right)}{n_{\bbQ}} \log \left( \frac{1}{\delta} \right) + 2 \left( 1 + C_{\textsf{rf}} \right) \sqrt{\frac{2 \log \left( \frac{1}{\delta} \right)}{n_{\bbQ}}} \\
    &+ \frac{2}{\sqrt{n_{\bbQ}}} \sqrt{2 \left( 1 + C_{\textsf{rf}} \right) \bbE_{\bfX_{1:n_{\bbQ}}^{\bbQ} \sim \bbQ_{X}^{\otimes n_{\bbQ}}} \left[ \sup \left\{ \bbG_{n_{\bbQ}}^{\bbQ, (1)} (\varphi) \left( \bfX_{1:n_{\bbQ}}^{\bbQ} \right) : \varphi \in \calF^* \cup \left( - \calF^* \right) \right\} \right] \log \left( \frac{1}{\delta} \right)} \nonumber \\
    \stackrel{\textnormal{(a)}}{\leq} \ & \frac{6 \left( 1 + C_{\textsf{rf}} \right)}{n_{\bbQ}} \log \left( \frac{1}{\delta} \right) + 2 \left( 1 + C_{\textsf{rf}} \right) \sqrt{\frac{2 \log \left( \frac{1}{\delta} \right)}{n_{\bbQ}}} + 4 \left( 1 + C_{\textsf{rf}} \right) \sqrt{\frac{2 \calR_{n_{\bbQ}}^{\bbQ_{X}} \left( \calF^* \right) \log \left( \frac{1}{\delta} \right)}{n_{\bbQ}}}, \nonumber
\end{align}
where the step (a) follows by Lemma \ref{lemma:control_Q_empirical_process_v1}, and this finishes the proof of Lemma \ref{lemma:control_Q_empirical_process_v2}.

\subsection{Proof of Lemma \ref{lemma:control_Q_empirical_process_v3}}
\label{subsec:proof_lemma:control_Q_empirical_process_v3}

Similar to the proofs for Lemmas \ref{lemma:control_P_empirical_process_v1} and \ref{lemma:control_Q_empirical_process_v1}, we embark on the proof using the standard symmetrization argument from empirical processes theory, which yields
\begin{equation}
    \label{eqn:subsec:proof_lemma:control_Q_empirical_process_v3_v1}
    \begin{split}
        &\bbE_{\bfX_{1:n_{\bbQ}}^{\bbQ} \sim \bbQ_{X}^{\otimes n_{\bbQ}}} \left[ \sup \left\{ \bbG_{n_{\bbQ}}^{\bbQ, (2)} (\varphi) \left( \bfX_{1:n_{\bbQ}}^{\bbQ} \right) : \varphi \in \calF^* \cup \left( - \calF^* \right) \right\} \right] \\ 
        = \ & \bbE_{\bfX_{1:n_{\bbQ}}^{\bbQ} \sim \bbQ_{X}^{\otimes n_{\bbQ}}} \left[ \sup \left\{ \left| \bbG_{n_{\bbQ}}^{\bbQ, (2)} (\varphi) \left( \bfX_{1:n_{\bbQ}}^{\bbQ} \right) \right| : \varphi \in \calF^* \right\} \right] \\
        \leq \ & 2 \bbE_{\left( \bfX_{1:n_{\bbQ}}^{\bbQ}, \bfsigma_{1:n_{\bbQ}} \right) \sim \bbQ_{X}^{\otimes n_{\bbQ}} \otimes \textsf{Unif} \left( \left\{ \pm 1 \right\}^{n_{\bbQ}} \right)} \left[ \sup \left\{ \left| \frac{1}{n_{\bbQ}} \sum_{j=1}^{n_{\bbQ}} \sigma_{j} \left\{ \varphi \left( X_{j}^{\bbQ} \right) \right\}^2 \right| : \varphi \in \calF^* \right\} \right].
    \end{split}
\end{equation}

\noindent We now consider the function $\beta^{\bbQ}: \left[ -2, 2 \right] \to \bbR$ defined to be
\[
    \beta^{\bbQ} (t) := t^2, \quad \forall t \in \left[ -2, 2 \right].
\]
It turns out that $\beta^{\bbQ} : \left[ -2, 2 \right] \to \bbR$ is a $4$-Lipschitz continuous function with $\beta^{\bbQ} (0) = 0$. Then, the Ledoux-Talagrand contraction principle (Lemma \ref{lemma:ledoux_talagrand_contraction}) tells us that while being conditioned on $\bfX_{1:n_{\bbQ}}^{\bbQ}$,
\begin{equation}
    \label{eqn:subsec:proof_lemma:control_Q_empirical_process_v3_v2}
    \begin{split}
        &\bbE_{\bfsigma_{1:n_{\bbQ}} \sim \textsf{Unif} \left( \left\{ \pm 1 \right\}^{n_{\bbQ}} \right)} \left[ \sup \left\{ \left| \frac{1}{n_{\bbQ}} \sum_{j=1}^{n_{\bbQ}} \sigma_{j} \left\{ \varphi \left( X_{j}^{\bbQ} \right) \right\}^2 \right| : \varphi \in \calF^* \right\} \right] \\
        = \ & \bbE_{\bfsigma_{1:n_{\bbQ}} \sim \textsf{Unif} \left( \left\{ \pm 1 \right\}^{n_{\bbQ}} \right)} \left[ \sup \left\{ \left| \frac{1}{n_{\bbQ}} \sum_{j=1}^{n_{\bbQ}} \sigma_{j} \beta^{\bbQ} \left( t_j \right) \right| : \bft_{1:n_{\bbQ}} \in \calT_{\bbQ} \left( \bfX_{1:n_{\bbQ}}^{\bbQ} \right) \right\} \right] \\
        \leq \ & 8 \bbE_{\bfsigma_{1:n_{\bbQ}} \sim \textsf{Unif} \left( \left\{ \pm 1 \right\}^{n_{\bbQ}} \right)} \left[ \sup \left\{ \left| \frac{1}{n_{\bbQ}} \sum_{j=1}^{n_{\bbQ}} \sigma_{j} t_j \right| : \bft_{1:n_{\bbQ}} \in \calT_{\bbQ} \left( \bfX_{1:n_{\bbQ}}^{\bbQ} \right) \right\} \right],
    \end{split}
\end{equation}
where the function $\calT_{\bbQ}: \bbX^{n_{\bbQ}} \to \calP \left( \bbR^{n_{\bbQ}} \right)$ is previously defined as \eqref{eqn:subsec:proof_lemma:control_Q_empirical_process_v1_v7}. Putting two pieces \eqref{eqn:subsec:proof_lemma:control_Q_empirical_process_v3_v1} and \eqref{eqn:subsec:proof_lemma:control_Q_empirical_process_v3_v2} together, we arrive at
\[
    \begin{split}
        &\bbE_{\bfX_{1:n_{\bbQ}}^{\bbQ} \sim \bbQ_{X}^{\otimes n_{\bbQ}}} \left[ \sup \left\{ \bbG_{n_{\bbQ}}^{\bbQ, (2)} (\varphi) \left( \bfX_{1:n_{\bbQ}}^{\bbQ} \right) : \varphi \in \calF^* \cup \left( - \calF^* \right) \right\} \right] \\
        \leq \ & 16 \cdot \bbE_{\left( \bfX_{1:n_{\bbQ}}^{\bbQ}, \bfsigma_{1:n_{\bbQ}} \right) \sim \bbQ_{X}^{\otimes n_{\bbQ}} \otimes \textsf{Unif} \left( \left\{ \pm 1 \right\}^{n_{\bbQ}} \right)} \left[ \sup \left\{ \left| \frac{1}{n_{\bbQ}} \sum_{j=1}^{n_{\bbQ}} \sigma_{j} t_j \right| : \bft_{1:n_{\bbQ}} \in \calT_{\bbQ} \left( \bfX_{1:n_{\bbQ}}^{\bbQ} \right) \right\} \right] \\
        = \ & 16 \cdot \bbE_{\bfX_{1:n_{\bbQ}}^{\bbQ} \sim \bbQ_{X}^{\otimes n_{\bbQ}}} \left[ \widehat{\calR}_{n_{\bbQ}} \left( \calF^* \right) \left( \bfX_{1:n_{\bbQ}}^{\bbQ} \right) \right] \\
        = \ & 16 \cdot \calR_{n_{\bbQ}}^{\bbQ} \left( \calF^* \right),
    \end{split}
\]
as desired.

\subsection{Proof of Lemma \ref{lemma:control_Q_empirical_process_v4}}
\label{subsec:proof_lemma:control_Q_empirical_process_v4}

\indent We begin the proof by introducing the function class $\left\{ \theta_{\bbQ}^{(2)} (\varphi) : \varphi \in \calF^* \cup \left( - \calF^* \right) \right\}$, where $\theta_{\bbQ}^{(2)} (\varphi): \bbX \to \bbR$ is a function defined as
\[
    \theta_{\bbQ}^{(2)} (\varphi) (x) := \left\{ \varphi (x) \right\}^2, \quad \forall x \in \bbX.
\]
Then, it is obvious that
\begin{equation}
    \label{eqn:subsec:proof_lemma:control_Q_empirical_process_v4_v1}
    \begin{split}
        &\sup \left\{ \bbG_{n_{\bbQ, (2)}}^{\bbQ} (\varphi) \left( \bfX_{1:n_{\bbQ}}^{\bbQ} \right) : \varphi \in \calF^* \cup \left( - \calF^* \right) \right\} \\
        = \ & \sup \left\{ \frac{1}{n_{\bbQ}} \sum_{j=1}^{n_{\bbQ}} \theta_{\bbQ}^{(2)} (\varphi) \left( X_{j}^{\bbQ} \right) - \bbE_{X \sim \bbQ_{X}} \left[ \theta_{\bbQ}^{(2)} (\varphi) (X) \right] : \varphi \in \calF^* \cup \left( - \calF^* \right) \right\}.
    \end{split}
\end{equation}
Also, one can easily find based on Assumption \ref{assumption:uniform_boundedness} that
\begin{enumerate} [label = (\roman*)]
    \item $\left| \theta_{\bbQ}^{(2)} (\varphi) (x) \right| \leq 4$ for every $\left( x, \varphi \right) \in \bbX \times \left\{ \calF^* \cup \left( - \calF^* \right) \right\}$.
    \item $\textsf{Var}_{X \sim \bbQ_{X}} \left[ \theta_{\bbQ}^{(2)} (\varphi) (X) \right] \leq \bbE_{X \sim \bbQ_{X}} \left[ \left\{ \theta_{\bbQ}^{(2)} (\varphi) (X) \right\}^2 \right] \leq 16$ for every $f \in \calF \cup \left( - \calF \right)$.
\end{enumerate}
By virtue of the classical Talagrand's concentration inequality (Lemma \ref{lemma:talagrand_concentration}) with $\left( B, v^2 \right) = \left( 4, 16 \right)$ together with the equation \eqref{eqn:subsec:proof_lemma:control_Q_empirical_process_v4_v1}, we establish with probability at least $1 - \delta$ that
\[
    \begin{split}
        &\sup \left\{ \bbG_{n_{\bbQ}}^{\bbQ, (2)} (\varphi) \left( \bfX_{1:n_{\bbQ}}^{\bbQ} \right) : \varphi \in \calF^* \cup \left( - \calF^* \right) \right\} \\
        &- \bbE_{\bfX_{1:n_{\bbQ}}^{\bbQ} \sim \bbQ_{X}^{\otimes n_{\bbQ}}} \left[ \sup \left\{ \bbG_{n_{\bbQ}}^{\bbQ, (2)} (\varphi) \left( \bfX_{1:n_{\bbQ}}^{\bbQ} \right) : \varphi \in \calF^* \cup \left( - \calF^* \right) \right\} \right] \\
        \leq \ & \frac{12}{n_{\bbQ}} \log \left( \frac{1}{\delta} \right) + 4 \sqrt{\frac{2 \log \left( \frac{1}{\delta} \right)}{n_{\bbQ}}} \\
        &+ \frac{4}{\sqrt{n_{\bbQ}}} \sqrt{\bbE_{\bfX_{1:n_{\bbQ}}^{\bbQ} \sim \bbQ_{X}^{\otimes n_{\bbQ}}} \left[ \sup \left\{ \bbG_{n_{\bbQ}}^{\bbQ, (2)} (\varphi) \left( \bfX_{1:n_{\bbQ}}^{\bbQ} \right) : \varphi \in \calF^* \cup \left( - \calF^* \right) \right\} \right] \log \left( \frac{1}{\delta} \right)} \\
        \stackrel{\textnormal{(a)}}{\leq} \ & \frac{12}{n_{\bbQ}} \log \left( \frac{1}{\delta} \right) + 4 \sqrt{\frac{2 \log \left( \frac{1}{\delta} \right)}{n_{\bbQ}}} + 16 \sqrt{\frac{\calR_{n_{\bbQ}}^{\bbQ_{X}} \left( \calF^* \right) \log \left( \frac{1}{\delta} \right)}{n_{\bbQ}}},
    \end{split}
\]
where the step (a) follows due to Lemma \ref{lemma:control_Q_empirical_process_v3}. This ends the proof of Lemma \ref{lemma:control_Q_empirical_process_v4}.

\section{Proof of auxiliary lemmas for the proof of Theorem \ref{thm:detailed_upper_bound_alg:dr_covariate_shift_adaptation_v2}}
\label{sec:proof_auxiliary_lemmas_subsec:proof_thm:upper_bound_alg:dr_covariate_shift_adaptation_v2}

\subsection{Proof of Lemma \ref{lemma:subsec:proof_thm:upper_bound_alg:dr_covariate_shift_adaptation_v2_v1}}
\label{subsec:proof_lemma:subsec:proof_thm:upper_bound_alg:dr_covariate_shift_adaptation_v2_v1}

\indent Before delving into the proof of our main Lemma (Lemma \ref{lemma:subsec:proof_thm:upper_bound_alg:dr_covariate_shift_adaptation_v2_v1}) that plays a key role in the proof for Theorem \ref{thm:detailed_upper_bound_alg:dr_covariate_shift_adaptation_v2}, we first establish key concentration properties of the gradient and the Hessian matrix of the DR empirical risk \eqref{eqn:dr_empirical_risk_v2} that holds under Assumptions \ref{assumption:covariate_shift}--\ref{assumption:smoothness}, whose proofs are postponed to the final part of this subsection.

\begin{lemma} [Concentration property of the gradient]
\label{lemma:concentration_gradient}
For any $\bfA \in \bbR^{d \times d}$, there exists a universal constant $\calC (A) \in \left( 0, +\infty \right)$ that obeys the following concentration property of the gradient of the DR empirical risk \eqref{eqn:dr_empirical_risk_v2} with respect to the parameter vector $\bftheta$: for any $\delta \in (0, 1]$, it holds that
\begin{equation}
    \label{eqn:lemma:concentration_gradient_v1}
    \begin{split}
        &\left\| \bfA \left\{ \nabla_{\bftheta} \widehat{\calR}_{\textnormal{\textsf{DR}}} \left( \bfO_{1:n_{\bbP}}^{\bbP}, \bfX_{1:n_{\bbQ}}^{\bbQ} \right) \left( \bftheta^* \right) - \bbE_{\left( \bfO_{1:n_{\bbP}}^{\bbP}, \bfX_{1:n_{\bbQ}}^{\bbQ} \right) \sim \bbP^{\otimes n_{\bbP}} \otimes \bbQ_{X}^{n_{\bbQ}}} \left[ \nabla_{\bftheta} \widehat{\calR}_{\textnormal{\textsf{DR}}} \left( \bfO_{1:n_{\bbP}}^{\bbP}, \bfX_{1:n_{\bbQ}}^{\bbQ} \right) \left( \bftheta^* \right) \right] \right\} \right\|_{2} \\
        \leq \ & \calC (\bfA) \left\{ \sqrt{\frac{\calV_{\bbP} (\bfA) \log \left( \frac{2d}{\delta} \right)}{n_{\bbP}}} + \sqrt{\frac{\calV_{\bbQ} (\bfA) \log \left( \frac{2d}{\delta} \right)}{n_{\bbQ}}} \right. \\
        &\left.+ \underbrace{4 \left( 1 + C_{\textnormal{\textsf{dr}}} \right) \left( 1 + C_{\textnormal{\textsf{rf}}} \right) b_1}_{= \ \calB_1 \ \left( \textnormal{defined in \eqref{eqn:key_absolute_constants_thm:upper_bound_alg:dr_covariate_shift_adaptation_v2}} \right)} \left\| \bfA \right\|_{\textnormal{\textnormal{\textsf{op}}}} \cdot \log \left( \frac{2d}{\delta} \right) \left( \frac{1}{n_{\bbP}} + \frac{1}{n_{\bbQ}} \right) \right\}
    \end{split}
\end{equation}
with probability higher than $1 - \delta$ under the data generating process $\left( \bfO_{1:n_{\bbP}}^{\bbP}, \bfX_{1:n_{\bbQ}}^{\bbQ} \right) \sim \bbP^{\otimes n_{\bbP}} \otimes \bbQ_{X}^{n_{\bbQ}}$, where the functions $\Phi_{\bbP} : \bbX \times \bbR \to \bbR^d$ and $\Phi_{\bbQ} : \bbX \to \bbR^d$ are defined as
\begin{equation}
    \label{eqn:lemma:concentration_gradient_v2}
    \begin{split}
        \Phi_{\bbP} (x, y) := \ & 2 \hat{\rho} (x) \left\{ \hat{f}_0 (x) - y \right\} \nabla_{\bftheta} f \left( x; \bftheta^* \right) \quad \textnormal{and} \\
        \Phi_{\bbQ} (x) := \ & 2 \left\{ f \left( x; \bftheta^* \right) - \hat{f}_0 (x) \right\} \nabla_{\bftheta} f \left( x; \bftheta^* \right),
    \end{split}
\end{equation}
respectively, and the quantities $\calV_{\bbP} (\bfA) \in \left( 0, +\infty \right)$ and $\calV_{\bbQ} (\bfA) \in \left( 0, +\infty \right)$ are defined by
\begin{equation}
    \label{eqn:lemma:concentration_gradient_v3}
    \begin{split}
        \calV_{\bbP} (\bfA) := \ & \bbE_{(X, Y) \sim \bbP} \left[ \left\| \bfA \left\{ \Phi_{\bbP} (X, Y) - \bbE_{(X, Y) \sim \bbP} \left[ \Phi_{\bbP} (X, Y) \right] \right\} \right\|_{2}^2 \right] \quad \textnormal{and} \\ \calV_{\bbQ} (\bfA) :=  \ & \bbE_{X \sim \bbQ_{X}} \left[ \left\| \bfA \left\{ \Phi_{\bbQ} (X) - \bbE_{X \sim \bbQ_{X}} \left[ \Phi_{\bbQ} (X) \right] \right\} \right\|_{2}^2 \right],
    \end{split}
\end{equation}
respectively.
\end{lemma}

\begin{lemma} [Concentration property of the Hessian]
\label{lemma:concentration_hessian_matrix}
The Hessian matrix of the DR empirical risk \eqref{eqn:dr_empirical_risk_v2} with respect to the parameter vector $\bftheta$ has the following concentration property: for any $\delta \in (0, 1]$, it holds that
\begin{align}
    \label{eqn:lemma:concentration_hessian_matrix_v1}
        &\left\| \nabla_{\bftheta}^2 \widehat{\calR}_{\textnormal{\textsf{DR}}} \left( \bfO_{1:n_{\bbP}}^{\bbP}, \bfX_{1:n_{\bbQ}}^{\bbQ} \right) \left( \bftheta^* \right) - \bbE_{\left( \bfO_{1:n_{\bbP}}^{\bbP}, \bfX_{1:n_{\bbQ}}^{\bbQ} \right) \sim \bbP^{\otimes n_{\bbP}} \otimes \bbQ_{X}^{n_{\bbQ}}} \left[ \nabla_{\bftheta}^2 \widehat{\calR}_{\textnormal{\textsf{DR}}} \left( \bfO_{1:n_{\bbP}}^{\bbP}, \bfX_{1:n_{\bbQ}}^{\bbQ} \right) \left( \bftheta^* \right) \right] \right\|_{\textnormal{op}} \nonumber \\
        \leq \ & \underbrace{8 \sqrt{2} \cdot \max \left\{ C_{\textnormal{\textsf{dr}}} \left( 1 + C_{\textnormal{\textsf{rf}}} \right) b_2, b_{1}^2 + \left( 1 + C_{\textnormal{\textsf{rf}}} \right) b_2 \right\}}_{= \ \calB_2 \ \left( \textnormal{defined in \eqref{eqn:key_absolute_constants_thm:upper_bound_alg:dr_covariate_shift_adaptation_v2}} \right)} \left\{ \sqrt{\frac{\log \left( \frac{4d}{\delta} \right)}{n_{\bbP}}} + \sqrt{\frac{\log \left( \frac{4d}{\delta} \right)}{n_{\bbQ}}} \right\}
\end{align}
with probability at least $1 - \delta$ under the data generating process $\left( \bfO_{1:n_{\bbP}}^{\bbP}, \bfX_{1:n_{\bbQ}}^{\bbQ} \right) \sim \bbP^{\otimes n_{\bbP}} \otimes \bbQ_{X}^{n_{\bbQ}}$.
\end{lemma}

\noindent Lastly, it is straightforward to see that
\begin{align}
    \label{eqn:derivatives_dr_risk}
    \nabla_{\bftheta} \widehat{\calR}_{\textsf{DR}} \left( \bfo_{1:n_{\bbP}}^{\bbP}, \bfx_{1:n_{\bbQ}}^{\bbQ} \right) \left( \bftheta \right) = \ & \frac{2}{n_{\bbP}} \sum_{i=1}^{n_{\bbP}} \hat{\rho} \left( x_{i}^{\bbP} \right) \left\{ \hat{f}_0 \left( x_{i}^{\bbP} \right) - y_{i}^{\bbP} \right\} \nabla_{\bftheta} f \left( x_{i}^{\bbP}; \bftheta \right) \nonumber \\
    &+ \frac{2}{n_{\bbQ}} \sum_{j=1}^{n_{\bbQ}} \left\{ f \left( x_{j}^{\bbQ}; \bftheta \right) - \hat{f}_0 \left( x_{j}^{\bbQ} \right) \right\} \nabla_{\bftheta} f \left( x_{j}^{\bbQ}; \bftheta \right), \nonumber \\
    \nabla_{\bftheta}^2 \widehat{\calR}_{\textsf{DR}} \left( \bfo_{1:n_{\bbP}}^{\bbP}, \bfx_{1:n_{\bbQ}}^{\bbQ} \right) \left( \bftheta \right) = \ & \frac{2}{n_{\bbP}} \sum_{i=1}^{n_{\bbP}} \hat{\rho} \left( x_{i}^{\bbP} \right) \left\{ \hat{f}_0 \left( x_{i}^{\bbP} \right) - y_{i}^{\bbP} \right\} \nabla_{\bftheta}^2 f \left( x_{i}^{\bbP}; \bftheta \right) \\
    &+ \frac{2}{n_{\bbQ}} \sum_{j=1}^{n_{\bbQ}} \left[ \nabla_{\bftheta} f \left( x_{j}^{\bbQ}; \bftheta \right) \left\{ \nabla_{\bftheta} f \left( x_{j}^{\bbQ}; \bftheta \right) \right\}^{\top} + \left\{ f \left( x_{j}^{\bbQ}; \bftheta \right) - \hat{f}_0 \left( x_{j}^{\bbQ} \right) \right\} \nabla_{\bftheta}^2 f \left( x_{j}^{\bbQ}; \bftheta \right) \right], \nonumber \\
    \nabla_{\bftheta}^3 \widehat{\calR}_{\textsf{DR}} \left( \bfo_{1:n_{\bbP}}^{\bbP}, \bfx_{1:n_{\bbQ}}^{\bbQ} \right) \left( \bftheta \right) = \ & \frac{2}{n_{\bbP}} \sum_{i=1}^{n_{\bbP}} \hat{\rho} \left( x_{i}^{\bbP} \right) \left\{ \hat{f}_0 \left( x_{i}^{\bbP} \right) - y_{i}^{\bbP} \right\} \nabla_{\bftheta}^3 f \left( x_{i}^{\bbP}; \bftheta \right) \nonumber \\
    &+ \frac{2}{n_{\bbQ}} \sum_{j=1}^{n_{\bbQ}} \left[ \nabla_{\bftheta}^2 f \left( x_{j}^{\bbQ}; \bftheta \right) \otimes \nabla_{\bftheta} f \left( x_{j}^{\bbQ}; \bftheta \right) + \nabla_{\bftheta}^2 f \left( x_{j}^{\bbQ}; \bftheta \right) \otimes \nabla_{\bftheta} f \left( x_{j}^{\bbQ}; \bftheta \right) \right. \nonumber \\
    &\left. + \left\{ f \left( x_{j}^{\bbQ}; \bftheta \right) - \hat{f}_0 \left( x_{j}^{\bbQ} \right) \right\} \nabla_{\bftheta}^3 f \left( x_{j}^{\bbQ}; \bftheta \right) \right]. \nonumber
\end{align}
By making use of the observation \eqref{eqn:derivatives_dr_risk}, it follows that
\begin{equation}
    \label{eqn:expectation_derivatives_dr_risk}
    \begin{split}
        \bbE \left[ \nabla_{\bftheta} \widehat{\calR}_{\textsf{DR}} \left( \bfO_{1:n_{\bbP}}^{\bbP}, \bfX_{1:n_{\bbQ}}^{\bbQ} \right) \left( \bftheta^* \right) \right] = \ & 2 \bbE_{X \sim \bbP_{X}} \left[ \left\{ \hat{\rho} (X) - \rho^* (X) \right\} \left\{ \hat{f}_0 (X) - f \left( X; \bftheta^* \right) \right\} \nabla_{\bftheta} f \left( X; \bftheta^* \right) \right], \\
        \bbE \left[ \nabla_{\bftheta}^2 \widehat{\calR}_{\textsf{DR}} \left( \bfO_{1:n_{\bbP}}^{\bbP}, \bfX_{1:n_{\bbQ}}^{\bbQ} \right) \left( \bftheta^* \right) \right] = \ & 2 \bbE_{X \sim \bbP_{X}} \left[ \left\{ \hat{\rho} (X) - \rho^* (X) \right\} \left\{ \hat{f}_0 (X) - f \left( X; \bftheta^* \right) \right\} \nabla_{\bftheta}^2 f \left( X; \bftheta^* \right) \right] \\
        &+ \calI_{\bbQ} \left( \bftheta^* \right),
    \end{split}
\end{equation}
and while being conditioned on the event $\calA_{\bbP} := \bigcap_{i=1}^{n_{\bbP}} \left\{ \left| Y_{i}^{\bbP} \right| \leq 1 \right\}$ that
\begin{equation}
    \label{eqn:bound_operator_norm_third_order_derivative_dr_risk}
    \begin{split}
        \left\| \nabla_{\bftheta}^3 \widehat{\calR}_{\textsf{DR}} \left( \bfO_{1:n_{\bbP}}^{\bbP}, \bfX_{1:n_{\bbQ}}^{\bbQ} \right) \left( \bftheta \right) \right\|_{\textsf{op}} \leq \ & \frac{2}{n_{\bbP}} \sum_{i=1}^{n_{\bbP}} \left\| \hat{\rho} \right\|_{\infty} \left( 1 + \left\| \hat{f}_0 \right\|_{\infty} \right) \left\| \nabla_{\bftheta}^2 f \left( X_{i}^{\bbP}; \bftheta \right) \right\|_{\textsf{op}} \\
        &+ \frac{2}{n_{\bbQ}} \sum_{i=1}^{n_{\bbQ}} \left\{ 3 \left\| \nabla_{\bftheta} f \left( X_{j}^{\bbQ}; \bftheta \right) \right\|_{2} \left\| \nabla_{\bftheta}^2 f \left( X_{j}^{\bbQ}; \bftheta \right) \right\|_{\textsf{op}} \right. \\
        &\left. + \left( 1 + \left\| \hat{f}_0 \right\|_{\infty} \right) \left\| \nabla_{\bftheta}^3 f \left( X_{j}^{\bbQ}; \bftheta \right) \right\|_{\textsf{op}} \right\} \\
        \stackrel{\textnormal{(a)}}{\leq} \ & 2 \left( 1 + C_{\textsf{dr}} \right) \left( 1 + C_{\textsf{rf}} \right) b_3 + 6 b_1 b_2 \leq \calB_3
    \end{split}
\end{equation}
for every $\bftheta \in \Theta$, where the step (a) comes from Assumption \ref{assumption:black_box_estimates} and the part (\romannumeral 2) of Assumption \ref{assumption:smoothness}.
\medskip

\indent With these preliminary results in our hand, we are now ready to prove Lemma \ref{lemma:subsec:proof_thm:upper_bound_alg:dr_covariate_shift_adaptation_v2_v1}. Hereafter, we focus on the case where $\Theta = \bbR^d$ for simplicity of presentation. Given any $\delta \in (0, 1)$ and any fixed matrix $\bfA \in \bbR^{d \times d}$, we define the events
\begin{align}
    \label{eqn:subsec:proof_lemma:subsec:proof_thm:upper_bound_alg:dr_covariate_shift_adaptation_v2_v1_v1}
        &\calE_1 \left( \delta; \bfA \right) := \left\{ \left( \bfo_{1:n_{\bbP}}^{\bbP}, \bfx_{1:n_{\bbQ}}^{\bbQ} \right) \in \bbO^{n_{\bbP}} \times \bbX^{n_{\bbQ}} : \right. \nonumber \\
        & \left\| \bfA \left\{ \nabla_{\bftheta} \widehat{\calR}_{\textsf{DR}} \left( \bfo_{1:n_{\bbP}}^{\bbP}, \bfx_{1:n_{\bbQ}}^{\bbQ} \right) \left( \bftheta^* \right) - \bbE \left[ \nabla_{\bftheta} \widehat{\calR}_{\textsf{DR}} \left( \bfO_{1:n_{\bbP}}^{\bbP}, \bfX_{1:n_{\bbQ}}^{\bbQ} \right) \left( \bftheta^* \right) \right] \right\} \right\|_{2} \\
        \leq \ & \left. \calC (\bfA) \left\{ \sqrt{\frac{\calV_{\bbP} (\bfA) \log \left( \frac{d}{\delta} \right)}{n_{\bbP}}} + \sqrt{\frac{\calV_{\bbQ} (\bfA) \log \left( \frac{d}{\delta} \right)}{n_{\bbQ}}} + \calB_1 \left\| \bfA \right\|_{\textsf{op}} \cdot \log \left( \frac{d}{\delta} \right) \left( \frac{1}{n_{\bbP}} + \frac{1}{n_{\bbQ}} \right) \right\} \right\}, \nonumber
\end{align}
and
\begin{equation}
    \label{eqn:subsec:proof_lemma:subsec:proof_thm:upper_bound_alg:dr_covariate_shift_adaptation_v2_v1_v2}
    \begin{split}
        &\calE_2 (\delta) := \left\{ \left( \bfo_{1:n_{\bbP}}^{\bbP}, \bfx_{1:n_{\bbQ}}^{\bbQ} \right) \in \bbO^{n_{\bbP}} \times \bbX^{n_{\bbQ}} : \right. \\
        &\left\| \nabla_{\bftheta}^2 \widehat{\calR}_{\textsf{DR}} \left( \bfO_{1:n_{\bbP}}^{\bbP}, \bfX_{1:n_{\bbQ}}^{\bbQ} \right) \left( \bftheta^* \right) - \bbE \left[ \nabla_{\bftheta}^2 \widehat{\calR}_{\textsf{DR}} \left( \bfO_{1:n_{\bbP}}^{\bbP}, \bfX_{1:n_{\bbQ}}^{\bbQ} \right) \left( \bftheta^* \right) \right] \right\|_{\textsf{op}} \\
        \leq \ & \left. \calB_2 \left\{ \sqrt{\frac{\log \left( \frac{d}{\delta} \right)}{n_{\bbP}}} + \sqrt{\frac{\log \left( \frac{d}{\delta} \right)}{n_{\bbQ}}} \right\} \right\},
    \end{split}
\end{equation}
so that we have $\left( \bbP^{\otimes n_{\bbP}} \otimes \bbQ_{X}^{\otimes n_{\bbQ}} \right) \left\{ \calE_1 \left( \delta; \bfA \right) \right\} \geq 1 - 2 \delta$ and $\left( \bbP^{\otimes n_{\bbP}} \otimes \bbQ_{X}^{\otimes n_{\bbQ}} \right) \left\{ \calE_2 (\delta) \right\} \geq 1 - 4 \delta$ for any $\left( \bfA, \delta \right) \in \bbR^{d \times d} \times \left( 0, \frac{1}{4} \right]$ due to Lemma \ref{lemma:concentration_gradient} and \ref{lemma:concentration_hessian_matrix}. For simplicity, we employ the notation $\widehat{\calR}_{\textsf{DR}} := \widehat{\calR}_{\textsf{DR}} \left( \bfO_{1:n_{\bbP}}^{\bbP}, \bfX_{1:n_{\bbQ}}^{\bbQ} \right) : \bbR^d \to \bbR$ for $\left( \bfO_{1:n_{\bbP}}^{\bbP}, \bfX_{1:n_{\bbQ}}^{\bbQ} \right) \sim \bbP^{\otimes n_{\bbP}} \otimes \bbQ_{X}^{\otimes n_{\bbQ}}$ as well as $\bfg := \nabla_{\bftheta} \widehat{\calR}_{\textsf{DR}} \left( \bftheta^* \right) - \bbE \left[ \nabla_{\bftheta} \widehat{\calR}_{\textsf{DR}} \left( \bftheta^* \right) \right]$ throughout this subsection. Owing to Assumption \ref{assumption:smoothness}, it turns out for every $\bftheta \in \bbR^d$ that while being on the event $\calA_{\bbP} \cap \calE_{2} (\delta)$, where $\calA_{\bbP} = \bigcap_{i=1}^{n_{\bbP}} \left\{ \left| Y_{i}^{\bbP} \right| \leq 1 \right\}$,
\begin{align}
    \label{eqn:subsec:proof_lemma:subsec:proof_thm:upper_bound_alg:dr_covariate_shift_adaptation_v2_v1_v3}
        &\widehat{\calR}_{\textsf{DR}} (\bftheta) - \widehat{\calR}_{\textsf{DR}} \left( \bftheta^* \right) \nonumber \\
        \stackrel{\textnormal{(a)}}{\leq} \ & \left( \bftheta - \bftheta^* \right)^{\top} \nabla_{\bftheta} \widehat{\calR}_{\textsf{DR}} \left( \bftheta^* \right) + \frac{1}{2!} \left( \bftheta - \bftheta^* \right)^{\top} \nabla_{\bftheta}^2 \widehat{\calR}_{\textsf{DR}} \left( \bftheta^* \right) \left( \bftheta - \bftheta^* \right) + \frac{\calB_3}{3!} \left\| \bftheta - \bftheta^* \right\|_{2}^3 \nonumber \\
        \stackrel{\textnormal{(b)}}{\leq} \ & \left( \bftheta - \bftheta^* \right)^{\top} \bbE \left[ \nabla_{\bftheta} \widehat{\calR}_{\textsf{DR}} \left( \bftheta^* \right) \right] + \left( \bftheta - \bftheta^* \right)^{\top} \bfg + \frac{1}{2} \left( \bftheta - \bftheta^* \right)^{\top} \bbE \left[ \nabla_{\bftheta}^2 \widehat{\calR}_{\textsf{DR}} \left( \bftheta^* \right) \right] \left( \bftheta - \bftheta^* \right) \nonumber \\
        &+ \frac{\calB_2}{2} \sqrt{\log \left( \frac{d}{\delta} \right)} \left( \frac{1}{\sqrt{n_{\bbP}}} + \frac{1}{\sqrt{n_{\bbQ}}} \right) \left\| \bftheta - \bftheta^* \right\|_{2}^2 + \frac{\calB_3}{6} \left\| \bftheta - \bftheta^* \right\|_{2}^3 \\
        \stackrel{\textnormal{(c)}}{=} \ & 2 \bbE_{X \sim \bbP_{X}} \left[ \left\{ \hat{\rho} (X) - \rho^* (X) \right\} \left\{ \hat{f}_0 (X) - f \left( X; \bftheta^* \right) \right\} \left( \bftheta - \bftheta^* \right)^{\top} \nabla_{\bftheta} f \left( X; \bftheta^* \right) \right] + \left( \bftheta - \bftheta^* \right)^{\top} \bfg \nonumber \\
        &+ \left( \bftheta - \bftheta^* \right)^{\top} \bbE_{X \sim \bbP_{X}} \left[ \left\{ \hat{\rho} (X) - \rho^* (X) \right\} \left\{ \hat{f}_0 (X) - f \left( X; \bftheta^* \right) \right\} \nabla_{\bftheta}^2 f \left( X; \bftheta^* \right) \right] \left( \bftheta - \bftheta^* \right) \nonumber \\
        &+ \frac{1}{2} \left( \bftheta - \bftheta^* \right)^{\top} \calI_{\bbQ} \left( \bftheta^* \right) \left( \bftheta - \bftheta^* \right) + \frac{\calB_2}{2} \sqrt{\log \left( \frac{d}{\delta} \right)} \left( \frac{1}{\sqrt{n_{\bbP}}} + \frac{1}{\sqrt{n_{\bbQ}}} \right) \left\| \bftheta - \bftheta^* \right\|_{2}^2 + \frac{\calB_3}{6} \left\| \bftheta - \bftheta^* \right\|_{2}^3, \nonumber
\end{align}
where the step (a) follows due to Taylor's theorem together with the fact \eqref{eqn:bound_operator_norm_third_order_derivative_dr_risk}, the step (b) invokes Lemma \ref{lemma:concentration_hessian_matrix}, and the step (c) holds by the observation \eqref{eqn:expectation_derivatives_dr_risk}. Thus, by letting $\Delta (\bftheta) := \bftheta - \bftheta^* \in \bbR^d$, we obtain from the inequality \eqref{eqn:subsec:proof_lemma:subsec:proof_thm:upper_bound_alg:dr_covariate_shift_adaptation_v2_v1_v3} that on the event $\calA_{\bbP} \cap \calE_{2} (\delta)$,
\begin{align}
    \label{eqn:subsec:proof_lemma:subsec:proof_thm:upper_bound_alg:dr_covariate_shift_adaptation_v2_v1_v4}
        &\widehat{\calR}_{\textsf{DR}} (\bftheta) - \widehat{\calR}_{\textsf{DR}} \left( \bftheta^* \right) \nonumber \\
        \leq \ & 2 \bbE_{X \sim \bbP_{X}} \left[ \left\{ \hat{\rho} (X) - \rho^* (X) \right\} \left\{ \hat{f}_0 (X) - f \left( X; \bftheta^* \right) \right\} \left\{ \Delta (\bftheta) \right\}^{\top} \nabla_{\bftheta} f \left( X; \bftheta^* \right) \right] + \left\{ \Delta (\bftheta) \right\}^{\top} \bfg \nonumber \\
        &+ \bbE_{X \sim \bbP_{X}} \left[ \left\{ \hat{\rho} (X) - \rho^* (X) \right\} \left\{ \hat{f}_0 (X) - f \left( X; \bftheta^* \right) \right\} \left\{ \Delta (\bftheta) \right\}^{\top} \nabla_{\bftheta}^2 f \left( X; \bftheta^* \right) \Delta (\bftheta) \right] \\
        &+ \frac{1}{2} \left\{ \Delta (\bftheta) \right\}^{\top} \calI_{\bbQ} \left( \bftheta^* \right) \Delta (\bftheta) + \frac{\calB_2}{2} \sqrt{\log \left( \frac{d}{\delta} \right)} \left( \frac{1}{\sqrt{n_{\bbP}}} + \frac{1}{\sqrt{n_{\bbQ}}} \right) \left\| \Delta (\bftheta) \right\|_{2}^2 + \frac{\calB_3}{6} \left\| \Delta (\bftheta) \right\|_{2}^3 \nonumber \\
        = \ & \bbE_{X \sim \bbP_{X}} \left[ \left\{ \hat{\rho} (X) - \rho^* (X) \right\} \left\{ \hat{f}_0 (X) - f \left( X; \bftheta^* \right) \right\} \left[ 2 \left\{ \Delta (\bftheta) \right\}^{\top} \nabla_{\bftheta} f \left( X; \bftheta^* \right) + \left\{ \Delta (\bftheta) \right\}^{\top} \nabla_{\bftheta}^2 f \left( X; \bftheta^* \right) \Delta (\bftheta) \right] \right] \nonumber \\
        &+ \frac{1}{2} \left\{ \Delta (\bftheta) - \bfz \right\}^{\top} \calI_{\bbQ} \left( \bftheta^* \right) \left\{ \Delta (\bftheta) - \bfz \right\} - \frac{1}{2} \bfz^{\top} \calI_{\bbQ} \left( \bftheta^* \right) \bfz \nonumber \\
        &+ \frac{\calB_2}{2} \sqrt{\log \left( \frac{d}{\delta} \right)} \left( \frac{1}{\sqrt{n_{\bbP}}} + \frac{1}{\sqrt{n_{\bbQ}}} \right) \left\| \Delta (\bftheta) \right\|_{2}^2 + \frac{\calB_3}{6} \left\| \Delta (\bftheta) \right\|_{2}^3, \nonumber
\end{align}
where $\bfz := - \calI_{\bbQ}^{-1} \left( \bftheta^* \right) \bfg \in \bbR^d$. By employing a similar argument, we reveal for any $\bftheta \in \bbR^d$ that while being conditioned on the event $\calA_{\bbP} \cap \calE_{2} (\delta)$, we have
\begin{align}
    \label{eqn:subsec:proof_lemma:subsec:proof_thm:upper_bound_alg:dr_covariate_shift_adaptation_v2_v1_v5}
        &\widehat{\calR}_{\textsf{DR}} (\bftheta) - \widehat{\calR}_{\textsf{DR}} \left( \bftheta^* \right) \nonumber \\
        \stackrel{\textnormal{(d)}}{\geq} \ & \left( \bftheta - \bftheta^* \right)^{\top} \nabla_{\bftheta} \widehat{\calR}_{\textsf{DR}} \left( \bftheta^* \right) + \frac{1}{2!} \left( \bftheta - \bftheta^* \right)^{\top} \nabla_{\bftheta}^2 \widehat{\calR}_{\textsf{DR}} \left( \bftheta^* \right) \left( \bftheta - \bftheta^* \right) - \frac{\calB_3}{3!} \left\| \bftheta - \bftheta^* \right\|_{2}^3 \nonumber \\
        \stackrel{\textnormal{(e)}}{\geq} \ & 2 \bbE_{X \sim \bbP_{X}} \left[ \left\{ \hat{\rho} (X) - \rho^* (X) \right\} \left\{ \hat{f}_0 (X) - f \left( X; \bftheta^* \right) \right\} \left( \bftheta - \bftheta^* \right)^{\top} \nabla_{\bftheta} f \left( X; \bftheta^* \right) \right] + \left( \bftheta - \bftheta^* \right)^{\top} \bfg \nonumber \\
        &+ \left( \bftheta - \bftheta^* \right)^{\top} \bbE_{X \sim \bbP_{X}} \left[ \left\{ \hat{\rho} (X) - \rho^* (X) \right\} \left\{ \hat{f}_0 (X) - f \left( X; \bftheta^* \right) \right\} \nabla_{\bftheta}^2 f \left( X; \bftheta^* \right) \right] \left( \bftheta - \bftheta^* \right) \\
        &+ \frac{1}{2} \left( \bftheta - \bftheta^* \right)^{\top} \calI_{\bbQ} \left( \bftheta^* \right) \left( \bftheta - \bftheta^* \right) - \frac{\calB_2}{2} \sqrt{\log \left( \frac{d}{\delta} \right)} \left( \frac{1}{\sqrt{n_{\bbP}}} + \frac{1}{\sqrt{n_{\bbQ}}} \right) \left\| \bftheta - \bftheta^* \right\|_{2}^2 - \frac{\calB_3}{6} \left\| \bftheta - \bftheta^* \right\|_{2}^3 \nonumber \\
        = \ & \bbE_{X \sim \bbP_{X}} \left[ \left\{ \hat{\rho} (X) - \rho^* (X) \right\} \left\{ \hat{f}_0 (X) - f \left( X; \bftheta^* \right) \right\} \left[ 2 \left\{ \Delta (\bftheta) \right\}^{\top} \nabla_{\bftheta} f \left( X; \bftheta^* \right) + \left\{ \Delta (\bftheta) \right\}^{\top} \nabla_{\bftheta}^2 f \left( X; \bftheta^* \right) \Delta (\bftheta) \right] \right] \nonumber \\
        &+ \frac{1}{2} \left\{ \Delta (\bftheta) - \bfz \right\}^{\top} \calI_{\bbQ} \left( \bftheta^* \right) \left\{ \Delta (\bftheta) - \bfz \right\} - \frac{1}{2} \bfz^{\top} \calI_{\bbQ} \left( \bftheta^* \right) \bfz \nonumber \\
        &- \frac{\calB_2}{2} \sqrt{\log \left( \frac{d}{\delta} \right)} \left( \frac{1}{\sqrt{n_{\bbP}}} + \frac{1}{\sqrt{n_{\bbQ}}} \right) \left\| \Delta (\bftheta) \right\|_{2}^2 - \frac{\calB_3}{6} \left\| \Delta (\bftheta) \right\|_{2}^3, \nonumber
\end{align}
where the step (d) invokes Taylor's theorem as well as the observation \eqref{eqn:bound_operator_norm_third_order_derivative_dr_risk}, and the step (e) follows due to Lemma \ref{lemma:concentration_hessian_matrix} and the fact \eqref{eqn:expectation_derivatives_dr_risk}.
\medskip

\indent Now, we leverage Lemma \ref{lemma:concentration_gradient} with $\bfA = \calI_{\bbQ}^{-1} \left( \bftheta^* \right)$. While being conditioned on the event $\calE_1 \left( \delta; \calI_{\bbQ}^{-1} \left( \bftheta^* \right) \right)$, we have
\begin{equation}
    \label{eqn:subsec:proof_lemma:subsec:proof_thm:upper_bound_alg:dr_covariate_shift_adaptation_v2_v1_v6}
    \begin{split}
        &\left\| \bfz \right\|_{2} = \left\| \calI_{\bbQ}^{-1} \left( \bftheta^* \right) \bfg \right\|_{2} \\
        \leq \ & \calC \left\{ \calI_{\bbQ}^{-1} \left( \bftheta^* \right) \right\} \left[ \sqrt{\frac{\calV_{\bbP} \left\{ \calI_{\bbQ}^{-1} \left( \bftheta^* \right) \right\} \log \left( \frac{d}{\delta} \right)}{n_{\bbP}}} + \sqrt{\frac{\calV_{\bbQ} \left\{ \calI_{\bbQ}^{-1} \left( \bftheta^* \right) \right\} \log \left( \frac{d}{\delta} \right)}{n_{\bbQ}}} \right. \\
        &\left.+ \calB_1 \left\| \calI_{\bbQ}^{-1} \left( \bftheta^* \right) \right\|_{\textsf{op}} \log \left( \frac{d}{\delta} \right) \left( \frac{1}{n_{\bbP}} + \frac{1}{n_{\bbQ}} \right) \right] \\
        \stackrel{\textnormal{(f)}}{\leq} \ & \calC \left\{ \calI_{\bbQ}^{-1} \left( \bftheta^* \right) \right\} \left[ \sqrt{2} \left( 1 + C_{\textsf{dr}} \right) \left( 1 + C_{\textsf{rf}} \right) \sqrt{\log \left( \frac{d}{\delta} \right)} \right. \\
        &\cdot \left\{ \sqrt{\frac{\textsf{Trace} \left\{ \calI_{\bbP} \left( \bftheta^* \right) \calI_{\bbQ}^{-2} \left( \bftheta^* \right) \right\}}{n_{\bbP}}} + \sqrt{\frac{\textsf{Trace} \left\{ \calI_{\bbQ}^{-1} \left( \bftheta^* \right) \right\}}{n_{\bbQ}}} \right\} \\
        &\left. + \ \calB_1 \left\| \calI_{\bbQ}^{-1} \left( \bftheta^* \right) \right\|_{\textsf{op}} \log \left( \frac{d}{\delta} \right) \left( \frac{1}{n_{\bbP}} + \frac{1}{n_{\bbQ}} \right) \right],
    \end{split}
\end{equation}
where the step (f) can be obtained by letting $\bfA = \calI_{\bbQ}^{-1} \left( \bftheta^* \right)$ in the following facts: for any given fixed matrix $\bfA \in \bbR^{d \times d}$, it holds that
\begin{equation}
    \label{eqn:subsec:proof_lemma:subsec:proof_thm:upper_bound_alg:dr_covariate_shift_adaptation_v2_v1_v7}
    \begin{split}
        &\calV_{\bbP} (\bfA) := \bbE_{(X, Y) \sim \bbP} \left[ \left\| \bfA \left\{ \Phi_{\bbP} (X, Y) - \bbE_{(X, Y) \sim \bbP} \left[ \Phi_{\bbP} (X, Y) \right] \right\} \right\|_{2}^2 \right] \\
        = \ & \textsf{Trace} \left\{ \bfA \textsf{Cov}_{(X, Y) \sim \bbP} \left[ \Phi_{\bbP} (X, Y) \right] \bfA^{\top} \right\} \\
        \leq \ & \textsf{Trace} \left\{ \bfA \bbE_{(X, Y) \sim \bbP} \left[ \Phi_{\bbP} (X, Y) \left\{ \Phi_{\bbP} (X, Y) \right\}^{\top} \right] \bfA^{\top} \right\} \\
        = \ & 4 \cdot \textsf{Trace} \left\{ \bfA \bbE_{(X, Y) \sim \bbP} \left[ \hat{\rho}^2 (X) \left\{ \hat{f}_0 (X) - Y \right\}^2 \nabla_{\bftheta} f \left( X ; \bftheta^* \right) \left\{ \nabla_{\bftheta} f \left( X ; \bftheta^* \right) \right\}^{\top} \right] \bfA^{\top} \right\} \\
        \stackrel{\textnormal{(g)}}{\leq} \ & 2 C_{\textsf{dr}}^2 \left( 1 + C_{\textsf{rf}} \right)^2 \textsf{Trace} \left\{ \bfA \underbrace{\bbE_{X \sim \bbP_{X}} \left[ 2 \nabla_{\bftheta} f \left( X ; \bftheta^* \right) \left\{ \nabla_{\bftheta} f \left( X ; \bftheta^* \right) \right\}^{\top} \right]}_{= \ \calI_{\bbP} \left( \bftheta^* \right)} \bfA^{\top} \right\} \\ 
        = \ & 2 C_{\textsf{dr}}^2 \left( 1 + C_{\textsf{rf}} \right)^2 \textsf{Trace} \left\{ \bfA \calI_{\bbP} \left( \bftheta^* \right) \bfA^{\top} \right\},
    \end{split}
\end{equation}
and
\begin{equation}
    \label{eqn:subsec:proof_lemma:subsec:proof_thm:upper_bound_alg:dr_covariate_shift_adaptation_v2_v1_v8}
    \begin{split}
        &\calV_{\bbQ} (\bfA) := \bbE_{X \sim \bbQ_{X}} \left[ \left\| \bfA \left\{ \Phi_{\bbQ} (X) - \bbE_{X \sim \bbQ_{X}} \left[ \Phi_{\bbQ} (X) \right] \right\} \right\|_{2}^2 \right] \\
        = \ & \textsf{Trace} \left\{ \bfA \textsf{Cov}_{X \sim \bbQ_{X}} \left[ \Phi_{\bbQ} (X) \right] \bfA^{\top} \right\} \\
        \leq \ & \textsf{Trace} \left\{ \bfA \bbE_{X \sim \bbQ_{X}} \left[ \Phi_{\bbQ} (X) \left\{ \Phi_{\bbQ} (X) \right\}^{\top} \right] \bfA^{\top} \right\} \\
        = \ & 4 \cdot \textsf{Trace} \left\{ \bfA \bbE_{X \sim \bbQ_{X}} \left[ \left\{ f \left( X ; \bftheta^* \right) - \hat{f}_0 (X) \right\}^2 \nabla_{\bftheta} f \left( X ; \bftheta^* \right) \left\{ \nabla_{\bftheta} f \left( X ; \bftheta^* \right) \right\}^{\top} \right] \bfA^{\top} \right\} \\
        \stackrel{\textnormal{(h)}}{\leq} \ & 2 \left( 1 + C_{\textsf{rf}} \right)^2 \textsf{Trace} \left\{ \bfA \underbrace{\bbE_{X \sim \bbQ_{X}} \left[ 2 \nabla_{\bftheta} f \left( X ; \bftheta^* \right) \left\{ \nabla_{\bftheta} f \left( X ; \bftheta^* \right) \right\}^{\top} \right]}_{= \ \calI_{\bbQ} \left( \bftheta^* \right)} \bfA^{\top} \right\} \\ 
        = \ & 2 \left( 1 + C_{\textsf{rf}} \right)^2 \textsf{Trace} \left\{ \bfA \calI_{\bbQ} \left( \bftheta^* \right) \bfA^{\top} \right\},
    \end{split}
\end{equation}
where the steps (g) and (h) follow by Assumption \ref{assumption:black_box_estimates} on the pilot black-box ML estimates $\hat{\rho}: \bbX \to \bbR_{+}$ and $\hat{f}_0 : \bbX \to \bbR$. Now, we let $K := \max \left\{ \calC \left\{ \calI_{\bbQ}^{-1} \left( \bftheta^* \right) \right\}, \calC \left\{ \calI_{\bbQ}^{-\frac{1}{2}} \left( \bftheta^* \right) \right\} \right\} \in \left( 0, +\infty \right)$. By noting that $\Delta \left( \bftheta^* + \bfz \right) = \bfz = - \calI_{\bbQ}^{-1} \left( \bftheta^* \right) \bfg$, it follows from the inequality \eqref{eqn:subsec:proof_lemma:subsec:proof_thm:upper_bound_alg:dr_covariate_shift_adaptation_v2_v1_v4} that on the event $\calA_{\bbP} \cap \calE_1 \left( \delta; \calI_{\bbQ}^{-1} \left( \bftheta^* \right) \right) \cap \calE_2 (\delta)$,
\begin{align}
    \label{eqn:subsec:proof_lemma:subsec:proof_thm:upper_bound_alg:dr_covariate_shift_adaptation_v2_v1_v9}
        &\widehat{\calR}_{\textsf{DR}} \left( \bftheta^* + \bfz \right) - \widehat{\calR}_{\textsf{DR}} \left( \bftheta^* \right) \nonumber \\
        \leq \ & \bbE_{X \sim \bbP_{X}} \left[ \left\{ \hat{\rho} (X) - \rho^* (X) \right\} \left\{ \hat{f}_0 (X) - f \left( X; \bftheta^* \right) \right\} \left[ 2 \left\{ \Delta (\bftheta) \right\}^{\top} \nabla_{\bftheta} f \left( X; \bftheta^* \right) + \left\{ \Delta (\bftheta) \right\}^{\top} \nabla_{\bftheta}^2 f \left( X; \bftheta^* \right) \Delta (\bftheta) \right] \right] \nonumber \\
        &- \frac{1}{2} \bfz^{\top} \calI_{\bbQ} \left( \bftheta^* \right) \bfz + \frac{\calB_2}{2} \sqrt{\log \left( \frac{d}{\delta} \right)} \left( \frac{1}{\sqrt{n_{\bbP}}} + \frac{1}{\sqrt{n_{\bbQ}}} \right) \left\| \bfz \right\|_{2}^2 + \frac{\calB_3}{6} \left\| \bfz \right\|_{2}^3 \nonumber \\
        \stackrel{\textnormal{(i)}}{\leq} \ & \bbE_{X \sim \bbP_{X}} \left[ \left\{ \hat{\rho} (X) - \rho^* (X) \right\} \left\{ \hat{f}_0 (X) - f \left( X; \bftheta^* \right) \right\} \left[ 2 \left\{ \Delta (\bftheta) \right\}^{\top} \nabla_{\bftheta} f \left( X; \bftheta^* \right) + \left\{ \Delta (\bftheta) \right\}^{\top} \nabla_{\bftheta}^2 f \left( X; \bftheta^* \right) \Delta (\bftheta) \right] \right] \nonumber \\
        &- \frac{1}{2} \bfz^{\top} \calI_{\bbQ} \left( \bftheta^* \right) \bfz + 2 \calB_2 \cdot K^2 \left( 1 + C_{\textsf{dr}} \right)^2 \left( 1 + C_{\textsf{rf}} \right)^2 \log^{\frac{3}{2}} \left( \frac{d}{\delta} \right) \left( \frac{1}{\sqrt{n_{\bbP}}} + \frac{1}{\sqrt{n_{\bbQ}}} \right) \nonumber \\
        &\cdot \left\{ \sqrt{\frac{\textsf{Trace} \left\{ \calI_{\bbP} \left( \bftheta^* \right) \calI_{\bbQ}^{-2} \left( \bftheta^* \right) \right\}}{n_{\bbP}}} + \sqrt{\frac{\textsf{Trace} \left\{ \calI_{\bbQ}^{-1} \left( \bftheta^* \right) \right\}}{n_{\bbQ}}} \right\}^2 \\
        &+ \calB_{1}^2 \calB_2 \cdot K^2 \left\| \calI_{\bbQ}^{-1} \left( \bftheta^* \right) \right\|_{\textsf{op}}^2 \log^{\frac{5}{2}} \left( \frac{d}{\delta} \right) \left( \frac{1}{\sqrt{n_{\bbP}}} + \frac{1}{\sqrt{n_{\bbQ}}} \right) \left( \frac{1}{n_{\bbP}} + \frac{1}{n_{\bbQ}} \right)^2 \nonumber \\
        &+ \frac{8}{3} \calB_3 \cdot K^3 \left( 1 + C_{\textsf{dr}} \right)^3 \left( 1 + C_{\textsf{rf}} \right)^3 \log^{\frac{3}{2}} \left( \frac{d}{\delta} \right) \left\{ \sqrt{\frac{\textsf{Trace} \left\{ \calI_{\bbP} \left( \bftheta^* \right) \calI_{\bbQ}^{-2} \left( \bftheta^* \right) \right\}}{n_{\bbP}}} + \sqrt{\frac{\textsf{Trace} \left\{ \calI_{\bbQ}^{-1} \left( \bftheta^* \right) \right\}}{n_{\bbQ}}} \right\}^3 \nonumber \\
        &+ \frac{2}{3} \calB_{1}^3 \calB_3 \cdot K^3 \left\| \calI_{\bbQ}^{-1} \left( \bftheta^* \right) \right\|_{\textsf{op}}^3 \log^{3} \left( \frac{d}{\delta} \right) \left( \frac{1}{n_{\bbP}} + \frac{1}{n_{\bbQ}} \right)^3, \nonumber
\end{align}
where the step (i) utilizes the bound \eqref{eqn:subsec:proof_lemma:subsec:proof_thm:upper_bound_alg:dr_covariate_shift_adaptation_v2_v1_v6} together with the following simple inequality:
\[
    \left( x + y \right)^n \leq 2^{n-1} \left( x^n + y^n \right), \quad \forall \left( x, y, n \right) \in \bbR_{+} \times \bbR_{+} \times \bbN.
\]

\indent On the other hand, for every $\bftheta \in \bbB_{r (\delta)} \left( \bftheta^* \right)$, one can find by taking advantage of the lower bound \eqref{eqn:subsec:proof_lemma:subsec:proof_thm:upper_bound_alg:dr_covariate_shift_adaptation_v2_v1_v5} that while being on the event $\calA_{\bbP} \cap \calE_2 (\delta)$,
\begin{align}
    \label{eqn:subsec:proof_lemma:subsec:proof_thm:upper_bound_alg:dr_covariate_shift_adaptation_v2_v1_v10}
        &\widehat{\calR}_{\textsf{DR}} (\bftheta) - \widehat{\calR}_{\textsf{DR}} \left( \bftheta^* \right) \nonumber \\
        \geq \ & \bbE_{X \sim \bbP_{X}} \left[ \left\{ \hat{\rho} (X) - \rho^* (X) \right\} \left\{ \hat{f}_0 (X) - f \left( X; \bftheta^* \right) \right\} \left[ 2 \left\{ \Delta (\bftheta) \right\}^{\top} \nabla_{\bftheta} f \left( X; \bftheta^* \right) + \left\{ \Delta (\bftheta) \right\}^{\top} \nabla_{\bftheta}^2 f \left( X; \bftheta^* \right) \Delta (\bftheta) \right] \right] \nonumber \\
        &+ \frac{1}{2} \left\{ \Delta (\bftheta) - \bfz \right\}^{\top} \calI_{\bbQ} \left( \bftheta^* \right) \left\{ \Delta (\bftheta) - \bfz \right\} - \frac{1}{2} \bfz^{\top} \calI_{\bbQ} \left( \bftheta^* \right) \bfz \nonumber \\
        &- \frac{9}{2} \calB_2 \cdot K^2 \left( 1 + C_{\textsf{dr}} \right)^2 \left( 1 + C_{\textsf{rf}} \right)^2 \log^{\frac{3}{2}} \left( \frac{d}{\delta} \right) \left( \frac{1}{\sqrt{n_{\bbP}}} + \frac{1}{\sqrt{n_{\bbQ}}} \right) \\
        &\cdot \left[ \sqrt{\frac{\textsf{Trace} \left\{ \calI_{\bbP} \left( \bftheta^* \right) \calI_{\bbQ}^{-2} \left( \bftheta^* \right) \right\}}{n_{\bbP}}} + \sqrt{\frac{\textsf{Trace} \left\{ \calI_{\bbQ}^{-1} \left( \bftheta^* \right) \right\}}{n_{\bbQ}}} \right]^2 \nonumber \\
        &- \frac{9}{2} \calB_3 \cdot K^3 \left( 1 + C_{\textsf{dr}} \right)^3 \left( 1 + C_{\textsf{rf}} \right)^3 \log^{\frac{3}{2}} \left( \frac{d}{\delta} \right) \left[ \sqrt{\frac{\textsf{Trace} \left\{ \calI_{\bbP} \left( \bftheta^* \right) \calI_{\bbQ}^{-2} \left( \bftheta^* \right) \right\}}{n_{\bbP}}} + \sqrt{\frac{\textsf{Trace} \left\{ \calI_{\bbQ}^{-1} \left( \bftheta^* \right) \right\}}{n_{\bbQ}}} \right]^3. \nonumber
\end{align}
Subtracting the bound \eqref{eqn:subsec:proof_lemma:subsec:proof_thm:upper_bound_alg:dr_covariate_shift_adaptation_v2_v1_v9} from \eqref{eqn:subsec:proof_lemma:subsec:proof_thm:upper_bound_alg:dr_covariate_shift_adaptation_v2_v1_v11} yields that on the event $\calA_{\bbP} \cap \calE_1 \left( \delta; \calI_{\bbQ}^{-1} \left( \bftheta^* \right) \right) \cap \calE_2 (\delta)$,
\begin{equation}
    \label{eqn:subsec:proof_lemma:subsec:proof_thm:upper_bound_alg:dr_covariate_shift_adaptation_v2_v1_v11}
    \begin{split}
        &\widehat{\calR}_{\textsf{DR}} (\bftheta) - \widehat{\calR}_{\textsf{DR}} \left( \bftheta^* + \bfz \right) \\
        \geq \ & \frac{1}{2} \left\{ \Delta (\bftheta) - \bfz \right\}^{\top} \calI_{\bbQ} \left( \bftheta^* \right) \left\{ \Delta (\bftheta) - \bfz \right\} \\
        &- \left[ \frac{13}{2} \calB_2 \cdot K^2 \left( 1 + C_{\textsf{dr}} \right)^2 \left( 1 + C_{\textsf{rf}} \right)^2 \log^{\frac{3}{2}} \left( \frac{d}{\delta} \right) \left( \frac{1}{\sqrt{n_{\bbP}}} + \frac{1}{\sqrt{n_{\bbQ}}} \right) \right. \\
        &\cdot \left[ \sqrt{\frac{\textsf{Trace} \left\{ \calI_{\bbP} \left( \bftheta^* \right) \calI_{\bbQ}^{-2} \left( \bftheta^* \right) \right\}}{n_{\bbP}}} + \sqrt{\frac{\textsf{Trace} \left\{ \calI_{\bbQ}^{-1} \left( \bftheta^* \right) \right\}}{n_{\bbQ}}} \right]^2 \\
        &+ \calB_{1}^2 \calB_2 \cdot K^2 \left\| \calI_{\bbQ}^{-1} \left( \bftheta^* \right) \right\|_{\textsf{op}}^2 \log^{\frac{5}{2}} \left( \frac{d}{\delta} \right) \left( \frac{1}{\sqrt{n_{\bbP}}} + \frac{1}{\sqrt{n_{\bbQ}}} \right) \left( \frac{1}{n_{\bbP}} + \frac{1}{n_{\bbQ}} \right)^2 \\
        &+ \frac{43}{6} \calB_3 \cdot K^3 \left( 1 + C_{\textsf{dr}} \right)^3 \left( 1 + C_{\textsf{rf}} \right)^3 \log^{\frac{3}{2}} \left( \frac{d}{\delta} \right) \\
        &\cdot \left[ \sqrt{\frac{\textsf{Trace} \left\{ \calI_{\bbP} \left( \bftheta^* \right) \calI_{\bbQ}^{-2} \left( \bftheta^* \right) \right\}}{n_{\bbP}}} + \sqrt{\frac{\textsf{Trace} \left\{ \calI_{\bbQ}^{-1} \left( \bftheta^* \right) \right\}}{n_{\bbQ}}} \right]^3 \\
        &\left. + \frac{2}{3} \calB_{1}^3 \calB_3 \cdot K^3 \left\| \calI_{\bbQ}^{-1} \left( \bftheta^* \right) \right\|_{\textsf{op}}^3 \log^{3} \left( \frac{d}{\delta} \right) \left( \frac{1}{n_{\bbP}} + \frac{1}{n_{\bbQ}} \right)^3 \right]
    \end{split}
\end{equation}
for every $\bftheta \in \bbB_{r (\delta)} \left( \bftheta^* \right)$. At this point, we consider the $d$-dimensional ellipsoid
\begin{equation}
    \label{eqn:subsec:proof_lemma:subsec:proof_thm:upper_bound_alg:dr_covariate_shift_adaptation_v2_v1_v12}
    \begin{split}
        \Gamma (\delta) := \ & \left\{ \bftheta \in \bbR^d : \frac{1}{2} \left\{ \Delta (\bftheta) - \bfz \right\}^{\top} \calI_{\bbQ} \left( \bftheta^* \right) \left\{ \Delta (\bftheta) - \bfz \right\} \right. \\
        \leq \ & \frac{13}{2} \calB_2 \cdot K^2 \left( 1 + C_{\textsf{dr}} \right)^2 \left( 1 + C_{\textsf{rf}} \right)^2 \log^{\frac{3}{2}} \left( \frac{d}{\delta} \right) \left( \frac{1}{\sqrt{n_{\bbP}}} + \frac{1}{\sqrt{n_{\bbQ}}} \right) \\
        &\cdot \left[ \sqrt{\frac{\textsf{Trace} \left\{ \calI_{\bbP} \left( \bftheta^* \right) \calI_{\bbQ}^{-2} \left( \bftheta^* \right) \right\}}{n_{\bbP}}} + \sqrt{\frac{\textsf{Trace} \left\{ \calI_{\bbQ}^{-1} \left( \bftheta^* \right) \right\}}{n_{\bbQ}}} \right]^2 \\
        &+ \calB_{1}^2 \calB_2 \cdot K^2 \left\| \calI_{\bbQ}^{-1} \left( \bftheta^* \right) \right\|_{\textsf{op}}^2 \log^{\frac{5}{2}} \left( \frac{d}{\delta} \right) \left( \frac{1}{\sqrt{n_{\bbP}}} + \frac{1}{\sqrt{n_{\bbQ}}} \right) \left( \frac{1}{n_{\bbP}} + \frac{1}{n_{\bbQ}} \right)^2 \\
        &+ \frac{43}{6} \calB_3 \cdot K^3 \left( 1 + C_{\textsf{dr}} \right)^3 \left( 1 + C_{\textsf{rf}} \right)^3 \log^{\frac{3}{2}} \left( \frac{d}{\delta} \right) \\
        &\cdot \left[ \sqrt{\frac{\textsf{Trace} \left\{ \calI_{\bbP} \left( \bftheta^* \right) \calI_{\bbQ}^{-2} \left( \bftheta^* \right) \right\}}{n_{\bbP}}} + \sqrt{\frac{\textsf{Trace} \left\{ \calI_{\bbQ}^{-1} \left( \bftheta^* \right) \right\}}{n_{\bbQ}}} \right]^3 \\
        &\left. + \frac{2}{3} \calB_{1}^3 \calB_3 \cdot K^3 \left\| \calI_{\bbQ}^{-1} \left( \bftheta^* \right) \right\|_{\textsf{op}}^3 \log^{3} \left( \frac{d}{\delta} \right) \left( \frac{1}{n_{\bbP}} + \frac{1}{n_{\bbQ}} \right)^3 \right\}.
    \end{split}
\end{equation}
Then, using the inequality \eqref{eqn:subsec:proof_lemma:subsec:proof_thm:upper_bound_alg:dr_covariate_shift_adaptation_v2_v1_v11}, it follows for every $\bftheta \in \bbB_{r (\delta)} \left( \bftheta^* \right) \setminus \Gamma (\delta)$ that $\widehat{\calR}_{\textsf{DR}} (\bftheta) - \widehat{\calR}_{\textsf{DR}} \left( \bftheta^* + \bfz \right) > 0$ on the event $\calA_{\bbP} \cap \calE_1 \left( \delta; \calI_{\bbQ}^{-1} \left( \bftheta^* \right) \right) \cap \calE_2 (\delta)$. Together with $\left\| \calI_{\bbQ}^{-1} \left( \bftheta^* \right) \right\|_{\textsf{op}} = \lambda_{\max} \left\{ \calI_{\bbQ}^{-1} \left( \bftheta^* \right) \right\} = \frac{1}{\lambda_{\min} \left\{ \calI_{\bbQ}^{-1} \left( \bftheta^* \right) \right\}}$, we obtain for any $\bftheta \in \Gamma (\delta)$ that
\begin{equation}
    \label{eqn:subsec:proof_lemma:subsec:proof_thm:upper_bound_alg:dr_covariate_shift_adaptation_v2_v1_v13}
    \begin{split}
        &\left\| \Delta (\bftheta) - \bfz \right\|_{2}^2 \\
        \leq \ & 13 \calB_2 \cdot K^2 \left( 1 + C_{\textsf{dr}} \right)^2 \left( 1 + C_{\textsf{rf}} \right)^2 \left\| \calI_{\bbQ}^{-1} \left( \bftheta^* \right) \right\|_{\textsf{op}} \log^{\frac{3}{2}} \left( \frac{d}{\delta} \right) \left( \frac{1}{\sqrt{n_{\bbP}}} + \frac{1}{\sqrt{n_{\bbQ}}} \right) \\
        &\cdot \left[ \sqrt{\frac{\textsf{Trace} \left\{ \calI_{\bbP} \left( \bftheta^* \right) \calI_{\bbQ}^{-2} \left( \bftheta^* \right) \right\}}{n_{\bbP}}} + \sqrt{\frac{\textsf{Trace} \left\{ \calI_{\bbQ}^{-1} \left( \bftheta^* \right) \right\}}{n_{\bbQ}}} \right]^2 \\
        &+ 2 \calB_{1}^2 \calB_2 \cdot K^2 \left\| \calI_{\bbQ}^{-1} \left( \bftheta^* \right) \right\|_{\textsf{op}}^3 \log^{\frac{5}{2}} \left( \frac{d}{\delta} \right) \left( \frac{1}{\sqrt{n_{\bbP}}} + \frac{1}{\sqrt{n_{\bbQ}}} \right) \left( \frac{1}{n_{\bbP}} + \frac{1}{n_{\bbQ}} \right)^2 \\
        &+ \frac{43}{3} \calB_3 \cdot K^3 \left( 1 + C_{\textsf{dr}} \right)^3 \left( 1 + C_{\textsf{rf}} \right)^3 \left\| \calI_{\bbQ}^{-1} \left( \bftheta^* \right) \right\|_{\textsf{op}} \log^{\frac{3}{2}} \left( \frac{d}{\delta} \right) \\
        &\cdot \left[ \sqrt{\frac{\textsf{Trace} \left\{ \calI_{\bbP} \left( \bftheta^* \right) \calI_{\bbQ}^{-2} \left( \bftheta^* \right) \right\}}{n_{\bbP}}} + \sqrt{\frac{\textsf{Trace} \left\{ \calI_{\bbQ}^{-1} \left( \bftheta^* \right) \right\}}{n_{\bbQ}}} \right]^3 \\
        &+ \frac{4}{3} \calB_{1}^3 \calB_3 \cdot K^3 \left\| \calI_{\bbQ}^{-1} \left( \bftheta^* \right) \right\|_{\textsf{op}}^4 \log^{3} \left( \frac{d}{\delta} \right) \left( \frac{1}{n_{\bbP}} + \frac{1}{n_{\bbQ}} \right)^3.
    \end{split}
\end{equation}
Thus, the triangle inequality implies on the event $\calE_1 \left( \delta; \calI_{\bbQ}^{-1} \left( \bftheta^* \right) \right)$ that for any $\bftheta \in \Gamma (\delta)$, one has
\begin{equation}
    \label{eqn:subsec:proof_lemma:subsec:proof_thm:upper_bound_alg:dr_covariate_shift_adaptation_v2_v1_v14}
    \begin{split}
        \left\| \Delta (\bftheta) \right\|_{2}^2 \leq \ & 2 \left\| \Delta (\bftheta) - \bfz \right\|_{2}^2 + 2 \left\| \bfz \right\|_{2}^2 \\
        \stackrel{\textnormal{(j)}}{\leq} \ & 26 \calB_2 \cdot K^2 \left( 1 + C_{\textsf{dr}} \right)^2 \left( 1 + C_{\textsf{rf}} \right)^2 \left\| \calI_{\bbQ}^{-1} \left( \bftheta^* \right) \right\|_{\textsf{op}} \log^{\frac{3}{2}} \left( \frac{d}{\delta} \right) \left( \frac{1}{\sqrt{n_{\bbP}}} + \frac{1}{\sqrt{n_{\bbQ}}} \right) \\
        &\cdot \left[ \sqrt{\frac{\textsf{Trace} \left\{ \calI_{\bbP} \left( \bftheta^* \right) \calI_{\bbQ}^{-2} \left( \bftheta^* \right) \right\}}{n_{\bbP}}} + \sqrt{\frac{\textsf{Trace} \left\{ \calI_{\bbQ}^{-1} \left( \bftheta^* \right) \right\}}{n_{\bbQ}}} \right]^2 \\
        &+ 4 \calB_{1}^2 \calB_2 \cdot K^2 \left\| \calI_{\bbQ}^{-1} \left( \bftheta^* \right) \right\|_{\textsf{op}}^3 \log^{\frac{5}{2}} \left( \frac{d}{\delta} \right) \left( \frac{1}{\sqrt{n_{\bbP}}} + \frac{1}{\sqrt{n_{\bbQ}}} \right) \left( \frac{1}{n_{\bbP}} + \frac{1}{n_{\bbQ}} \right)^2 \\
        &+ \frac{86}{3} \calB_3 \cdot K^3 \left( 1 + C_{\textsf{dr}} \right)^3 \left( 1 + C_{\textsf{rf}} \right)^3 \left\| \calI_{\bbQ}^{-1} \left( \bftheta^* \right) \right\|_{\textsf{op}} \log^{\frac{3}{2}} \left( \frac{d}{\delta} \right) \\
        &\cdot \left[ \sqrt{\frac{\textsf{Trace} \left\{ \calI_{\bbP} \left( \bftheta^* \right) \calI_{\bbQ}^{-2} \left( \bftheta^* \right) \right\}}{n_{\bbP}}} + \sqrt{\frac{\textsf{Trace} \left\{ \calI_{\bbQ}^{-1} \left( \bftheta^* \right) \right\}}{n_{\bbQ}}} \right]^3 \\
        &+ \frac{8}{3} \calB_{1}^3 \calB_3 \cdot K^3 \left\| \calI_{\bbQ}^{-1} \left( \bftheta^* \right) \right\|_{\textsf{op}}^4 \log^{3} \left( \frac{d}{\delta} \right) \left( \frac{1}{n_{\bbP}} + \frac{1}{n_{\bbQ}} \right)^3 \\
        &+ 8 K^2 \left( 1 + C_{\textsf{dr}} \right)^2 \left( 1 + C_{\textsf{rf}} \right)^2 \log \left( \frac{d}{\delta} \right) \\
        &\cdot \left[ \sqrt{\frac{\textsf{Trace} \left\{ \calI_{\bbP} \left( \bftheta^* \right) \calI_{\bbQ}^{-2} \left( \bftheta^* \right) \right\}}{n_{\bbP}}} + \sqrt{\frac{\textsf{Trace} \left\{ \calI_{\bbQ}^{-1} \left( \bftheta^* \right) \right\}}{n_{\bbQ}}} \right]^2 \\
        &+ 4 \calB_{1}^2 \cdot K^2 \left\| \calI_{\bbQ}^{-1} \left( \bftheta^* \right) \right\|_{\textsf{op}}^2 \log^{2} \left( \frac{d}{\delta} \right) \left( \frac{1}{n_{\bbP}} + \frac{1}{n_{\bbQ}} \right)^2,
    \end{split}
\end{equation}
where the step (j) utilize the consequence \eqref{eqn:subsec:proof_lemma:subsec:proof_thm:upper_bound_alg:dr_covariate_shift_adaptation_v2_v1_v6} from Lemma \ref{lemma:concentration_gradient} with $\bfA = \calI_{\bbQ}^{-1} \left( \bftheta^* \right)$. To guarantee that
\[
    \log \left( \frac{d}{\delta} \right) \left[ \sqrt{\frac{\textsf{Trace} \left\{ \calI_{\bbP} \left( \bftheta^* \right) \calI_{\bbQ}^{-2} \left( \bftheta^* \right) \right\}}{n_{\bbP}}} + \sqrt{\frac{\textsf{Trace} \left\{ \calI_{\bbQ}^{-1} \left( \bftheta^* \right) \right\}}{n_{\bbQ}}} \right]^2
\]
is the leading term, we only need it to dominate the remaining terms. In particular, whenver $\min \left\{ n_{\bbP}, n_{\bbQ} \right\} \geq \kappa \cdot \calN_1 \log \left( \frac{d}{\delta} \right)$, then one can conclude from the bound \eqref{eqn:subsec:proof_lemma:subsec:proof_thm:upper_bound_alg:dr_covariate_shift_adaptation_v2_v1_v14} that on the event $\calE_1 \left( \delta; \calI_{\bbQ}^{-1} \left( \bftheta^* \right) \right)$,
\begin{equation}
    \label{eqn:subsec:proof_lemma:subsec:proof_thm:upper_bound_alg:dr_covariate_shift_adaptation_v2_v1_v15}
    \begin{split}
        \left\| \Delta (\bftheta) \right\|_{2}^2 \leq \ & 9 K^2 \left( 1 + C_{\textsf{dr}} \right)^2 \left( 1 + C_{\textsf{rf}} \right)^2 \log \left( \frac{d}{\delta} \right) \left[ \sqrt{\frac{\textsf{Trace} \left\{ \calI_{\bbP} \left( \bftheta^* \right) \calI_{\bbQ}^{-2} \left( \bftheta^* \right) \right\}}{n_{\bbP}}} + \sqrt{\frac{\textsf{Trace} \left\{ \calI_{\bbQ}^{-1} \left( \bftheta^* \right) \right\}}{n_{\bbQ}}} \right]^2 \\
        = \ & r^2 (\delta)
    \end{split}
\end{equation}
for every $\bftheta \in \Gamma (\delta)$. To sum up, we have established the following conclusions so far:
\begin{enumerate} [label = (FACT \Alph*)]
    \item On the event $\calA_{\bbP} \cap \calE_1 \left( \delta; \calI_{\bbQ}^{-1} \left( \bftheta^* \right) \right) \cap \calE_2 (\delta)$, we have $\widehat{\calR}_{\textsf{DR}} (\bftheta) - \widehat{\calR}_{\textsf{DR}} \left( \bftheta^* + \bfz \right) > 0$ for every $\bftheta \in \bbB_{r (\delta)} \left( \bftheta^* \right) \setminus \Gamma (\delta)$.
    \item On the event $\calE_1 \left( \delta; \calI_{\bbQ}^{-1} \left( \bftheta^* \right) \right)$, we have $\Gamma (\delta) \subseteq \bbB_{r (\delta)} \left( \bftheta^* \right)$ if $\min \left\{ n_{\bbP}, n_{\bbQ} \right\} \geq \kappa \cdot \calN_1 \log \left( \frac{d}{\delta} \right)$.
\end{enumerate}

\indent Lastly, it is time to put everything (FACT A and B) together to establish the part (\romannumeral 1) of Lemma \ref{lemma:subsec:proof_thm:upper_bound_alg:dr_covariate_shift_adaptation_v2_v1}. Towards this end, let's claim that being on the event $\calA_{\bbP} \cap \calE_1 \left( \delta; \calI_{\bbQ}^{-1} \left( \bftheta^* \right) \right) \cap \calE_2 (\delta)$, $\widehat{\calR}_{\textsf{DR}} : \bbR^d \to \bbR$ attains a local minimum in the ellipsoid $\Gamma (\delta)$ when $\min \left\{ n_{\bbP}, n_{\bbQ} \right\} \geq \kappa \cdot \calN_1 \log \left( \frac{d}{\delta} \right)$. Due to its continuity together with the compactness of the $d$-dimensional closed ball $\bbB_{r (\delta)} \left( \bftheta^* \right) \subseteq \bbR^d$, $\widehat{\calR}_{\textsf{DR}} : \bbB_{r (\delta)} \left( \bftheta^* \right) \to \bbR$ achieves a global minimum, and so it becomes a local minimum of $\widehat{\calR}_{\textsf{DR}} : \bbR^d \to \bbR$. Let $\overline{\bftheta} \in \argmin \left\{ \widehat{\calR}_{\textsf{DR}} (\bftheta) : \bftheta \in \bbB_{r (\delta)} \left( \bftheta^* \right) \right\}$. Then if $\overline{\bftheta} \in \bbB_{r (\delta)} \left( \bftheta^* \right) \setminus \Gamma (\delta)$, (FACT A) implies
\[
    \widehat{\calR}_{\textsf{DR}} \left( \overline{\bftheta} \right) > \widehat{\calR}_{\textsf{DR}} \left( \bftheta^* + \bfz \right) \stackrel{\textnormal{(k)}}{\geq} \widehat{\calR}_{\textsf{DR}} \left( \overline{\bftheta} \right),
\]
which yields a contradiction, where the step (k) follows since $\bftheta^* + \bfz \in \Gamma (\delta) \subseteq \bbB_{r (\delta)} \left( \bftheta^* \right)$, which holds by (FACT B). Hence, one can conclude that $\overline{\bftheta} \in \Gamma (\delta)$ as desired. Assumption \ref{assumption:simple_landscape_empirical_DR_risk} implies that the global minimum of the empirical DR risk $\widehat{\calR}_{\textsf{DR}} : \bbR^d \to \bbR$ belongs to the ellipsoid $\Gamma (\delta) \subseteq \bbR^d$ on the event $\calA_{\bbP} \cap \calE_1 \left( \delta; \calI_{\bbQ}^{-1} \left( \bftheta^* \right) \right) \cap \calE_2 (\delta)$, i.e., on the event $\calA_{\bbP} \cap \calE_1 \left( \delta; \calI_{\bbQ}^{-1} \left( \bftheta^* \right) \right) \cap \calE_2 (\delta)$,
\begin{equation}
    \label{eqn:subsec:proof_lemma:subsec:proof_thm:upper_bound_alg:dr_covariate_shift_adaptation_v2_v1_v16}
    \begin{split}
        \hat{\bftheta}_{\textsf{DR}} \in \Gamma (\delta) \subseteq \bbB_{r (\delta)} \left( \bftheta^* \right),
    \end{split}
\end{equation}
provided that $\min \left\{ n_{\bbP}, n_{\bbQ} \right\} \geq \kappa \cdot \calN_1 \log \left( \frac{d}{\delta} \right)$.
\medskip

\indent In the sequel, we shall work with the high-probability event $\calA_{\bbP} \cap \calE_1 \left( \delta; \calI_{\bbQ}^{-1} \left( \bftheta^* \right) \right) \cap \calE_2 (\delta)$ in order to further establish the part (\romannumeral 2) of Lemma \ref{lemma:subsec:proof_thm:upper_bound_alg:dr_covariate_shift_adaptation_v2_v1}. Because $\hat{\bftheta}_{\textsf{DR}} \in \Gamma (\delta)$ while being on the event $\calA_{\bbP} \cap \calE_1 \left( \delta; \calI_{\bbQ}^{-1} \left( \bftheta^* \right) \right) \cap \calE_2 (\delta)$ provided that $\min \left\{ n_{\bbP}, n_{\bbQ} \right\} \geq \kappa \cdot \calN_1 \log \left( \frac{d}{\delta} \right)$, it follows that
\begin{equation}
    \label{eqn:subsec:proof_lemma:subsec:proof_thm:upper_bound_alg:dr_covariate_shift_adaptation_v2_v1_v17}
    \begin{split}
        &\left\| \calI_{\bbQ}^{\frac{1}{2}} \left( \bftheta^* \right) \left\{ \Delta \left(  \hat{\bftheta}_{\textsf{DR}} \right) - \bfz \right\} \right\|_{2}^2 \\
        = \ & \left\{ \Delta \left(  \hat{\bftheta}_{\textsf{DR}} \right) - \bfz \right\}^{\top} \calI_{\bbQ} \left( \bftheta^* \right) \left\{ \Delta \left(  \hat{\bftheta}_{\textsf{DR}} \right) - \bfz \right\} \\
        \leq \ & 13 \calB_2 \cdot K^2 \left( 1 + C_{\textsf{dr}} \right)^2 \left( 1 + C_{\textsf{rf}} \right)^2 \log^{\frac{3}{2}} \left( \frac{d}{\delta} \right) \left( \frac{1}{\sqrt{n_{\bbP}}} + \frac{1}{\sqrt{n_{\bbQ}}} \right) \\
        &\cdot \left[ \sqrt{\frac{\textsf{Trace} \left\{ \calI_{\bbP} \left( \bftheta^* \right) \calI_{\bbQ}^{-2} \left( \bftheta^* \right) \right\}}{n_{\bbP}}} + \sqrt{\frac{\textsf{Trace} \left\{ \calI_{\bbQ}^{-1} \left( \bftheta^* \right) \right\}}{n_{\bbQ}}} \right]^2 \\
        &+ 2 \calB_{1}^2 \calB_2 \cdot K^2 \left\| \calI_{\bbQ}^{-1} \left( \bftheta^* \right) \right\|_{\textsf{op}}^2 \log^{\frac{5}{2}} \left( \frac{d}{\delta} \right) \left( \frac{1}{\sqrt{n_{\bbP}}} + \frac{1}{\sqrt{n_{\bbQ}}} \right) \left( \frac{1}{n_{\bbP}} + \frac{1}{n_{\bbQ}} \right)^2 \\
        &+ \frac{43}{3} \calB_3 \cdot K^3 \left( 1 + C_{\textsf{dr}} \right)^3 \left( 1 + C_{\textsf{rf}} \right)^3 \log^{\frac{3}{2}} \left( \frac{d}{\delta} \right) \\
        &\cdot \left[ \sqrt{\frac{\textsf{Trace} \left\{ \calI_{\bbP} \left( \bftheta^* \right) \calI_{\bbQ}^{-2} \left( \bftheta^* \right) \right\}}{n_{\bbP}}} + \sqrt{\frac{\textsf{Trace} \left\{ \calI_{\bbQ}^{-1} \left( \bftheta^* \right) \right\}}{n_{\bbQ}}} \right]^3 \\
        &+ \frac{4}{3} \calB_{1}^3 \calB_3 \cdot K^3 \left\| \calI_{\bbQ}^{-1} \left( \bftheta^* \right) \right\|_{\textsf{op}}^3 \log^{3} \left( \frac{d}{\delta} \right) \left( \frac{1}{n_{\bbP}} + \frac{1}{n_{\bbQ}} \right)^3
    \end{split}
\end{equation}
On the other hand, one can readily apply Lemma \ref{lemma:concentration_gradient} by letting $\bfA = \calI_{\bbQ}^{- \frac{1}{2}} \left( \bftheta^* \right)$ to obtain
\begin{align}
    \label{eqn:subsec:proof_lemma:subsec:proof_thm:upper_bound_alg:dr_covariate_shift_adaptation_v2_v1_v18}
        &\left\| \calI_{\bbQ}^{\frac{1}{2}} \left( \bftheta^* \right) \bfz \right\|_{2}^2 = \left\| \calI_{\bbQ}^{- \frac{1}{2}} \left( \bftheta^* \right) \bfg \right\|_{2}^2 \nonumber \\
        \leq \ & \calC \left\{ \calI_{\bbQ}^{- \frac{1}{2}} \left( \bftheta^* \right) \right\}^2 \left[ \sqrt{\frac{\calV_{\bbP} \left\{ \calI_{\bbQ}^{- \frac{1}{2}} \left( \bftheta^* \right) \right\} \log \left( \frac{d}{\delta} \right)}{n_{\bbP}}} + \sqrt{\frac{\calV_{\bbQ} \left\{ \calI_{\bbQ}^{- \frac{1}{2}} \left( \bftheta^* \right) \right\} \log \left( \frac{d}{\delta} \right)}{n_{\bbQ}}} \right. \nonumber \\
        &\left. + \calB_1 \left\| \calI_{\bbQ}^{- \frac{1}{2}} \left( \bftheta^* \right) \right\|_{\textsf{op}} \log \left( \frac{d}{\delta} \right) \left( \frac{1}{n_{\bbP}} + \frac{1}{n_{\bbQ}} \right) \right]^2 \nonumber \\
        \stackrel{\textnormal{(l)}}{\leq} \ & K^2 \left[ \sqrt{2} \left( 1 + C_{\textsf{dr}} \right) \left( 1 + C_{\textsf{rf}} \right) \sqrt{\log \left( \frac{d}{\delta} \right)} \left\{ \sqrt{\frac{\textsf{Trace} \left\{ \calI_{\bbP} \left( \bftheta^* \right) \calI_{\bbQ}^{-1} \left( \bftheta^* \right) \right\}}{n_{\bbP}}} + \sqrt{\frac{d}{n_{\bbQ}}} \right\} \right. \\
        &\left. + \ \calB_1 \left\| \calI_{\bbQ}^{- \frac{1}{2}} \left( \bftheta^* \right) \right\|_{\textsf{op}} \log \left( \frac{d}{\delta} \right) \left( \frac{1}{n_{\bbP}} + \frac{1}{n_{\bbQ}} \right) \right]^2 \nonumber \\
        \leq \ & 4 K^2 \left( 1 + C_{\textsf{dr}} \right)^2 \left( 1 + C_{\textsf{rf}} \right)^2 \log \left( \frac{d}{\delta} \right) \left\{ \sqrt{\frac{\textsf{Trace} \left\{ \calI_{\bbP} \left( \bftheta^* \right) \calI_{\bbQ}^{-1} \left( \bftheta^* \right) \right\}}{n_{\bbP}}} + \sqrt{\frac{d}{n_{\bbQ}}} \right\}^2 \nonumber \\
        &+ 2 \calB_{1}^2 \cdot K^2 \left\| \calI_{\bbQ}^{-1} \left( \bftheta^* \right) \right\|_{\textsf{op}} \log^2 \left( \frac{d}{\delta} \right) \left( \frac{1}{n_{\bbP}} + \frac{1}{n_{\bbQ}} \right)^2 \nonumber
\end{align}
on the event $\calE_{1} \left( \delta; \calI_{\bbQ}^{- \frac{1}{2}} \left( \bftheta^* \right) \right)$, where the step (l) makes use of \eqref{eqn:subsec:proof_lemma:subsec:proof_thm:upper_bound_alg:dr_covariate_shift_adaptation_v2_v1_v7} and \eqref{eqn:subsec:proof_lemma:subsec:proof_thm:upper_bound_alg:dr_covariate_shift_adaptation_v2_v1_v8} with $\bfA = \calI_{\bbQ}^{- \frac{1}{2}} \left( \bftheta^* \right)$. Thus, we obtain by taking two pieces \eqref{eqn:subsec:proof_lemma:subsec:proof_thm:upper_bound_alg:dr_covariate_shift_adaptation_v2_v1_v17} and \eqref{eqn:subsec:proof_lemma:subsec:proof_thm:upper_bound_alg:dr_covariate_shift_adaptation_v2_v1_v18} collectively that while being on the event
\[
    \Lambda (\delta) := \left\{ \calA_{\bbP} \cap \calE_1 \left( \delta; \calI_{\bbQ}^{-1} \left( \bftheta^* \right) \right) \cap \calE_2 (\delta) \right\} \cap \calE_{1} \left( \delta; \calI_{\bbQ}^{- \frac{1}{2}} \left( \bftheta^* \right) \right),
\]
one has
\begin{equation}
    \label{eqn:subsec:proof_lemma:subsec:proof_thm:upper_bound_alg:dr_covariate_shift_adaptation_v2_v1_v19}
    \begin{split}
        &\left\| \calI_{\bbQ}^{\frac{1}{2}} \left( \bftheta^* \right) \left( \hat{\bftheta}_{\textsf{DR}} - \bftheta^* \right) \right\|_{2}^2 \\
        = \ & \left\| \calI_{\bbQ}^{\frac{1}{2}} \left( \bftheta^* \right) \Delta \left(  \hat{\bftheta}_{\textsf{DR}} \right) \right\|_{2}^2 \\
        \stackrel{\textnormal{(m)}}{\leq} \ & 2 \left\| \calI_{\bbQ}^{\frac{1}{2}} \left( \bftheta^* \right) \left\{ \Delta \left(  \hat{\bftheta}_{\textsf{DR}} \right) - \bfz \right\} \right\|_{2}^2 + 2 \left\| \calI_{\bbQ}^{\frac{1}{2}} \left( \bftheta^* \right) \bfz \right\|_{2}^2 \\
        \leq \ & 26 \calB_2 \cdot K^2 \left( 1 + C_{\textsf{dr}} \right)^2 \left( 1 + C_{\textsf{rf}} \right)^2 \log^{\frac{3}{2}} \left( \frac{d}{\delta} \right) \left( \frac{1}{\sqrt{n_{\bbP}}} + \frac{1}{\sqrt{n_{\bbQ}}} \right) \\
        &\left[ \sqrt{\frac{\textsf{Trace} \left\{ \calI_{\bbP} \left( \bftheta^* \right) \calI_{\bbQ}^{-2} \left( \bftheta^* \right) \right\}}{n_{\bbP}}} + \sqrt{\frac{\textsf{Trace} \left\{ \calI_{\bbQ}^{-1} \left( \bftheta^* \right) \right\}}{n_{\bbQ}}} \right]^2 \\
        &+ 4 \calB_{1}^2 \calB_2 \cdot K^2 \left\| \calI_{\bbQ}^{-1} \left( \bftheta^* \right) \right\|_{\textsf{op}}^2 \log^{\frac{5}{2}} \left( \frac{d}{\delta} \right) \left( \frac{1}{\sqrt{n_{\bbP}}} + \frac{1}{\sqrt{n_{\bbQ}}} \right) \left( \frac{1}{n_{\bbP}} + \frac{1}{n_{\bbQ}} \right)^2 \\
        &+ \frac{86}{3} \calB_3 \cdot K^3 \left( 1 + C_{\textsf{dr}} \right)^3 \left( 1 + C_{\textsf{rf}} \right)^3 \log^{\frac{3}{2}} \left( \frac{d}{\delta} \right) \\
        &\left[ \sqrt{\frac{\textsf{Trace} \left\{ \calI_{\bbP} \left( \bftheta^* \right) \calI_{\bbQ}^{-2} \left( \bftheta^* \right) \right\}}{n_{\bbP}}} + \sqrt{\frac{\textsf{Trace} \left\{ \calI_{\bbQ}^{-1} \left( \bftheta^* \right) \right\}}{n_{\bbQ}}} \right]^3 \\
        &+ \frac{8}{3} \calB_{1}^3 \calB_3 \cdot K^3 \left\| \calI_{\bbQ}^{-1} \left( \bftheta^* \right) \right\|_{\textsf{op}}^3 \log^{3} \left( \frac{d}{\delta} \right) \left( \frac{1}{n_{\bbP}} + \frac{1}{n_{\bbQ}} \right)^3 \\
        &+ 8 K^2 \left( 1 + C_{\textsf{dr}} \right)^2 \left( 1 + C_{\textsf{rf}} \right)^2 \log \left( \frac{d}{\delta} \right) \left\{ \sqrt{\frac{\textsf{Trace} \left\{ \calI_{\bbP} \left( \bftheta^* \right) \calI_{\bbQ}^{-1} \left( \bftheta^* \right) \right\}}{n_{\bbP}}} + \sqrt{\frac{d}{n_{\bbQ}}} \right\}^2 \\
        &+ 4 \calB_{1}^2 \cdot K^2 \left\| \calI_{\bbQ}^{-1} \left( \bftheta^* \right) \right\|_{\textsf{op}} \log^2 \left( \frac{d}{\delta} \right) \left( \frac{1}{n_{\bbP}} + \frac{1}{n_{\bbQ}} \right)^2
    \end{split}
\end{equation}
if $\min \left\{ n_{\bbP}, n_{\bbQ} \right\} \geq \kappa \cdot \calN_1 \log \left( \frac{d}{\delta} \right)$, where the step (m) follows due to the triangle inequality. To guarantee that 
\[
    \log \left( \frac{d}{\delta} \right) \left\{ \sqrt{\frac{\textsf{Trace} \left\{ \calI_{\bbP} \left( \bftheta^* \right) \calI_{\bbQ}^{-1} \left( \bftheta^* \right) \right\}}{n_{\bbP}}} + \sqrt{\frac{d}{n_{\bbQ}}} \right\}^2
\]
becomes the leading term, it suffices to make it dominate the remaining terms. In particular, provided that $\min \left\{ n_{\bbP}, n_{\bbQ} \right\} \geq \kappa \cdot \max \left\{ \calN_1, \calN_2 \right\} \log \left( \frac{d}{\delta} \right)$, then it follows directly from the inequality \eqref{eqn:subsec:proof_lemma:subsec:proof_thm:upper_bound_alg:dr_covariate_shift_adaptation_v2_v1_v19} that while being on the event $\Lambda (\delta)$, we have
\begin{equation}
    \label{eqn:subsec:proof_lemma:subsec:proof_thm:upper_bound_alg:dr_covariate_shift_adaptation_v2_v1_v20}
    \begin{split}
        &\left\| \calI_{\bbQ}^{\frac{1}{2}} \left( \bftheta^* \right) \left( \hat{\bftheta}_{\textsf{DR}} - \bftheta^* \right) \right\|_{2}^2 \\
        \leq \ & 9 K^2 \left( 1 + C_{\textsf{dr}} \right)^2 \left( 1 + C_{\textsf{rf}} \right)^2 \log \left( \frac{d}{\delta} \right) \left\{ \sqrt{\frac{\textsf{Trace} \left\{ \calI_{\bbP} \left( \bftheta^* \right) \calI_{\bbQ}^{-1} \left( \bftheta^* \right) \right\}}{n_{\bbP}}} + \sqrt{\frac{d}{n_{\bbQ}}} \right\}^2.
    \end{split}
\end{equation}
Hence, on the event $\Lambda (\delta)$, the desired results (\romannumeral 1) and (\romannumeral 2) both hold if $\min \left\{ n_{\bbP}, n_{\bbQ} \right\} \geq \kappa \cdot \max \left\{ \calN_1, \calN_2 \right\} \log \left( \frac{d}{\delta} \right)$. This completes the proof of Lemma \ref{lemma:subsec:proof_thm:upper_bound_alg:dr_covariate_shift_adaptation_v2_v1} since
\[
    \begin{split}
        &\left( \bbP^{\otimes n_{\bbP}} \otimes \bbQ_{X}^{\otimes n_{\bbQ}} \right) \left\{ \Lambda (\delta) \right\} \\
        = \ & \left( \bbP^{\otimes n_{\bbP}} \otimes \bbQ_{X}^{\otimes n_{\bbQ}} \right) \left\{ \calE_1 \left( \delta; \calI_{\bbQ}^{-1} \left( \bftheta^* \right) \right) \cap \calE_2 (\delta) \cap \calE_{1} \left( \delta; \calI_{\bbQ}^{- \frac{1}{2}} \left( \bftheta^* \right) \right) \right\} \\
        \stackrel{\textnormal{(n)}}{\geq} \ & 1 - \left( \bbP^{\otimes n_{\bbP}} \otimes \bbQ_{X}^{\otimes n_{\bbQ}} \right) \left\{ \left( \bbO^{n_{\bbP}} \times \bbX_{n_{\bbQ}} \right) \setminus \calE_1 \left( \delta; \calI_{\bbQ}^{-1} \left( \bftheta^* \right) \right) \right\} \\
        &- \left( \bbP^{\otimes n_{\bbP}} \otimes \bbQ_{X}^{\otimes n_{\bbQ}} \right) \left\{ \left( \bbO^{n_{\bbP}} \times \bbX_{n_{\bbQ}} \right) \setminus \calE_2 (\delta) \right\} \\
        &- \left( \bbP^{\otimes n_{\bbP}} \otimes \bbQ_{X}^{\otimes n_{\bbQ}} \right) \left\{ \left( \bbO^{n_{\bbP}} \times \bbX_{n_{\bbQ}} \right) \setminus \calE_1 \left( \delta; \calI_{\bbQ}^{- \frac{1}{2}} \left( \bftheta^* \right) \right) \right\} \\
        \stackrel{\textnormal{(o)}}{\geq} \ & 1 - 8 \delta,
    \end{split}
\]
where the step (n) arises from the union bound, and the step (o) holds true due to Lemma \ref{lemma:concentration_gradient} and \ref{lemma:concentration_hessian_matrix}.

\subsubsection{Proof of Lemma \ref{lemma:concentration_gradient}}
\label{subsubsec:proof_lemma:concentration_gradient}

\indent For simplicity, let $\bfg := \nabla_{\bftheta} \widehat{\calR}_{\textsf{DR}} \left( \bfO_{1:n_{\bbP}}^{\bbP}, \bfX_{1:n_{\bbQ}}^{\bbQ} \right) \left( \bftheta^* \right) - \bbE \left[ \nabla_{\bftheta} \widehat{\calR}_{\textsf{DR}} \left( \bfO_{1:n_{\bbP}}^{\bbP}, \bfX_{1:n_{\bbQ}}^{\bbQ} \right) \left( \bftheta^* \right) \right]$ under $\left( \bfO_{1:n_{\bbP}}^{\bbP}, \bfX_{1:n_{\bbQ}}^{\bbQ} \right) \sim \bbP^{\otimes n_{\bbP}} \otimes \bbQ_{X}^{n_{\bbQ}}$. Then, we obtain from the equation \eqref{eqn:subsec:proof_thm:upper_bound_alg:dr_covariate_shift_adaptation_v2_v6} that
\begin{equation}
    \label{eqn:subsubsec:proof_lemma:concentration_gradient_v1}
    \begin{split}
        \bfg = \left( \widehat{\bbP} - \bbP \right) \left[ \Phi_{\bbP} (X, Y) \right] + \left( \widehat{\bbQ}_X - \bbQ_X \right) \left[ \Phi_{\bbQ} (X) \right] \quad \left( \bbP^{\otimes n_{\bbP}} \otimes \bbQ_{X}^{n_{\bbQ}} \right)\textnormal{-almost surely,}
    \end{split}
\end{equation}
where $\widehat{\bbP} \in \Delta \left( \bbX \times \bbR \right)$ and $\widehat{\bbQ}_{X} \in \Delta (\bbX)$ denote the empirical distributions for the $n_{\bbP}$ labeled source samples $\bfO_{1:n_{\bbP}}^{\bbP}$ and $n_{\bbQ}$ unlabeled target samples $\bfX_{1:n_{\bbQ}}^{\bbQ}$, i.e., $\widehat{\bbP} := \frac{1}{n_{\bbP}} \sum_{i=1}^{n_{\bbP}} \delta_{\left( X_{i}^{\bbP}, Y_{i}^{\bbP} \right)}$ and $\widehat{\bbQ}_{X} := \frac{1}{n_{\bbQ}} \sum_{j=1}^{n_{\bbQ}} \delta_{X_{j}^{\bbQ}}$, respectively, and the functions $\Phi_{\bbP} : \bbX \times \bbR \to \bbR^d$ and $\Phi_{\bbQ} : \bbX \to \bbR^d$ are defined as \eqref{eqn:lemma:concentration_gradient_v2}.
\medskip

\indent Now, we fix any matrix $\bfA \in \bbR^{d \times d}$. One can readily find that
\[
    \begin{split}
        \left\| \bfA \Phi_{\bbP} (X, Y) \right\|_{2} = 2 \hat{\rho} (X) \left| Y - \hat{f}_0 (X) \right| \left\| \bfA \nabla_{\bftheta} f \left( X; \bftheta^* \right) \right\|_{2} \leq 2 C_{\textsf{dr}} \left( 1 + C_{\textsf{rf}} \right) b_1 \left\| \bfA \right\|_{\textsf{op}} \quad \bbP\textnormal{-almost surely,}
    \end{split}
\]
which immediately yields
\begin{align}
    \label{eqn:subsubsec:proof_lemma:concentration_gradient_v2}
        \calB_{\bbP} (\alpha) := \ & \inf \left\{ t \in \left( 0, +\infty \right): \bbE_{(X, Y) \sim \bbP} \left[ \exp \left\{ \left( \frac{\left\| \bfA \left\{ \Phi_{\bbP} (X, Y) - \bbE_{(X, Y) \sim \bbP} \left[ \Phi_{\bbP} (X, Y) \right] \right\} \right\|_{2}}{t} \right)^{\alpha} \right\} \right] \leq 2 \right\} \nonumber \\
        \leq \ & \inf \left\{ t \in \left( 0, +\infty \right): \bbE_{(X, Y) \sim \bbP} \left[ \exp \left\{ \left( \frac{4 C_{\textsf{dr}} \left( 1 + C_{\textsf{rf}} \right) b_1 \left\| \bfA \right\|_{\textsf{op}}}{t} \right)^{\alpha} \right\} \right] \leq 2 \right\} \\
        \leq \ & \calB_1 \left\| \bfA \right\|_{\textsf{op}} \left( \log 2 \right)^{- \frac{1}{\alpha}} \nonumber
\end{align}
for any $\alpha \in \left[ 1, +\infty \right)$. By virtue of Lemma \ref{lemma:generic_Bernstein_inequality}, there exists an absolute constant $\calC_{\bbP} (\bfA) \in \left( 0, +\infty \right)$ such that for any $\delta \in (0, 1]$, it holds by utilizing the observation \eqref{eqn:subsubsec:proof_lemma:concentration_gradient_v2} that
\begin{equation}
    \label{eqn:subsubsec:proof_lemma:concentration_gradient_v3}
    \begin{split}
        \left\| \left( \widehat{\bbP} - \bbP \right) \left[ \bfA \Phi_{\bbP} (X, Y) \right] \right\|_{2} \leq \ & \calC_{\bbP} (\bfA) \left[ \sqrt{\frac{\calV_{\bbP} (\bfA) \log \left( \frac{2d}{\delta} \right)}{n_{\bbP}}} + \calB_1 \left\| \bfA \right\|_{\textsf{op}} \left( \log 2 \right)^{- \frac{1}{\alpha}} \right. \\
        &\left. \cdot \ \log^{\frac{1}{\alpha}} \left\{ \frac{\calB_1 \left\| \bfA \right\|_{\textsf{op}} \left( \log 2 \right)^{- \frac{1}{\alpha}}}{\sqrt{\calV_{\bbP} (\bfA)}} \right\} \cdot \frac{\log \left( \frac{2d}{\delta} \right)}{n_{\bbP}} \right]
    \end{split}
\end{equation}
with probability greater than $1 - \frac{\delta}{2}$ under $\bfO_{1 : n_{\bbP}} \sim \bbP^{\otimes n_{\bbP}}$, for every $\alpha \in \left[ 1, +\infty \right)$. By taking $\alpha \to +\infty$ in the inequality \eqref{eqn:subsubsec:proof_lemma:concentration_gradient_v3}, one can conclude that
\begin{equation}
    \label{eqn:subsubsec:proof_lemma:concentration_gradient_v4}
    \begin{split}
        &\bbP^{\otimes n_{\bbP}} \left( \left\{ \left\| \left( \widehat{\bbP} - \bbP \right) \left[ \bfA \Phi_{\bbP} (X, Y) \right] \right\|_{2} \leq \calC_{\bbP} (\bfA) \right. \right. \\
        &\left. \left. \left\{ \sqrt{\frac{\calV_{\bbP} (\bfA) \log \left( \frac{2d}{\delta} \right)}{n_{\bbP}}} + \calB_1 \left\| \bfA \right\|_{\textsf{op}} \cdot \frac{\log \left( \frac{2d}{\delta} \right)}{n_{\bbP}} \right\} \right\} \right) \geq 1 - \frac{\delta}{2},
    \end{split}
\end{equation}
where $\calV_{\bbP} (\bfA) := \bbE_{(X, Y) \sim \bbP} \left[ \left\| \bfA \left\{ \Phi_{\bbP} (X, Y) - \bbE_{(X, Y) \sim \bbP} \left[ \Phi_{\bbP} (X, Y) \right] \right\} \right\|_{2}^2 \right]$.
\medskip

\indent On the other hand, it can be observed that
\[
    \begin{split}
        \left\| \bfA \Phi_{\bbQ} (X) \right\|_{2} = 2 \left| f \left( X; \bftheta^* \right) - \hat{f}_0 (X) \right| \left\| \bfA \nabla_{\bftheta} f \left( X; \bftheta^* \right) \right\|_{2} \leq 2 \left( 1 + C_{\textsf{rf}} \right) b_1 \left\| \bfA \right\|_{\textsf{op}},
    \end{split}
\]
which directly implies
\begin{align}
    \label{eqn:subsubsec:proof_lemma:concentration_gradient_v5}
        \calB_{\bbQ} (\alpha) := \ & \inf \left\{ t \in \left( 0, +\infty \right): \bbE_{X \sim \bbQ_{X}} \left[ \exp \left\{ \left( \frac{\left\| \bfA \left\{ \Phi_{\bbQ} (X) - \bbE_{X \sim \bbQ_{X}} \left[ \Phi_{\bbQ} (X) \right] \right\} \right\|_{2}}{t} \right)^{\alpha} \right\} \right] \leq 2 \right\} \nonumber \\
        \leq \ & \inf \left\{ t \in \left( 0, +\infty \right): \bbE_{(X, Y) \sim \bbP} \left[ \exp \left\{ \left( \frac{4 \left( 1 + C_{\textnormal{rf}} \right) b_1 \left\| \bfA \right\|_{\textnormal{op}}}{t} \right)^{\alpha} \right\} \right] \leq 2 \right\} \\
        = \ & 4 \left( 1 + C_{\textsf{rf}} \right) b_1 \left\| \bfA \right\|_{\textsf{op}} \cdot \left( \log 2 \right)^{- \frac{1}{\alpha}} \nonumber \\
        \leq \ & \calB_1 \left\| \bfA \right\|_{\textsf{op}} \left( \log 2 \right)^{- \frac{1}{\alpha}} \nonumber
\end{align}
for all $\alpha \in \left[ 1, +\infty \right)$. Applying Lemma \ref{lemma:generic_Bernstein_inequality}, one can see that there exists a universal constant $\calC_{\bbQ} (\bfA) \in \left( 0, +\infty \right)$ such that for any $\delta \in \left( 0, 1 \right]$, we have from the fact \eqref{eqn:subsubsec:proof_lemma:concentration_gradient_v5} that
\begin{equation}
    \label{eqn:subsubsec:proof_lemma:concentration_gradient_v6}
    \begin{split}
        \left\| \left( \widehat{\bbQ}_{X} - \bbQ_{X} \right) \left[ \bfA \Phi_{\bbQ} (X) \right] \right\|_{2} \leq \ & \calC_{\bbQ} (\bfA) \left[ \sqrt{\frac{\calV_{\bbQ} (\bfA) \log \left( \frac{2d}{\delta} \right)}{n_{\bbQ}}} + \calB_1 \left\| \bfA \right\|_{\textsf{op}} \ \left( \log 2 \right)^{- \frac{1}{\alpha}} \right. \\
        &\left. \cdot \ \log^{\frac{1}{\alpha}} \left\{ \frac{\calB_1 \left\| \bfA \right\|_{\textsf{op}} \left( \log 2 \right)^{- \frac{1}{\alpha}}}{\sqrt{\calV_{\bbP} (\bfA)}} \right\} \cdot \frac{\log \left( \frac{2d}{\delta} \right)}{n_{\bbQ}} \right]
    \end{split}
\end{equation}
with probability at least $1 - \frac{\delta}{2}$ under $\bfX_{1 : n_{\bbQ}} \sim \bbQ_{X}^{\otimes n_{\bbQ}}$, for every $\alpha \in \left[ 1, +\infty \right)$. Taking $\alpha \to +\infty$ in the bound \eqref{eqn:subsubsec:proof_lemma:concentration_gradient_v6}, it follows that
\begin{equation}
    \label{eqn:subsubsec:proof_lemma:concentration_gradient_v7}
    \begin{split}
        &\bbQ_{X}^{\otimes n_{\bbQ}} \left( \left\{ \left\| \left( \widehat{\bbQ}_{X} - \bbQ_{X} \right) \left[ \bfA \Phi_{\bbQ} (X) \right] \right\|_{2} \leq \calC_{\bbQ} (\bfA) \right. \right. \\
        &\left. \left. \left\{ \sqrt{\frac{\calV_{\bbQ} (\bfA) \log \left( \frac{2d}{\delta} \right)}{n_{\bbQ}}} + \calB_1 \left\| \bfA \right\|_{\textsf{op}} \cdot \frac{\log \left( \frac{2d}{\delta} \right)}{n_{\bbQ}} \right\} \right\} \right) \geq 1 - \frac{\delta}{2},
    \end{split}
\end{equation}
where $\calV_{\bbQ} (\bfA) := \bbE_{X \sim \bbQ_{X}} \left[ \left\| \bfA \left\{ \Phi_{\bbQ} (X) - \bbE_{X \sim \bbQ_{X}} \left[ \Phi_{\bbQ} (X) \right] \right\} \right\|_{2}^2 \right]$. 
\medskip

\indent Lastly, it is time to put all pieces together. Making use of the union bound together with two conclusions \eqref{eqn:subsubsec:proof_lemma:concentration_gradient_v4} and \eqref{eqn:subsubsec:proof_lemma:concentration_gradient_v7} and setting $\calC (\bfA) := \max \left\{ \calC_{\bbP} (\bfA), \calC_{\bbQ} (\bfA) \right\} \in \left( 0, +\infty \right)$, one has
\[
    \begin{split}
        &\left\| \bfA \left\{ \nabla_{\bftheta} \widehat{\calR}_{\textsf{DR}} \left( \bfO_{1:n_{\bbP}}^{\bbP}, \bfX_{1:n_{\bbQ}}^{\bbQ} \right) \left( \bftheta^* \right) - \bbE_{\left( \bfO_{1:n_{\bbP}}^{\bbP}, \bfX_{1:n_{\bbQ}}^{\bbQ} \right) \sim \bbP^{\otimes n_{\bbP}} \otimes \bbQ_{X}^{n_{\bbQ}}} \left[ \nabla_{\bftheta} \widehat{\calR}_{\textsf{DR}} \left( \bfO_{1:n_{\bbP}}^{\bbP}, \bfX_{1:n_{\bbQ}}^{\bbQ} \right) \left( \bftheta^* \right) \right] \right\} \right\|_{2} \\
        \stackrel{\textnormal{a.s.}}{=} \ & \left\| \left( \widehat{\bbP} - \bbP \right) \left[ \bfA \Phi_{\bbP} (X, Y) \right] + \left( \widehat{\bbQ}_X - \bbQ_X \right) \left[ \bfA \Phi_{\bbQ} (X) \right] \right\|_{2} \\
        \stackrel{\textnormal{(a)}}{\leq} \ & \left\| \left( \widehat{\bbP} - \bbP \right) \left[ \bfA \Phi_{\bbP} (X, Y) \right] \right\|_{2} + \left\| \left( \widehat{\bbQ}_X - \bbQ_X \right) \left[ \bfA \Phi_{\bbQ} (X) \right] \right\|_{2} \\
        \leq \ & \calC (\bfA) \left\{ \sqrt{\frac{\calV_{\bbP} (\bfA) \log \left( \frac{2d}{\delta} \right)}{n_{\bbP}}} + \sqrt{\frac{\calV_{\bbQ} (\bfA) \log \left( \frac{2d}{\delta} \right)}{n_{\bbQ}}} + \calB_1 \left\| \bfA \right\|_{\textsf{op}} \cdot \log \left( \frac{2d}{\delta} \right) \left( \frac{1}{n_{\bbP}} + \frac{1}{n_{\bbQ}} \right) \right\}
    \end{split}
\]
with probability at least $1 - \delta$ under $\left( \bfO_{1:n_{\bbP}}^{\bbP}, \bfX_{1:n_{\bbQ}}^{\bbQ} \right) \sim \bbP^{\otimes n_{\bbP}} \otimes \bbQ_{X}^{n_{\bbQ}}$, where the step (a) follows by the triangle inequality, as desired.

\subsubsection{Proof of Lemma \ref{lemma:concentration_hessian_matrix}}
\label{subsubsec:proof_lemma:concentration_hessian_matrix}

\indent For brevity, we let $\bfH := \nabla_{\bftheta}^2 \widehat{\calR}_{\textsf{DR}} \left( \bfO_{1:n_{\bbP}}^{\bbP}, \bfX_{1:n_{\bbQ}}^{\bbQ} \right) \left( \bftheta^* \right) - \bbE \left[ \nabla_{\bftheta}^2 \widehat{\calR}_{\textsf{DR}} \left( \bfO_{1:n_{\bbP}}^{\bbP}, \bfX_{1:n_{\bbQ}}^{\bbQ} \right) \left( \bftheta^* \right) \right]$ for $\left( \bfO_{1:n_{\bbP}}^{\bbP}, \bfX_{1:n_{\bbQ}}^{\bbQ} \right) \sim \bbP^{\otimes n_{\bbP}} \otimes \bbQ_{X}^{n_{\bbQ}}$. Then, one may express $\bfH$ by using the fact \eqref{eqn:subsec:proof_thm:upper_bound_alg:dr_covariate_shift_adaptation_v2_v6} that
\begin{equation}
    \label{eqn:subsubsec:proof_lemma:concentration_hessian_matrix_v1}
    \begin{split}
        \bfH = \frac{1}{n_{\bbP}} \sum_{i=1}^{n_{\bbP}} \bfU_{i}^{\bbP} + \frac{1}{n_{\bbQ}} \sum_{j=1}^{n_{\bbQ}} \bfV_{j}^{\bbQ},
    \end{split}
\end{equation}
where
\begin{align}
    \label{eqn:subsubsec:proof_lemma:concentration_hessian_matrix_v2}
        \bfU_{i}^{\bbP} := \ & 2 \hat{\rho} \left( X_{i}^{\bbP} \right) \left\{ \hat{f}_0 \left( X_{i}^{\bbP} \right) - Y_{i}^{\bbP} \right\} \nabla_{\bftheta}^{2} f \left( X_{i}^{\bbP}; \bftheta^* \right) - 2 \bbE_{(X, Y) \sim \bbP} \left[ \hat{\rho} (X) \left\{ \hat{f}_0 (X) - Y \right\} \nabla_{\bftheta}^{2} f \left( X; \bftheta^* \right) \right], \nonumber \\
        \bfV_{j}^{\bbQ} := \ & 2 \left[ \nabla_{\bftheta} f \left( X_{j}^{\bbQ}; \bftheta^* \right) \left\{ \nabla_{\bftheta} f \left( X_{j}^{\bbQ}; \bftheta^* \right) \right\}^{\top} + \left\{ f \left( X_{j}^{\bbQ}; \bftheta^* \right) - \hat{f}_0 \left( X_{j}^{\bbQ} \right) \right\} \nabla_{\bftheta}^2 f \left( X_{j}^{\bbQ}; \bftheta^* \right) \right] \\
        &- 2 \bbE_{X \sim \bbQ_{X}} \left[ \nabla_{\bftheta} f \left( X; \bftheta^* \right) \left\{ \nabla_{\bftheta} f \left( X; \bftheta^* \right) \right\}^{\top} + \left\{ f \left( X; \bftheta^* \right) - \hat{f}_0 (X) \right\} \nabla_{\bftheta}^2 f \left( X; \bftheta^* \right) \right], \nonumber
\end{align}
for each $(i, j) \in \left[ n_{\bbP} \right] \times \left[ n_{\bbQ} \right]$. At this moment, one can readily realize the following facts of the $d \times d$ random matrices $\left\{ \bfU_{i}^{\bbP} : i \in \left[ n_{\bbP} \right] \right\}$:
\begin{itemize}
    \item The operator norm of $\bfU_{i}^{\bbP}$ can be bounded as
    \begin{equation}
        \label{eqn:subsubsec:proof_lemma:concentration_hessian_matrix_v3}
        \begin{split}
            \left\| \bfU_{i}^{\bbP} \right\|_{\textsf{op}} \leq \ & 2 \left\| \hat{\rho} \left( X_{i}^{\bbP} \right) \left\{ \hat{f}_0 \left( X_{i}^{\bbP} \right) - Y_{i}^{\bbP} \right\} \nabla_{\bftheta}^{2} f \left( X_{i}^{\bbP}; \bftheta^* \right) \right\|_{\textsf{op}} \\
            &+ 2 \bbE_{(X, Y) \sim \bbP} \left[ \left\| \hat{\rho} (X) \left\{ \hat{f}_0 (X) - Y \right\} \nabla_{\bftheta}^{2} f \left( X; \bftheta^* \right) \right\|_{\textsf{op}} \right] \\
            \stackrel{\textnormal{(a)}}{\leq} \ & 4 C_{\textsf{dr}} \left( 1 + C_{\textsf{rf}} \right) \left\| \nabla_{\bftheta}^{2} f \left( X_{i}^{\bbP}; \bftheta^* \right) \right\|_{\textsf{op}} \\
            \stackrel{\textnormal{(b)}}{\leq} \ & 4 C_{\textsf{dr}} \left( 1 + C_{\textsf{rf}} \right) b_2
        \end{split}
    \end{equation}
    for every $i \in \left[ n_{\bbP} \right]$, where the step (a) follows $\bbP$-almost surely by Assumptions \ref{assumption:uniform_boundedness} and \ref{assumption:black_box_estimates}, and the step (b) comes from Assumption \ref{assumption:smoothness}.
    \item Using the upper bound \eqref{eqn:subsubsec:proof_lemma:concentration_hessian_matrix_v3}, one can obtain $\left( \bfU_{i}^{\bbP} \right)^2 \preceq \sigma_{\bbP}^2 \bfI_d$ for all $i \in \left[ n_{\bbP} \right]$ $\bbP$-almost surely, where
    \[
        \sigma_{\bbP}^2 := 16 C_{\textsf{dr}}^2 \left( 1 + C_{\textsf{rf}} \right)^2 b_{2}^2.
    \]
\end{itemize}

\noindent We can combine the above properties on the $d \times d$ random matrices $\left\{ \bfU_{i}^{\bbP} : i \in \left[ n_{\bbP} \right] \right\}$ together with the matrix Hoeffding inequality (\emph{Theorem 1.3} in \cite{tropp2012user}) to reach
\begin{equation}
    \label{eqn:subsubsec:proof_lemma:concentration_hessian_matrix_v4}
    \begin{split}
        \left( \bbP^{\otimes n_{\bbP}} \right) \left( \left\{ \left\| \frac{1}{n_{\bbP}} \sum_{i=1}^{n_{\bbP}} \bfU_{i}^{\bbP} \right\|_{\textsf{op}} > t \right\} \right) \leq 2 d \exp \left( - \frac{n_{\bbP} t^2}{8 \sigma_{\bbP}^2} \right) = 2d \exp \left\{ - \frac{n_{\bbP} t^2}{128 C_{\textsf{dr}}^2 \left( 1 + C_{\textsf{rf}} \right)^2 b_{2}^2} \right\}
    \end{split}
\end{equation}
for any $t \in \bbR_{+}$. Thus, it follows for any $\delta \in (0, 1]$ that
\begin{equation}
    \label{eqn:subsubsec:proof_lemma:concentration_hessian_matrix_v5}
    \begin{split}
        \left\| \frac{1}{n_{\bbP}} \sum_{i=1}^{n_{\bbP}} \bfU_{i}^{\bbP} \right\|_{\textsf{op}} \leq 8 \sqrt{2} \cdot C_{\textsf{dr}} \left( 1 + C_{\textsf{rf}} \right) b_{2} \sqrt{\frac{\log \left( \frac{2d}{\delta} \right)}{n_{\bbP}}}
    \end{split}
\end{equation}
with probability at least $1 - \delta$ under the probability measure $\bbP^{\otimes n_{\bbP}}$.
\medskip

\indent Likewise, one can make the following observations on the $d \times d$ random matrices $\left\{ \bfV_{j}^{\bbQ} : j \in \left[ n_{\bbQ} \right] \right\}$:
\begin{itemize}
    \item The operator norm of $\bfV_{j}^{\bbQ}$ can be bounded as
    \begin{equation}
        \label{eqn:subsubsec:proof_lemma:concentration_hessian_matrix_v6}
        \begin{split}
            \left\| \bfV_{j}^{\bbQ} \right\|_{\textsf{op}} \leq \ & 2 \left\| \nabla_{\bftheta} f \left( X_{j}^{\bbQ}; \bftheta^* \right) \left\{ \nabla_{\bftheta} f \left( X_{j}^{\bbQ}; \bftheta^* \right) \right\}^{\top} \right\|_{\textsf{op}} \\
            &+ 2 \left\| \left\{ f \left( X_{j}^{\bbQ}; \bftheta^* \right) - \hat{f}_0 \left( X_{j}^{\bbQ} \right) \right\} \nabla_{\bftheta}^2 f \left( X_{j}^{\bbQ}; \bftheta^* \right) \right\|_{\textsf{op}} \\
            &+ 2 \bbE_{X \sim \bbQ_{X}} \left[ \left\| \nabla_{\bftheta} f \left( X; \bftheta^* \right) \left\{ \nabla_{\bftheta} f \left( X; \bftheta^* \right) \right\}^{\top} \right\|_{\textsf{op}} \right] \\
            &+ 2 \bbE_{X \sim \bbQ_{X}} \left[ \left\| \left\{ f \left( X; \bftheta^* \right) - \hat{f}_0 (X) \right\} \nabla_{\bftheta}^2 f \left( X; \bftheta^* \right) \right\|_{\textsf{op}} \right] \\
            = \ & 2 \left\| \nabla_{\bftheta} f \left( X_{j}^{\bbQ}; \bftheta^* \right) \right\|_{2}^2 \\
            &+ 2 \left| f \left( X_{j}^{\bbQ}; \bftheta^* \right) - \hat{f}_0 \left( X_{j}^{\bbQ} \right) \right| \left\| \nabla_{\bftheta}^2 f \left( X_{j}^{\bbQ}; \bftheta^* \right) \right\|_{\textsf{op}} \\
            &+ 2 \bbE_{X \sim \bbQ_{X}} \left[ \left\| \nabla_{\bftheta} f \left( X; \bftheta^* \right) \right\|_{2}^2 \right] \\
            &+ 2 \bbE_{X \sim \bbQ_{X}} \left[ \left| f \left( X; \bftheta^* \right) - \hat{f}_0 (X) \right| \left\| \nabla_{\bftheta}^2 f \left( X; \bftheta^* \right) \right\|_{\textsf{op}} \right] \\
            \stackrel{\textnormal{(c)}}{\leq} \ & 4 \left\{ b_{1}^2 + \left( 1 + C_{\textsf{rf}} \right) b_2 \right\}
        \end{split}
    \end{equation}
    for every $j \in \left[ n_{\bbQ} \right]$, where the step (c) holds due to Assumptions \ref{assumption:uniform_boundedness} and \ref{assumption:black_box_estimates}.
    \item Using the upper bound \eqref{eqn:subsubsec:proof_lemma:concentration_hessian_matrix_v6}, one can obtain $\left( \bfV_{j}^{\bbQ} \right)^2 \preceq \sigma_{\bbQ}^2 \bfI_d$ for every $j \in \left[ n_{\bbQ} \right]$, where
    \[
        \sigma_{\bbQ}^2 := 16 \left\{ b_{1}^2 + \left( 1 + C_{\textsf{rf}} \right) b_2 \right\}^2.
    \]
\end{itemize}

\noindent Making use of the above findings regarding the $d \times d$ random matrices $\left\{ \bfV_{j}^{\bbQ} : j \in \left[ n_{\bbQ} \right] \right\}$, the matrix Hoeffding inequality then reveals that
\begin{equation}
    \label{eqn:subsubsec:proof_lemma:concentration_hessian_matrix_v7}
    \begin{split}
        \left( \bbQ_{X}^{\otimes n_{\bbQ}} \right) \left( \left\{ \left\| \frac{1}{n_{\bbQ}} \sum_{j=1}^{n_{\bbQ}} \bfV_{j}^{\bbQ} \right\|_{\textsf{op}} > t \right\} \right) \leq 2d \exp \left( - \frac{n_{\bbQ} t^2}{8 \sigma_{\bbQ}^2} \right) = 2d \exp \left\{ - \frac{n_{\bbQ} t^2}{128 \left\{ b_{1}^2 + \left( 1 + C_{\textsf{rf}} \right) b_2 \right\}^2} \right\}
    \end{split}
\end{equation}
for every $t \in \bbR_{+}$. The inequality \eqref{eqn:subsubsec:proof_lemma:concentration_hessian_matrix_v7} tells us for any $\delta \in (0, 1]$ that
\begin{equation}
    \label{eqn:subsubsec:proof_lemma:concentration_hessian_matrix_v8}
    \begin{split}
        \left\| \frac{1}{n_{\bbQ}} \sum_{j=1}^{n_{\bbQ}} \bfV_{j}^{\bbQ} \right\|_{\textsf{op}} \leq 8 \sqrt{2} \left\{ b_{1}^2 + \left( 1 + C_{\textsf{rf}} \right) b_2 \right\} \sqrt{\frac{\log \left( \frac{2d}{\delta} \right)}{n_{\bbQ}}}
    \end{split}
\end{equation}
with probability at least $1 - \delta$ under the probability measure $\bbQ_{X}^{\otimes n_{\bbQ}}$. By combining two inequalities \eqref{eqn:subsubsec:proof_lemma:concentration_hessian_matrix_v5} and \eqref{eqn:subsubsec:proof_lemma:concentration_hessian_matrix_v8} together with the union bound and replacing $\delta$ by $\frac{\delta}{2}$ completes the proof of Lemma \ref{lemma:concentration_hessian_matrix}.

\end{document}